\documentclass[11pt,reqno]{amsart} 
\usepackage{amsmath}
\usepackage{amsthm}
\usepackage{amssymb}
\usepackage{amsfonts}
\usepackage{latexsym}

\usepackage{graphicx,color}
\usepackage{enumerate,fancybox}

\renewcommand{\Re}[1]{\operatorname{Re} #1 }

\newcommand{\conj}[1]{\overline{#1}}

\newcommand{\jeqN}{j=1,2,\ldots,N}

\newcommand{\capGm}{{\mathrm{cap}}(\Gamma)}

\newcommand{\grz}{{g_\Omega(z,\infty)}}
\newcommand{\Cbar}{{\overline{\mathbb{C}}}}

\def\Gbar{\overline{G}}

\newcommand{\sta}{\stackrel{\ast}{\longrightarrow}}

\newtheorem{theorem}{Theorem}[section]
\newtheorem{proposition}{Proposition}[section]
\newtheorem{lemma}{Lemma}[section]

\newtheorem{example}{Example}[section]
\newtheorem{corollary}{Corollary}[section]

\theoremstyle{remark}
\newtheorem{remark}{Remark}[section]

\begin{document}
    \title[Bergman polynomials on an archipelago]{Bergman polynomials on
     an Archipelago: Estimates, Zeros and Shape Reconstruction}
\dedicatory{To Christine Chodkiewicz-Putinar, who has enriched and
inspired us by adding a second viola to our quartet}
\date{\today}

\author[Gustafsson]{Bj\"orn Gustafsson} \thanks{{\it Acknowledgements.}
The first author was partially supported by grants from the
Swedish Research Council, the G\"oran Gustafsson Foundation and
the European Network HCAA. The second and third authors were
partially supported by the National Science Foundation, USA, under
grants DMS-0701094, DMS-0603828 and DMS-0808093. The fourth author was
supported by a University of Cyprus research grant. All authors
are indebted to the Mathematical Research Institute at
Oberwolfach, Germany, which provided exceptional working
conditions for our collaborative efforts.}
\address{Department of Mathematics,
          The Royal Institute of Technology,
           S-10044, Stockholm, Sweden}
    \email{gbjorn@math.kth.se}
    \urladdr{http://www.math.kth.se/\textasciitilde gbjorn/}

\author[Putinar]{Mihai Putinar}
    \address{Department of Mathematics,
          University of California at Santa Barbara,
           Santa Barbara, California,
          93106-3080}
    \email{mputinar@math.ucsb.edu}
    \urladdr{http://math.ucsb.edu/\textasciitilde mputinar}

\author[Saff]{Edward\ B.\ Saff}
    \address{Center for Constructive Approximation\\
         Department of Mathematics\\
         Vanderbilt University\\
         1326 Stevenson Center\\
         37240 Nashville\\
         USA}
    \email{Edward.B.Saff@Vanderbilt.Edu}
    \urladdr{http://atlas.math.vanderbilt.edu/~esaff/}

\author[Stylianopoulos]{Nikos Stylianopoulos}
    \address{Department of Mathematics and Statistics,
    University of Cyprus,
    P.O. Box 20537,
            1678 Nicosia,
            Cyprus}
    \email{nikos@ucy.ac.cy
    }
    \urladdr{http://ucy.ac.cy/\textasciitilde nikos}

\begin{abstract}
Growth estimates of complex orthogonal polynomials with respect to
the area measure supported by a disjoint union of planar Jordan
domains (called, in short, an archipelago) are obtained by a
combination of methods of potential theory and rational
approximation theory. The study of the asymptotic behavior of the
roots of these polynomials reveals a surprisingly rich geometry,
which reflects three characteristics: the relative position of an
island in the archipelago, the analytic continuation picture of
the Schwarz function of every individual boundary and the singular
points of the exterior Green function. By way of explicit example,
fine asymptotics are obtained for the lemniscate archipelago
$|z^m-1|<r^m, 0<r<1,$ which consists of $m$ islands. The
asymptotic analysis of the Christoffel functions associated to the
same orthogonal polynomials leads to a very accurate
reconstruction algorithm of the shape of the archipelago, knowing
only finitely many of its power moments. This work naturally
complements a 1969 study by H.\ Widom of Szeg\H{o} orthogonal
polynomials on an archipelago and the more recent asymptotic
analysis of Bergman orthogonal polynomials unveiled by the last
two authors and their collaborators.
\end{abstract}

\maketitle {\bf Archipelago} n. (pl.\ archipelagos or
archipelagoes) an extensive group of islands.


\allowdisplaybreaks

\section{Introduction}\label{section:intro}
The study of orthogonal polynomials, resurrected recently by many groups of scientists, some
departing from the classical framework of constructive approximation to fields as far as quantum
computing or number theory, does not need an introduction. Maybe only our predilection in
the present work for complex analytic orthogonal polynomials on disconnected open sets needs
some justification.

Complex orthogonal polynomials naturally came into focus quite a few decades ago in connection with
problems in rational approximation theory and conformal mapping. The major result, providing strong
asymptotics for Bergman orthogonal polynomials in a domain  with analytic Jordan boundary, goes back
to 1923 to a landmark article by T.\ Carleman \cite{Ca23}. About the same time
S.\ Bernstein discovered that the analogue of Taylor series in non-circular domains (specifically
ellipses in his case) is a Fourier expansion in terms of orthogonal polynomials that are well adapted
to the boundary shape, a phenomenon later elucidated in full generality by J.L.\ Walsh \cite{Wa}.
Then, it came as no surprise that good approximations of conformal mappings of simply-connected
planar domains bear on the Bergman orthogonal polynomials, that is those with respect to the area
measure supported by these domains. By contrast, the theory of orthogonal polynomials on the line
or on the circle has a longer and glorious history, a much wider area of applications and has
attracted an order of magnitude more attention. For history and details the reader can consult
the surveys \cite{Sa90} and \cite{To08} or the monographs \cite{Ga, ST, StTo, Su74}.

Bergman orthogonal polynomials provide a canonical orthonormal basis in the Bergman space of square summable analytic functions associated to a bounded Jordan domain of the complex plane. Contrary to the Hardy-Smirnov space, that is roughly speaking the closure of polynomials in the $L^2$ space with respect to the arc-length measure on a rectifiable Jordan curve, the functions belonging to the Bergman space do not possess non-tangential values on the boundary. This makes their study much more challenging, and less complete as of today. For instance, it is of recent date that the analogues of Blaschke products associated to the Hardy space of the disk have been discovered: the so-called contractive divisors in the Bergman space of the disk, see the monograph by Hedenmalm, Korenblum and Zhu \cite{HKZ}.

It is our aim to discuss in the present work $n$th-root and strong estimates for Bergman orthogonal polynomials on an archipelago, the asymptotics of their zero distribution, and a reconstruction algorithm of the archipelago from a finite set of the associated power moments. The specific choice of the above problems and degree of generality were dictated by the present status of the
theory of complex orthogonal polynomials.

A brief description of the subjects touched in this article follows. Let $G= \cup_{j=1}^N G_j$ be an archipelago, that is a finite union of mutually disjoint bounded Jordan domains of the complex plane. The Bergman orthonormal polynomials with respect to the area measure supported on $G$:
$$
P_n(z) = \lambda_n z^n+ \cdots, \quad \lambda_n>0,\quad
n=0,1,2,\ldots,
$$
carry in a refined (one would be inclined to say, aristocratic) manner the information about $G$.
For instance, simple linear algebra provides a constructive bijection between the sequence
$\{P_n\}_{n=0}^\infty$ and the power moments (correlation matrix entries)
\begin{equation}\label{eq:muG}
\mu_{mn}(G):=\int_G z^n \overline{z}^m dA, \ \ m,n \geq 0,
\end{equation}
where $dA$ stands for the area measure on $\mathbb{C}$.
Three major features distinguish Bergman orthogonal polynomials:
\begin{itemize}
\item[(i)]
An extremality property: ${P_n}/{\lambda_n}$ is the minimum $L^2(G, dA)$-norm monic polynomial of degree $n$,
\item [(ii)]
the Bergman kernel $K(z,\zeta)=\sum_{j=0}^\infty \overline{P_j(\zeta)}P_j(z)$ collects into a condensed form the (derivatives of the) conformal mappings from the disk to every connected component $G_j$,
\item [(iii)]
the square root of the Christoffel function
$\Lambda_n(z):=1/\sqrt{\sum_{j=0}^n |P_j(z)|^2}$
is the extremum value $\min \| q\|_{L^2(G,dA)}, \  q(z)=1, \ {\rm deg}\  q \leq n.$
\end{itemize}

We repeatedly use the above characteristic properties, by combining them with general
methods of potential theory and function theory. An important object in our work is the
multi-valued function
$$
\Phi(z) = \exp\{ g_\Omega(z,\infty) +i g_\Omega^\ast(z)\}, \ \ z \in \mathbb C
\setminus G,
$$
where $g_\Omega(z,\infty)$ is the Green function of the exterior domain
$\Omega:=\overline{\mathbb{C}}\setminus\overline{G}$, with a pole at infinity, and $g_\Omega^\ast$ is any harmonic conjugate of $g_\Omega$. We designate the name \textit{Walsh function} for $\Phi$. At a critical moment in our proofs, we rely on the pioneering work of Widom~\cite{Wi69} that refers to Szeg\H{o}'s orthogonal polynomials on $G$ and their intimate relation to the Walsh function $\Phi$. Our
Bergman space setting, however, departs in quite a few essential points from the Hardy-Smirnov space scenario. Both estimates of the growth of $P_n(z)$ and the limiting distribution of the zero sets of $\{P_n\}_{n=1}^\infty$ depend heavily on $\Phi$ and its analytic continuation across $\partial G$.

While the estimates for $P_n(z)$ are more or less expected, and only how to prove them might bring new turns, the zero distribution picture on an archipelago is full of surprises. The uncovering of this rich geometry began a few years ago, in the work of two of us and collaborators, on the zero distribution of
Bergman orthogonal polynomials on specific Jordan domains, cf.~\cite{LSS, M-DSS, SaSt08}. For example, for the single Jordan~region consisting of the interior of a regular $m$-gon, all the zeros of $P_n$, $n=1,2,\ldots$, lie on the $m$ radial lines joining the center to the vertices, for $m=3$ and $m=4$ (see \cite{MS}), while for $m \geq 5$ every boundary point of the $m$-gon attracts zeros of $P_n$, as $n\to\infty$ (see \cite[Thm.~5]{AB}).

As a byproduct of the estimates we have obtained for $\Lambda_n(z)$, we propose a very accurate
reconstruction-from-moments algorithm. In general, moment data can be regarded as the archetypal, indirect discrete measurements available to an observer, of a complex structure. To give a simple
indication how moments appear in geometric tomography, consider a density function $\rho(x,y)$ with compact support in the complex plane. When performing parallel tomography along a fixed direction $\theta$, one encounters the values of the Radon transform along the fixed screen
$$
R(\rho)(t,\theta)=(\rho(x,y) , \delta(x \cos \theta + y \sin \theta - t))
$$
where $\delta$ stands for Dirac's distribution and $(\cdot,\cdot)$ is the pairing between
test functions and distributions. Computing then the moments with respect to $t$ yields
$$
a_k(\theta) = \int_{\mathbb R} t^k R(\rho)(t, \theta) dt =
 \int_{\mathbb R^2} (x\cos\theta + y\sin\theta)^k \rho(x,y) dx dy.
$$
Denoting the power moments (with respect to the real variables) by
$$
\sigma_{j,k} = \int_{\mathbb R^2} x^j y^k \rho(x,y) dx dy,\ \ i,j \geq 0,
$$
we obtain a linear system
$$
a_k(\theta) = \sum_{i=0}^k \left( \begin{array}{c} k\\
i\end{array} \right) \cos^i\theta \sin^{k-i}\theta \ \sigma_{i,k-i}.
$$
After giving $\theta$ a number of distinct values, and noticing that the
determinant of the system is non-zero, one finds by linear algebra the values
$\{\sigma_{j,k}\}_{j,k=0}^n$.
This procedure was used by the first two authors of this paper in an image
reconstruction algorithm based on a different integral transform of the original measure,
see \cite{GGMPV} and \cite{GHMP}. In a forthcoming work we plan to compare, both computationally and theoretically, these two different reconstruction-from-moments algorithms.

The paper is organized as follows: Sections~\ref{section:basic} and \ref{section:prelim} are devoted to necessary background information. We introduce there the notation, conventions and recall certain facts from potential theory and function theory of a complex variable that needed for the rest of the work. Sections~\ref{Estimates} (Growth Estimates), \ref{section:reconsrt}
(Reconstruction of the Archipelago from Moments) and \ref{section:zeros} (Asymptotic Behavior of Zeros) contain the statements of the main results. In Section~\ref{section:lemniscate} we enter into the only
computational details available among all archipelagoes: disconnected lemniscates with central symmetry. Finally, Section~\ref{section:proofs} contains proofs of previously stated lemmata, propositions and theorems.

\section{Basic concepts}\label{section:basic}
\setcounter{equation}{0}
\subsection{General notations and definitions}
The unit disk, the exterior disk and the extended complex plane are
denoted, respectively,
$$
\mathbb{D}:=\{z\in {\mathbb{C}}:|z|<1\}, \quad \Delta:=\{z\in
{\mathbb{C}}:|z|>1\}\cup \{\infty\}, \quad
\overline{\mathbb{C}}:=\mathbb{C}\cup \{\infty\}.
$$
For the area measure in the complex plane we use $dA = dA(z)=dxdy$,
and for the arc-length measure on a curve we use $|dz|$. By a measure in
general, we always understand a positive Borel measure which is
finite on compact sets. The closed support of a measure $\mu$ is
denoted by ${\rm supp}\,\mu$.

As to curves in the complex plane, we shall use the following
terminology: a \textit{Jordan curve} is a homeomorphic image of
the unit circle into $\mathbb{C}$. (Thus, every Jordan curve in
the present work will be \textbf{bounded}.) An \textit{analytic
Jordan curve} is the image of the unit circle under an analytic
function, defined and univalent in a neighborhood of the circle.
Thus an analytic Jordan curve is by definition smooth. We shall
sometimes need to discuss also Jordan curves which are only
piecewise analytic. The appropriate definitions will then be
introduced as needed.

If $L$ is a Jordan curve, we denote by $\mathrm{int}(L)$ and $\mathrm{ext}(L)$ the
bounded and unbounded, respectively, components of $\overline{\mathbb{C}}\setminus L$.
By a \textit{Jordan domain} we mean the interior of a Jordan curve. If $E\subset\mathbb{C}$
is any set, ${\rm Co\,}(E)$ denotes its convex hull.

The set of polynomials of degree at most $n$ is denoted by $\mathcal{P}_n$.

\subsection{Bergman spaces and polynomials}\label{subsec:bergman}
The main characters in our story are the Bergman orthogonal
polynomials associated to an archipelago $G:=\cup_{j=1}^N G_j$,
where $G_1,...,G_N$ are Jordan domains with mutually disjoint
closures.  Set $\Gamma_j:=\partial G_j$ and $\Gamma:=\cup_{j=1}^N
\Gamma_j$. For later use we introduce also the exterior domain
$\Omega:=\overline{\mathbb{C}}\setminus\overline{G}$. Note that
$\Gamma=\partial G=\partial\Omega$.

Let $\{P_n\}_{n=0}^{\infty}$ denote the sequence of {\em Bergman
orthogonal polynomials} associated with $G$. This is defined as the
sequence of polynomials
$$
P_n(z) = \lambda_n z^n+ \cdots, \quad \lambda_n>0,\quad
n=0,1,2,\ldots,
$$
that are obtained by orthonormalizing the sequence $1,z,z^2,\dots$, with respect to the inner product
$$
\langle f,g\rangle := \int_G f(z) \overline{g(z)} dA.
$$
Equivalently, the corresponding monic polynomials
$P_n(z)/\lambda_n$, can be defined as the unique monic polynomials
of minimal $L^2$-norm over $G$:
\begin{equation}\label{eq:minimal1}
\|\frac{1}{\lambda_n}{P_n}\|_{L^2(G)}= m_n(G,dA):=
\min_{r\in \mathcal{P}_{n-1}}\|z^n+r(z)\|_{L^2(G)},
\end{equation}
where $\|f\|_{{L}^2(G)}:=\langle f,f\rangle^{1/2}$. Thus,
$$
\frac{1}{\lambda_n}=m_n(G,dA).
$$

Let $L_a^2(G)$ denote the Bergman space associated with $G$ and  $\langle \cdot,\cdot\rangle$:
\begin{equation*}
{L}_a^2(G):=\left\{ f\ \mathrm{analytic\ on}\ G\ \mathrm{and}\
\|f\|_{{L}^2(G)}<\infty \right\}.
\end{equation*}
Note that $L_a^2(G)$ is a Hilbert space that possesses a reproducing kernel which we denote by $K(z,\zeta)$. That is, $K(z,\zeta)$ is the unique function $K(z,\zeta):G\times G\to\mathbb{C}$ such that, for all $\zeta\in G$, $K(\cdot,\zeta)\in {L}_a^2(G)$ and
\begin{equation}\label{eq:repro}
f(\zeta)=\langle f,K(\cdot,\zeta)\rangle, \quad \forall\,  f\in {L}_a^2(G).
\end{equation}
Furthermore, due to the reproducing property and the completeness of polynomials in $L_a^2(G)$ (see Lemma~\ref{lem:dense} below), the kernel $K(z,\zeta)$ is given in terms of the Bergman polynomials by
$$
K(z,\zeta)=\sum_{j=0}^\infty \overline{P_j(\zeta)}P_j(z).
$$
We single out the square root of the inverse of the diagonal of the reproducing kernel of $G$
$$
\Lambda(z):=\frac{1}{\sqrt{K(z,z)}}, \quad z\in G,
$$
and the finite sections of $K(z,\zeta)$ and $\Lambda(z)$:
\begin{equation}\label{eq:KLndef}
K_n(z,\zeta):=\sum_{j=0}^n\overline{P_j(\zeta)}P_j(z), \qquad
\Lambda_n(z):=\frac{1}{\sqrt{K_n(z,z)}}.
\end{equation}

We note that the $\Lambda_n(z)$'s  are square roots of the so-called \textit{Christoffel functions} of $G$.

\subsection{Potential theoretic preliminaries}\label{subsec:potential}
Let $Q$ be a polynomial of degree $n$ with zeros
$z_1,z_2,\ldots,z_n$. The \emph{normalized counting measure of the
zeros} of $Q$ is defined by
\begin{equation}
\nu_Q:=\frac{1}{n}\sum_{k=1}^{n}\delta_{z_k},
\end{equation}
where $\delta_z$ denotes the unit point mass at the point $z$. In
other words, for any subset $A$ of $\mathbb{C}$,
$$
\nu_Q(A)=\frac{\mbox{number of zeros of }Q\mbox{ in }A}{n}.
$$

Next, given a sequence $\{\sigma_n\}$ of Borel measures, we say that
$\{\sigma_n\}$ {\it converges in the weak$^*$ sense} to a measure
$\sigma$, symbolically $\sigma_n\sta\sigma$, if
$$
\int f d\sigma_n \longrightarrow \int f d\sigma,\quad n\to\infty,
$$
for every function $f$ continuous on $\Cbar$.

For any finite positive Borel measure $\sigma$ of compact support in
$\mathbb{C}$, we define its \emph{logarithmic potential} by
$$
U^\sigma(z):=\int\log\frac{1}{|z-t|}d\sigma(t).
$$
In particular, if $Q_n$ is a monic polynomial of degree $n$, then
$$
U^{\nu_{Q_n}}(z)=-\frac{1}{n}\log|Q_n(z)|.
$$

Let $\Sigma\subset\mathbb{C}$ be a compact set. Then there is a
smallest number $\gamma\in\mathbb{R}\cup\{+\infty\}$ such that
there exists a probability measure $\mu_\Sigma$ on $\Sigma$ with
$U^{\mu_\Sigma}\leq\gamma$ in $\mathbb{C}$. The
\emph{(logarithmic) capacity} of $\Sigma$ is defined as ${\rm
cap\,}(\Sigma) =e^{-\gamma}$ (interpreted as zero if
$\gamma=+\infty$). If ${\rm cap\,}(\Sigma)>0$, then $\mu_\Sigma$
is unique and is called the \textit{equilibrium measure of
$\Sigma$}. For the definition of capacity of more general sets
than compact sets see, e.g., \cite{Ra} and \cite{ST}. A property
that holds everywhere, except on a set of capacity zero, is said
to hold \textit{quasi everywhere (q.e.)}. For example, it is known
that $U^{\mu_\Sigma} =\gamma$, q.e.\ on $\Sigma$.

Let $W$ denote the unbounded component of
$\overline{\mathbb{C}}\setminus\Sigma$. It is known that
$\mathrm{supp}(\mu_\Sigma)\subset
\partial W$, $\mu_\Sigma=\mu_{\partial W}$ and ${\rm cap\,}(\Sigma)=
{\rm cap\,}(\partial W)$. If ${\rm cap\,}(\Sigma) >0$, then the
equilibrium potential is related to the Green function
$g_W(z,\infty)$  of $W$, with pole at infinity,
by
\begin{equation}\label{UmuinK}
U^{\mu_\Sigma}(z)=\log\frac{1}{\mathrm{cap}(\Sigma)}-g_W
(z,\infty),\quad z\in W.
\end{equation}

In our applications $\partial W$ will be a finite disjoint union of
mutually exterior Jordan curves (typically $\Sigma =\overline{G}$ or
$\Sigma =\Gamma$, $W=\Omega$, $\partial W=\Gamma=\partial\Sigma$, in the notations
of Subsection~\ref{subsec:bergman}). Then, every point of $\partial
W$ is regular for the Dirichlet problem in $W$ \cite[Thm 4.2.2]{Ra} and therefore:
\begin{enumerate}[(i)]
\itemsep=5pt
\item
\begin{equation}\label{eq:musup}
\mathrm{supp}\,\mu_\Sigma=\partial W,
\end{equation}
\item
\begin{equation}\label{UmuonE}
U^{\mu_\Sigma}(z)=\log\frac{1}{\mathrm{cap}(\Sigma)},\quad z\in\Sigma.
\end{equation}
\end{enumerate}

If $\mu$ is a measure on a compact set $\Sigma$ with ${\rm
cap\,}(\Sigma)
>0$, the \emph{balayage} (or ``swept measure") of $\mu$ onto
$\partial\Sigma$ is the unique measure $\nu$ on $\partial\Sigma$
having the same exterior potential as $\mu$, i.e., satisfying
\begin{equation}\label{eq:balayage}
U^\nu=U^{\mu} \quad {\rm in\,\,}\mathbb{C}\setminus\Sigma.
\end{equation}
The potential $U^\nu$ of $\nu$ can be constructed as the smallest
superharmonic function in $\mathbb{C}$ satisfying $U^\nu\geq U^\mu$
in $\mathbb{C}\setminus\Sigma$. Since $U^\mu$ itself is competing it
follows that, in addition to (\ref{eq:balayage}), $U^\nu\leq U^\mu$
in all $\mathbb{C}$.


\subsection{The Green function and its level
curves}\label{subsec:green}

Returning to the archipelago, let $g_\Omega(z,\infty)$ denote the
Green function of
$\Omega=\overline{\mathbb{C}}\setminus\overline{G}$ with pole at
infinity. That is, $g_\Omega(z,\infty)$ is harmonic in
$\Omega\setminus\{\infty\}$, vanishes on the boundary $\Gamma$ of
$G$ and near $\infty$ satisfies
\begin{equation}\label{green_inf}
g_\Omega(z,\infty)=\log|z|+\log\frac{1}{\capGm}+O(\frac{1}{|z|}),
\quad |z|\to\infty.
\end{equation}

We consider next what we call the \textit{Walsh function} associated with $\Omega$. This is defined as the exponential of the complex Green function,
\begin{equation}\label{eq:Phi=}
\Phi(z):=\exp\{g_\Omega(z,\infty)+ig_\Omega^* (z,\infty)\},
\end{equation}
where $g_\Omega^* (z,\infty)$ is a (locally) harmonic conjugate of $\grz$ in $\Omega$. In the single-component case $N=1$, (\ref{eq:Phi=}) defines a conformal mapping from $\Omega$ onto $\Delta$. In the multiple-component case $N\ge 2$, $\Phi$ is a multi-valued analytic function in $\Omega$. However, $|\Phi(z)|$ is single-valued. We refer to Walsh~\cite[\S 4.1]{Wa} and Widom~\cite[\S~4]{Wi69} for comprehensive accounts of the properties of $\Phi$. We note in particular that $\Phi$ is single-valued near infinity and, since $g_\Omega^\ast$ is unique apart from a constant, that it can be chosen so that $\Phi$ has near infinity a Laurent series expansion of the form
\begin{equation}\label{PhiLaurent}
\Phi(z)=\frac{1}{\capGm}\,z+\alpha_0+\frac{\alpha_1}{z}+\frac{\alpha_2}{z^2}+\cdots.
\end{equation}
We also note that $\Phi^\prime(z)/\Phi(z)=2{\partial g_\Omega(\cdot,\infty)}/{\partial z}$
is single-valued and analytic in $\Omega$, with periods
\begin{equation}\label{eq:bj}
b_j:=\frac{1}{2\pi i}\int_{\Gamma_j}\frac{\Phi^\prime(z)}{\Phi(z)}\,dz
=\frac{1}{2\pi }\int_{\Gamma_j}\frac{\partial g_\Omega(z,\infty)}{\partial n}\,ds,\quad \jeqN.
\end{equation}
Here $\Gamma_j$ is oriented as the boundary of $G_j$ and the normal derivative is directed into
$\Omega$. If $\Gamma_j$ is not smooth the path of integration in (\ref{eq:bj}) is understood to
be moved slightly into $\Omega$. Note that
$
\sum_{j=1}^N b_j =1.
$

Next we consider for $R\geq 1$ the level curves (or equipotential
loci) of the Green function,
\begin{equation}\label{eq:LR}
L_R:=\{z\in\Omega:\,\grz=\log R\}=\{z\in\Omega:\,|\Phi(z)|= R\}
\end{equation}
and the open sets
$$
\Omega_R:=\{z\in\Omega:\,g_\Omega(z,\infty)>\log
R\}=\{z\in\Omega:\,|\Phi(z)|> R\}=\mathrm{ext}(L_R),
$$
$$
\mathcal{G}_R:=\overline{\mathbb{C}}\setminus\overline{\Omega}_R
=\mathrm{int}(L_R).
$$
Note that $L_1= \Gamma$, $\Omega_1=\Omega$, $\mathcal{G}_1=G$. It
follows from the maximum principle that $\Omega_R$ is always
connected. The Green function for $\Omega_R$ is given by
\begin{equation}\label{eq:gOmROm}
g_{\Omega_R} (z,\infty)=g_{\Omega} (z,\infty)-\log R,
\end{equation}
hence the capacity of $L_R$ (or $\overline{\mathcal{G}}_R$) is
\begin{equation}\label{eq:CapLRGam}
\mathrm{cap}(L_R)= R\,\mathrm{cap} (\Gamma).
\end{equation}

Unless stated to the contrary, we hereafter assume that $N\ge 2$, i.e. $G$ consists of more than one island. For small values of $R>1$, $\mathcal{G}_R$ consists of $N$ components, each of which contains exactly one component of $G$, while for large values of $R$,
$\mathcal{G}_R$ is connected (with $\overline{G}\subset\mathcal{G}_R$).
Consequently, we introduce the following sets and numbers:
\begin{itemize}
\item[{$\,$}]
$
\mathcal{G}_{j,R}:=\mbox{ the component of }\mathcal{G}_R\mbox{ that contains }G_j,\quad\jeqN.
$
\item[{$\,$}]
$
L_{j,R}:=\partial \mathcal{G}_{j,R},\quad\jeqN.
$
\item[{$\,$}]
$
R_j:=\sup\{R:\, \mathcal{G}_{j,R} \mbox{ contains no other island than }G_j\}.
$
\item[{$\,$}]
$
R^\prime:=\min\{R_1,\dots,R_N\}
=\sup\{R:\, \mathcal{G}_R \mbox{ has }N \mbox{ exactly components }\}.
$
\item[{$\,$}]
$
R^{\prime\prime}:=\inf\{R:\, \mathcal{G}_R \mbox{ is connected }\}
=\inf\{R:\, \Omega_R \mbox{ is simply connected }\}.
$
\end{itemize}
Thus, when $1<R<R'$, $\mathcal{G}_{R}$ is the disjoint union of the domains $\mathcal{G}_{j,R}$, $\jeqN$
and $L_R$ consists of the $N$ mutually exterior analytic Jordan curves $L_{j,R}$, $\jeqN$,
while for $R>R^{\prime\prime}$, we have
$\mathcal{G}_{1,R}=\mathcal{G}_{2,R}=\cdots=\mathcal{G}_{N,R}=\mathcal{G}_{R}$ and $L_R$ is
a single analytic curve.

It is well-known that $g_\Omega(z,\infty)$ has exactly $N-1$ critical points (multiplicities counted),
i.e., points where the gradient $\nabla g_\Omega(z,\infty)$, or equivalently
$\displaystyle{{\Phi^\prime}/{\Phi}=2{\partial g_\Omega(\cdot,\infty)}/{\partial z}}$, vanishes.
These critical points show up as singularities of $L_R$, which are points of self-intersection.
Such singularities must appear when $L_R$ changes topology.
It follows that there are no critical points in $\mathcal{G}_{R'}\setminus\overline{G}$, at least one
critical point on each $L_{R_j}$, $\jeqN$ (one of them is $L_{R'}$), at least one on $L_{R''}$ and no
critical point in $\Omega_{R''}$. Any $L_{j,R}$ that does not contain a critical point is an analytic Jordan curve. In particular, this applies whenever $1<R<R^\prime$ or $R^{\prime\prime}<R<\infty$.

When $R\ge R^{\prime\prime}$, $\Phi$ is the unique conformal map of $\Omega_R$ onto
$\Delta_R:=\{w:|w|>R\}$ that satisfies (\ref{PhiLaurent}) near infinity.

\begin{figure}[h]
\begin{center}
\includegraphics*[width=0.5\linewidth]{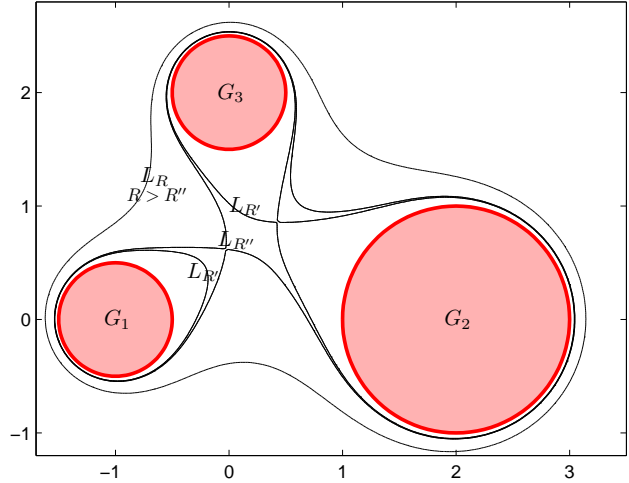}
\end{center}
\caption{Green level curves}\label{fig:greenlines}
\end{figure}

In Figure~\ref{fig:greenlines} we illustrate the three different types of level curves $L_{R'}$,
$L_{R''}$ and $L_{R}$ with $R>R''$, introduced above.

\begin{remark}\label{rem:manydisks}
The level curves in Figure~\ref{fig:greenlines} were computed by means
of Trefethen's MATLAB code {\tt manydisks.m} \cite{Tre}. This code provides an approximation to the
Green function $g_\Omega(z,\infty)$ in cases when $G$ consists of a finite number of disks,
realizing  an algorithm given in \cite{FCB}.
\end{remark}

Consider now the $N$ Hilbert spaces ${L}_a^2(G_j)$ defined by the components $G_j$, $\jeqN$,
and let $K^{G_j}(z,\zeta)$, $\jeqN$, denote their respective reproducing kernels. Then, it is easy to
verify that the kernel $K(z,\zeta)$ is related to $K^{G_j}(z,\zeta)$ as follows:
\begin{equation}\label{eq:KK1}
 K(z,\zeta)=\left\{
\begin{array}{cl}
K^{G_j}(z,\zeta)  &\mathrm{if}\ z,\ \zeta\in G_j,\\
 0            &\mathrm{if}\ z\in G_j,\ \zeta\in G_k,\ j\neq k.
\end{array}
\right.
\end{equation}

In view of (\ref{eq:KK1}), we can express $K(z,\zeta)$ in terms of
conformal mappings $\varphi_j:G_j\to\mathbb{D}$, $\jeqN$. This
will help us to determine the singularities of $K(\cdot,\zeta)$
and, in particular, whether or not this kernel has a singularity
on $\partial G_j$. This is so because, as it is well-known (see
e.g.\ \cite[p.\ 33]{Ga}),
\begin{equation}\label{eq:Kjphij}
K^{G_j}(z,\zeta)=\frac{\varphi'_j(z)\conj{\varphi'_j(\zeta)}}
{\pi\,\left[1-\varphi_j(z)\conj{\varphi_j(\zeta)}\right]^2}, \quad
z,\zeta\in G_j,\,\jeqN.
\end{equation}
By saying that a function analytic in $G_j$ has a singularity on
$\partial G_j$, we mean that there is no open neighborhood of
$\overline{G}_j$ in which the function has an analytic continuation.

\section{Preliminaries}\label{section:prelim}
\setcounter{equation}{0}
\subsection{The Schwarz function of an analytic curve and extension of harmonic functions}
Let $\Gamma$ be a Jordan curve. Then $\Gamma$ is analytic if and only if there exists an analytic function $S(z)$, the \textit{Schwarz function} of $\Gamma$, defined in a full neighborhood of $\Gamma$ and satisfying
$$
S(z)=\bar{z}  \quad \mbox{\rm for } \quad z\in \Gamma;
$$
see \cite{Da74} and \cite{Sh92}. The map $z\mapsto \overline{S(z)}$ is then the anticonformal reflection in $\Gamma$, which is an involution (i.e., is its own inverse) on a suitably defined neighborhood of $\Gamma$. In particular, $S'(z)\ne 0$ in such a neighborhood.

If $u$ is a harmonic function defined at one side of an analytic Jordan curve $\Gamma$ and $u$ has boundary values zero on $\Gamma$, then $u$ extends as a harmonic function across $\Gamma$ by reflection. In terms of the Schwarz function $S(z)$ of $\Gamma$ the extension is given by
\begin{equation}\label{uS}
u(z)=-u(\overline{S(z)})
\end{equation}
for $z$ on the other side of $\Gamma$ (and close to $\Gamma$). Conversely we have the following:

\begin{lemma}\label{lem:schwarz}
Let $\Gamma$ be a Jordan curve and let $u$ be a (real-valued) harmonic function defined in a domain $D$ containing $\Gamma$ such that, for some constant $c>0$, there holds:
\begin{itemize}
\item[(i)]
$
u=0 \quad \mbox{ on }\quad\Gamma,
$
\item[(ii)]
$
|u| \to c \quad\mbox{ as }\quad z\to \partial D,
$
\item[(iii)]
$
\nabla u \ne 0 \quad in\,\, D,
$
\end{itemize}
where $\nabla u$ denotes the gradient of $u$. Then $\Gamma$ is an analytic curve, the Schwarz function $S(z)$ of $\Gamma$ is defined in all $D$, and $z\mapsto \overline{S(z)}$ maps $D$ onto itself.
Moreover, $u$ and $S(z)$ are related by (\ref{uS}). In particular $z\mapsto \overline{S(z)}$ maps a level line $u=\alpha$ of $u$ onto the level line $u=-\alpha$.
\end{lemma}
We note that the Schwarz function is uniquely determined by $\Gamma$, but $u$ is not; there are many different harmonic functions that vanish on $\Gamma$. A domain which is mapped into itself by $z\mapsto \overline{S(z)}$ will be called a {\it domain of involution} for the Schwarz reflection.

\begin{example}\label{ex:Greentextened}
\emph{Assume that, under our main assumptions, one of the components of $\Gamma$, say $\Gamma_1$, is analytic. Then the Green function $u(z)=g_\Omega (z,\infty)$ extends harmonically, by the Schwarz reflection (\ref{uS}), from $\Omega$ into $G_1$. We keep the notation $g_\Omega(z,\infty)$ for so extended Green function. Recall that the level lines reflect to level lines, so that for $R>1$ close enough to one, $L_{1,R}$ is reflected to the level line
$$
L_{1,\frac{1}{R}}=\{z\in G_{1}: g_\Omega (z,\infty)= -\log{R}\}
=\{z\in G_{1}: |\Phi(z)|=\frac{1}{R}\}
$$
of the extended Green function (and extended $\Phi$). Generally speaking, whenever applicable we shall keep notations like $L_{j,\rho}$, $\mathcal{G}_{j,\rho}$, $L_\rho$, $\Omega_\rho$ etc. for values $\rho <1$ in case of analytic boundaries.}
\end{example}

As was previously remarked, $u(z)=g_\Omega (z,\infty)$ has no critical points in the region $\mathcal{G}_{1,R_1}\setminus G_1$. It follows then from (\ref{uS}) that if the Green function extends
harmonically to a region $G_1\setminus\overline{\mathcal{G}}_{1,\rho}$ with $\frac{1}{R_1}\leq \rho<1,$ then it has no critical points there, and the region $D={\mathcal{G}}_{1,\frac{1}{\rho}}
\setminus\overline{\mathcal{G}}_{1,\rho}$ is symmetric with respect to Schwarz reflection and is a region of the kind $D$ discussed in Lemma~\ref{lem:schwarz}.

\subsection{Regular measures}
The class {\bf Reg} of measures of orthogonality was introduced by Stahl and Totik \cite[Definition 3.1.2]{StTobo} and shown to have many desirable properties. Roughly speaking, $\mu\in {\bf Reg}$
means that in an \lq\lq$n$-th root sense", the $\sup$-norm on the support of $\mu$ and the $L^2$-norm generated by $\mu$ have the same asymptotic behavior (as $n\to\infty$) for any sequence of polynomials of respective degrees $n$. It is easy to see that area measure enjoys this property.

\begin{lemma}\label{lem:regular}
The area measure $dA|_G$ on $G$ belongs to the class {\bf Reg}.
\end{lemma}

\medskip
Lemma~\ref{lem:regular} yields the following $n$-th root asymptotic behavior for the Bergman polynomials $P_n$ in $\Omega$.
\begin{proposition}\label{pro:nthroot}
The following assertions hold:
\begin{itemize}
\item[(a)]
\begin{equation}\label{eq:nthlead}
\lim_{n\to\infty}\lambda_n^{1/n}=\frac{1}{\capGm}.
\end{equation}
\item[(b)] For every $z\in\Cbar\setminus{\rm{Co}}(\Gbar)$
and for any $z\in{\rm{Co}}(\Gbar)\setminus\Gbar$ not a limit point of zeros of the $P_n$'s, we have
\begin{equation}\label{eq:nthroot3}
\lim_{n \to\infty}|P_n(z)|^{1/n}=|\Phi(z)|.
\end{equation}
The convergence is uniform on compact subsets of
$\Cbar\setminus{\rm{Co}}(\Gbar)$.
\item[(c)]
\begin{equation}\label{eq:nthroot2}
\limsup_{n \to\infty}|P_n(z)|^{1/n}=|\Phi(z)|,\ z\in\overline{\Omega},
\end{equation}
locally uniformly in ${\Omega}$.
\item[(d)]
\begin{equation}\label{eq:PnprimePn}
\lim_{n \to\infty}\frac{1}{n}\frac{P^\prime_n(z)}{P_n(z)}=
\frac{\Phi^\prime(z)}{\Phi(z)},
\end{equation}
locally uniformly in $\Cbar\setminus{\rm{Co}}(\Gbar)$.
\end{itemize}
\end{proposition}
The first three parts of the proposition  follow from Theorems~3.1.1, 3.2.3 of \cite{StTobo} and from Theorem~III.4.7 of \cite{ST}, in combination with the results of \cite{Am95}, because $\Omega$ is regular with respect to the Dirichlet problem; see e.g.\ \cite[p.~92]{Ra}. The last assertion (d) is immediate from (b).

Another fundamental property of Bergman polynomials, whose proof involves a simple extension
of the simply-connected case treated in Theorem~1, Section~1.3 of Gaier \cite{Ga} is the following.
\begin{lemma}\label{lem:dense}
Polynomials are dense in the Hilbert space ${L}_a^2(G)$.
Consequently, for fixed $\zeta\in G$,
\begin{equation}\label{eq:sumPnPn}
K(z,\zeta)=\sum_{n=0}^\infty\overline{P_n(\zeta)} P_n(z),
\end{equation}
locally uniformly with respect to $z$ in $G$.
\end{lemma}

The analytic continuation properties of $K(z,\zeta)$ play an essential role in the analysis. The following notation will be useful in this regard. If $f$ is an analytic function in $G$, we
define
\begin{equation}\label{rho}
\rho(f):=\sup\left\{R:f\ \mathrm{\ has \ an\ analytic\ continuation\ to}\
\mathcal{G}_R\right\}.
\end{equation}
Note that $1\leq \rho(f)\leq \infty$. The following important lemma, which is an analogue of
the Cauchy-Hadamard formula, is due to Walsh.
\begin{lemma}\label{fouriercoefthm}
Let $f\in {L}^2_a(G)$ . Then,
\begin{equation}\label{fouriercoefequa}
\limsup_{n\to \infty}|\langle f,P_n\rangle|^{1/n}=\frac{1}{\rho(f)}.
\end{equation}
Moreover,
\[
f(z)=\sum_{n=0}^\infty \langle f,P_n\rangle P_n(z),
\]
locally uniformly in $\mathcal{G}_{\rho(f)}$.
\end{lemma}
\noindent
The result is given in Walsh \cite[pp.\ 130--131]{Wa} (see also \cite[Thm~2.1]{PSG}) for a single Jordan region and, as Walsh asserts, is immediately extendable to several Jordan regions.

By applying Lemma~\ref{fouriercoefthm} to $f=K(\cdot,\zeta)$, and by
using the reproducing property (\ref{eq:repro}), in conjunction with
(\ref{eq:KK1}) and (\ref{eq:sumPnPn}), we obtain:
\begin{corollary}\label{cor2}
Let $j$ be fixed, $1\leq j \leq N$. Then for any $\zeta\in G_j$,
\begin{equation}\label{cor2eq}
\limsup_{n\to\infty}|P_n(\zeta)|^{1/n}=\frac{1}{\rho\left(K(\cdot,\zeta)\right)}
=\frac{1}{\min\{\rho\left(K^{G_j}(\cdot,\zeta)\right),R_j\}}\,,
\end{equation}
where (as previously defined) $R_j>1$ is the largest $R$ such that
the component $\mathcal{G}_{j,R}$ of $\mathcal{G}_{R}$ containing
$G_j$ contains no other $G_k$. In particular,
\begin{equation}\label{cor3eq}
\limsup_{n\to \infty}|P_n(\zeta)|^{1/n}=1
\end{equation}
if and only if $K^{G_j}(\cdot,\zeta_0)$ has a singularity on
$\partial G_j$, for some (and then for every) $\zeta_0\in G_j$.
\end{corollary}
\noindent The last statement is based on the observation, from
(\ref{eq:Kjphij}), that the property of $K^{G_j}(\cdot,\zeta_0)$
having a singularity on $\partial G_j$  is independent of the
choice of $\zeta_0$ (within $G_j$). We remark also that the
appearance of $R_j$ in (\ref{cor2eq}) is essential since, for
$R>R_j$, the component $\mathcal{G}_{j,R}$ contains an open set
where $K(\cdot,\zeta)$ is identically zero (recall (\ref{eq:KK1}))
and hence not an analytic continuation of $K^{G_j}(\cdot,\zeta)$.
Corollary~\ref{cor2} will be further elaborated in
Theorem~\ref{thm:ASgeneral}.

Corollary~\ref{cor2} describes a basic relationship between the
orthogonal polynomials $\{P_n(\zeta)\}_{n=0}^\infty$ and the kernel
function $K(\cdot,\zeta)$ which will play an essential role in
deriving our zero distribution results for the sequence
$\{P_n\}_{n=1}^\infty$.

\begin{remark}\label{rem:Fejer}
A well-known result by Fej\'{e}r asserts that the zeros of orthogonal polynomials with respect to a
compactly supported measure $\sigma$ are contained in the closed convex hull of
$\mathrm{supp}\,\sigma$. This result was refined by Saff \cite{Sa90} to the interior of the convex
hull of $\mathrm{supp}\,\sigma$, provided this convex hull is not a line segment. Consequently,
we see that all the zeros of the sequence $\{P_n\}_{n=1}^\infty$ are contained in the interior of
convex hull of $\Gbar$. This fact should be coupled with a result of Widom \cite{Wi67} to the effect
that, on any compact subset $E$ of $\Omega$ and for any $n\in\mathbb{N}$, the number of zeros of
$P_n$ on $E$ is bounded independently of $n$.
\end{remark}

\subsection{Carleman estimates} We continue this section by recalling certain results due to T.~Carleman and P.K.~Suetin, regarding the asymptotic behavior of the Bergman polynomials in the case where $G$ consists of a single component (i.e.\ for $N=1$). In this case the Walsh function (\ref{eq:Phi=}) coincides with the unique conformal map $\Phi:\Omega\to\Delta$ satisfying (\ref{PhiLaurent}).

The first result requires the boundary $\Gamma$ to be analytic (hence the conformal map $\Phi$
has an analytic and univalent continuation across $\Gamma$ inside $G$) and is due to Carleman~\cite{Ca23}; see also \cite[p.~12]{Ga}.
\begin{theorem}\label{th:Carleman}
Assume that $\Gamma$ is an  \textbf{analytic} Jordan curve and let $\rho$, $0<\rho<1$, be the smallest index for which $\Phi$ is  conformal in $\Omega_\rho$. Then,
\begin{equation}\label{eq:Carleman_ln}
\lambda_n=\sqrt{\frac{n+1}{\pi}}\frac{1}{\capGm^{n+1}}\{1+O(\rho^{2n})\},
\end{equation}
and
\begin{equation}\label{eq:Carleman_Pn}
P_n(z)=\sqrt{\frac{n+1}{\pi}}\Phi^\prime(z)\Phi^n(z)\{1+A_n(z)\},
\end{equation}
where
\begin{equation}\label{eq:Carleman_an}
A_n(z)=\left\{
\begin{array}{cl}
O(\sqrt{n})\rho^n, &\mathrm{if}\ z\in\overline{\Omega},\\
O({1}/{\sqrt{n}})\left({\rho}/{r}\right)^n,&\mathrm{if}\ z\in L_r,\,\,\rho<r<1 .
\end{array}
\right.
\end{equation}
\end{theorem}
The second result which is due to Suetin~\cite[Thms 1.1 \& 1.2]{Su74}, requires that $\Gamma$ can be parameterized with respect to the arc-length, so that the defining function has a $p$-th order derivative (where $p$ is a positive integer) in a H\"older class of order $\alpha$. We express this by saying $\Gamma$ is $C^{p+\alpha}$-smooth. (In particular, this implies that $\Gamma$ can have no corners.)

\begin{theorem}\label{th:Suetin}
Assume that $\Gamma$ is $C^{(p+1)+\alpha}$-smooth,
with $p+\alpha>1/2$. Then,
\begin{equation}\label{eq:Suetin_ln}
\lambda_n=\sqrt{\frac{n+1}{\pi}}\frac{1}{\capGm^{n+1}}
\{1+O\left(\frac{1}{n^{2p+2\alpha}}\right)\},
\end{equation}
and
\begin{equation}\label{eq:Suetin_Pn}
P_n(z)=\sqrt{\frac{n+1}{\pi}}\Phi^\prime(z)\Phi^n(z)
\{1+O\left(\frac{\log n}{n^{p+\alpha}}\right)\},\quad
z\in\overline{\Omega}.
\end{equation}
\end{theorem}

We emphasize that the above two theorems concern \textbf{only the case when $N=1$}. We also remark that for the case when $\Gamma$ is analytic, E.\ Mi{\~n}a-D{\'{\i}}az \cite{MD08} has recently derived an improved version of Carleman's theorem for the special case when $L_{\rho}$ is a piecewise analytic Jordan curve without cusps.

\subsection{Comparison of area and line integrals of polynomials}
The following observation is due to Suetin \cite{Su66}; see also \cite[p.\ 38]{Su74}.

\begin{lemma}\label{lem:suetin}
Let $G$ be a Jordan domain with $C^{1+\alpha}$-smooth boundary. Then there
exists a positive constant $C$, with the property that, for every
polynomial $Q_n$ of degree at most $n$, there holds
$$
\|Q_n\|_{L^2(\Gamma)}\leq C\,\sqrt{n+1}\,\|Q_n\|_{L^2(G)},
$$
where $\|\cdot\|_{L^2(\Gamma)}$ denotes the $L^2$-norm on $\Gamma$ with respect to $|dz|$.
\end{lemma}
The proof in \cite{Su66} uses the following analogue of
Bernstein's inequality (a result Suetin attributes to S. Yu.\
Al'per):
$$
\|Q_n^\prime\|_{L^2(G)}\leq C\,n\,\|Q_n\|_{L^2(G)}
$$
and leads to similar inequalities for $L^p$, $1<p<\infty$, or uniform norms.

\section{Growth Estimates}\label{Estimates}
\setcounter{equation}{0} The main results of this article are stated  in this and the next three sections. Their proofs are given in Section~\ref{section:proofs}.

We recall the notation and definitions in Section~\ref{section:basic}, in particular the
definition of the archipelago $G:=\cup_{j=1}^N G_j$ consisting of the union of $N$ Jordan
domains in $\mathbb{C}$, with boundaries $\Gamma_j:=\partial G_j$. We also recall that
$\Omega:=\overline{\mathbb{C}}\setminus\overline{G}$ and note $\Gamma:=\cup_{j=1}^N \Gamma_j=\partial G=\partial\Omega$.
\begin{theorem}\label{th:GPSS_ln}
Assume that every curve $\Gamma_j$ constituting $\Gamma$ is $C^{2+\alpha}$-smooth. Then there exists a positive constant $C_1$ such that
\begin{equation}\label{eq:GPSS_lnle}
\lambda_n\le C_1\sqrt{\frac{n+1}{\pi}}
\frac{1}{\capGm^{n+1}},\quad n\in\mathbb{N}.
\end{equation}
In addition, if every $\Gamma_j$ is analytic, $j=1,2,\ldots,N$, then there exists a
positive constant $C_2$ such that
\begin{equation}\label{eq:GPSS_lnge}
C_2 \sqrt{\frac{n+1}{\pi}} \frac{1}{\capGm^{n+1}}\le \lambda_n \le
C_1 \sqrt{\frac{n+1}{\pi}} \frac{1}{\capGm^{n+1}},\quad
n\in\mathbb{N}.
\end{equation}
\end{theorem}

The following estimate for the diagonal $K(z,z)$, $z\in G$, of the reproducing kernel follows from
classical estimates for the boundary growth of the Bergman kernel of a simply-connected domain,
obtained via conformal mapping techniques. More precisely, by using the results for the hyperbolic
metric presented by Hayman in \cite[pp.\ 682--692]{Ha89}, which require no smoothness for the boundary curve, and recalling (\ref{eq:Kjphij}), it is easy to verify the following double inequality, holding for every $j$, $\jeqN$:
\begin{equation}\label{eq:K_Hayman}
\frac{1}{16\pi}\frac{1}{\textup{dist}^2(z,\Gamma_j)}\ \le\ K^{G_j}(z,z)\le\
\frac{1}{\pi}\frac{1} {\textup{dist}^2(z,\Gamma_j)},\quad z\in G_j,
\end{equation}
Thus $K(z,z)=K^{G_j}(z,z)$, $z\in G_j$, inherits the same estimates and, clearly, the function $\Lambda(z):=1/\sqrt{K(z,z)}$ satisfies
\begin{equation}\label{eq:LamdazGamma}
\sqrt{\pi}\,\textup{dist}(z,\Gamma_j)\, \le\, \Lambda(z)\le\, 4\sqrt{\pi}\,\textup{dist}(z,\Gamma_j),\quad z\in G_j.
\end{equation}
(Above and throughout this article $\textup{dist}(z,\Gamma_j)$ stands for the Euclidean distance of $z$ from the boundary $\Gamma_j$.)

It is always useful to recall that the monic orthogonal polynomials $P_n(z)/\lambda_n$, $n=0,1,\ldots$, satisfy a minimum distance condition with respect to the $L^2$-norm on $G$, in the sense that
\begin{equation}
\frac{1}{\lambda_n}=\|\frac{P_n}{\lambda_n}\|_{L^2(G)}=
\min_{a_0,a_1,...,a_{n-1}}\|z^n+a_{n-1}z^{n-1}+...+a_0\|_{L^2(G)}.
\end{equation}
Similarly, the square root of the Christoffel functions $\Lambda_n(z)$, $n=0,1,\ldots$, defined by by (\ref{eq:KLndef}), can be described as
\begin{equation}
{\Lambda_n}(z)=\min_{p\in\mathcal{P}_n,\atop{p(z)=1}}\|p\|_{L^2(G)};
\end{equation}
cf.\ Lemma~\ref{lem:maxpzet} below. Based on the above two simple extremal properties, we derive the following comparison between $\Lambda_n(z)$ and the functions $\Lambda_n^{G_j}(z):=1/\sqrt{K_n^{G_j}(z,z)}$ associated with each individual island $G_j$.

\begin{theorem}\label{pro:ChrinC}
For every $\jeqN$ and any $n\in\mathbb{N}$,
\begin{equation}\label{eq:ChrinC0}
{\Lambda_n^{G_j}}(z)\le\ {\Lambda_n}(z),\quad z\in \mathbb{C}.
\end{equation}
In addition, if $\Gamma_j$ is analytic, then there exist a sequence $\{\gamma_n\}_{n=0}^\infty$, with $0<\gamma_n<1$ and $\displaystyle{\lim_{n\to\infty}\gamma_n=0}$ geometrically, and a number $m\in\mathbb{N}$, $m\ge 1$, such that for any $ n\in\mathbb{N}$,
\begin{equation}\label{eq:ChrinC}
\frac{1-\gamma_n}{2}\,{\Lambda_{m n}}(z) \le\ {\Lambda_n^{G_j}}(z),\quad z\in \overline{G}_j.
\end{equation}
\end{theorem}

Let $\Phi_j$ denote the normalized, like (\ref{PhiLaurent}), exterior conformal map $\Phi_j:\overline{\mathbb{C}}\setminus\overline{G}_j\to\Delta$. The growth of ${\Lambda_n^{G_j}}(z)$ inside the island $G_j$ is described as follows.
\begin{theorem}\label{th:ChrinGj}
Fix $j$, $\jeqN,$ and assume that $\Gamma_j$ is analytic. Then there exist positive
constants $C_1$, $C_2$ and $\rho<1$  such that for any $ n\in\mathbb{N}$,
\begin{equation}\label{eq:ChrinGj}
0<\ {\Lambda_n^{G_j}}(z)-{\Lambda^{G_j}}(z)\le\
C_1\,|\Phi_j(z)|^n\left( \textup{dist}(z,\Gamma_j)+\frac{1}{n}\right),
\quad z\in\overline{G}_j\setminus\mathcal{G}_{j,\rho}.
\end{equation}
Moreover,
\begin{equation}\label{eq:limnChr}
\lim_{n\to\infty}n {\Lambda_n^{G_j}}(z)=\frac{\sqrt{2\pi}}{|\Phi_j^\prime(z)|},
\end{equation}
uniformly for $z\in\Gamma_j$.

Furthermore, if every curve constituting $\Gamma$ is analytic then
\begin{equation}\label{eq:CrhinGm}
C_3 \leq {n}{\Lambda_n(z)} \leq C_4,\ \ z \in \Gamma,
\end{equation}
and
\begin{equation}\label{eq:ChrinGjb}
 C_5~ {\rm dist} (z, \Gamma)~ \delta(z)\le\ {\Lambda_n}(z)\le \frac{C_6}{\sqrt{n}|\Phi(z)|^n}, \quad
n\in\mathbb{N},\ z \notin \overline{G},
\end{equation}
where $C_3,C_4,C_5, C_6$ are positive constants and
$$ \delta(z) = \frac{|\Phi(z)|^2-1}{|\Phi(z)|} \frac{1}{\sqrt{(n+1)|\Phi(z)|^{2n}(|\Phi(z)|^2-1)+1}}.$$
\end{theorem}
An estimate for $\Lambda_n(z)$ on $\Gamma$ which is finer than (\ref{eq:CrhinGm}), in the sense that it coincides  with (\ref{eq:limnChr}) for the case of a single island, and under weaker smoothness conditions on $\Gamma$, is presented in \cite{Topre}, where asymptotics of  Christoffel functions defined by more general measures on $\mathbb{C}$ are considered.

Finally, we derive the following exterior estimates for Bergman polynomials.
\begin{theorem}\label{th:PnOme}
Assume that every curve constituting $\Gamma$ is analytic. Then the following hold:
\begin{itemize}
\item[(i)]
There exists a positive constant $C$, so that
\begin{equation}
|P_n(z)|\le\ \frac{C}{\textup{dist}(z,\Gamma)}\sqrt{n}|\Phi(z)|^n,\ \ z
\notin \overline{G}.
\end{equation}
\item[(ii)]
For every $\epsilon>0$ there exist a constant $C_\epsilon>0$, such that
$$
|P_n(z)|\ge\ C_\epsilon \sqrt{n}|\Phi(z)|^n,\quad \textup{dist}(z,{\rm
Co}(\overline{G}))\ge\epsilon.
$$
\end{itemize}
\end{theorem}
Note that in the region ${\rm
Co}(\overline{G})\setminus\overline{G}$ the orthogonal polynomials
may have zeros (as the case of the lemniscates considered in Section~\ref{section:lemniscate} illustrates).

\section{Reconstruction of the archipelago from moments}\label{section:reconsrt}
\setcounter{equation}{0}
The present section contains a brief description of a shape reconstruction algorithm. This algorithm is motivated by the estimates established in the previous sections. The comparison of the speed of convergence and accuracy of this approximation scheme with other known ones (see e.g.~\cite{GHMP}) will be analyzed in a separate work.

The algorithm is based on the following observations:
\begin{remark}\label{rem:reconstruction}
$\, $
\begin{itemize}
\item[(i)]
From (\ref{eq:LamdazGamma}) we see that the function $\Lambda(z)$ is bounded from below and above in $G$ by constants times the distance of $z$ to the boundary. Consequently, its truncation
\begin{equation}\label{eq:LambdaPn}
\Lambda_n(z)=\frac{1}{\sqrt{\sum_{k=0}^n|P_k(z)|^2}}
\end{equation}
approximates the distance function to $\Gamma$ in $G$. Furthermore, on $\Gamma$ and in $\Omega$ $\Lambda_n$ decays to zero at certain rates, as $n\to\infty$. More precisely, a close inspection of the inequalities in Theorems~\ref{pro:ChrinC} and \ref{th:ChrinGj}, in conjunction with (\ref{eq:LamdazGamma}), reveals the following asymptotic behavior of $\Lambda_n(z)$ in $\mathbb{C}$:
\begin{itemize}
\item[(a)]
$\displaystyle{\sqrt{\pi}\,\textup{dist}(z,\Gamma)\, \le\, \Lambda_n(z), \quad z\in G}$;
\item[(b)]
$\displaystyle{\Lambda_n(z)\, \le\,C\, \textup{dist}(z,\Gamma),\quad z\in G\cap\Omega_\rho}\,$, for some $0<\rho< 1$ and $C\ge\sqrt{\pi}$;
\item[(c)]
$\displaystyle{\Lambda_n(z)\asymp\frac{1}{n}},\quad z\in\Gamma$;
\item[(d)]
$\displaystyle{\Lambda_n(z)\asymp\frac{1}{\sqrt{n}|\Phi(z)|^n}},\quad z\in\Omega$.
\end{itemize}
\item[(ii)]
In order to construct $\Lambda_n$ we need to have available the finite section $\{P_0,P_1,\ldots,P_n\}$ of Bergman polynomials, and this can be determined by means of the Gram-Schmidt process, requiring only the power moments (\ref{eq:muG}), of degree less or equal than $n$ in each variable.
\item[(iii)]
For any $n=1,2,\ldots$, all the zeros of $P_n(z)$ lie in the interior of the convex hull of $\overline{G}$; see Remark~\ref{rem:Fejer}.
\end{itemize}
The expression $A\asymp B$ means that $C_1 B\le A\le C_2 B$ for positive constants $C_1$ and $C_2$.
\end{remark}

Consequently, Remark~\ref{rem:reconstruction} supports the following algorithm for reconstructing the archipelago $G$, by using a given finite set of the associated power moments
$$
\mu_{ij}:=\langle z^i,z^j\rangle=\int_G
z^i\overline{z}^j\,dA(z),\quad i,j=0,1,\ldots,n.
$$

\noindent {\tt Reconstruction Algorithm}
\begin{enumerate}[1.]
\itemsep=5pt
\item
\textit{Use an Arnoldi version of the Gram-Schmidt process, in the way indicated in \cite{NSpre}, to construct the Bergman polynomials $\{P_k\}_{k=0}^n$
from $\mu_{ij}$, $i,j=0,1,\ldots,n$. This involves at the $k$-step the orthonormalization of the set $\{P_0,P_1,\ldots,P_{k-1},zP_{k-1}\}$, rather than the set of monomials $\{1,z,\ldots,z^{k-1},z^k\}$, as in the standard Gram-Schmidt process.}
\item
\textit{Plot the zeros of $P_n$, $n=1,2,\ldots,n$.}
\item
\textit{Form $\Lambda_n(z)$ as in (\ref{eq:LambdaPn}}).
\item
\textit{Plot the level curves of the function ${\Lambda_n}(x+iy)$ on a
suitable rectangular frame for $(x,y)$ that surrounds the plotted
zero set.}
\end{enumerate}

\medskip
Regarding the stability of the Gram-Schmidt process in the {Reconstruction Algorithm}, we note a fact that pointed out in \cite{NSpre}. That is the  Arnoldi version of the Gram-Schmidt does not suffer from the severe ill-conditioning associated with its ordinary use; see, for instance, the theoretical and numerical evidence reported in \cite{PaWa}. This feature of the Arnoldi Gram-Schmidt has enabled us to compute accurately Bergman polynomials for degrees as high as 160.
We also note that the use of orthogonal polynomials in a reconstruction-from-moments algorithm, was first employed in \cite{NSpre}. However, the algorithm of \cite{NSpre} is only suitable for the single island case $N=1$.

Applications of the Reconstruction Algorithm are illustrated in the following six examples. In each example, the only information used from the associated archipelago $G$ was its power moments. The resulting plots indicate a remarkable fitting, even in the case of non-smooth boundaries, for which our theory, as stated in Section~\ref{Estimates}, does not apply.
\begin{figure}[h]
\begin{center}
\includegraphics*[scale=0.6]{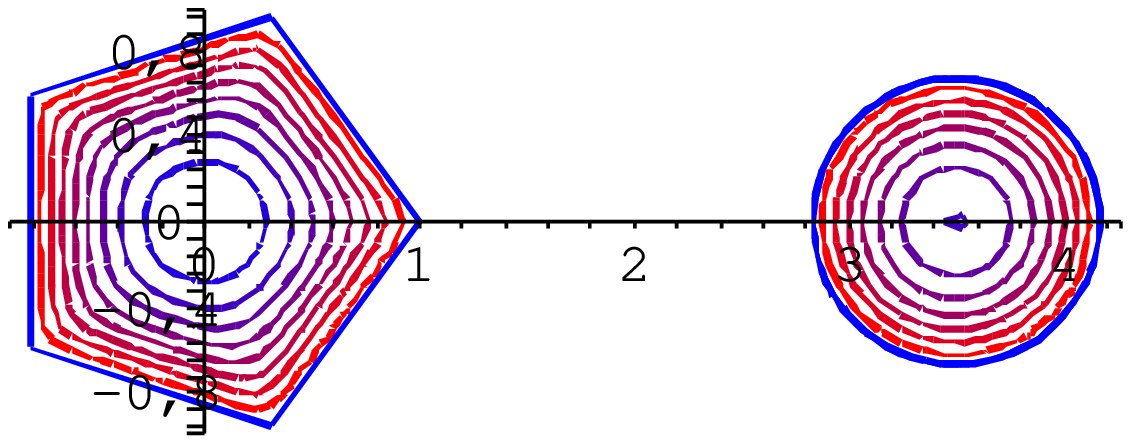}
\end{center}
\caption{Level curves of ${\Lambda_{100}}(x+iy)$, on $\{(x,y):-2\le
x\le 5,-2\le y\le 2\}$, with $G$ as in Example~\ref{ex:diskpen}.}
\label{fig:3.diskpen}
\end{figure}

\begin{figure}[h]
\begin{center}
\includegraphics*[scale=0.6]{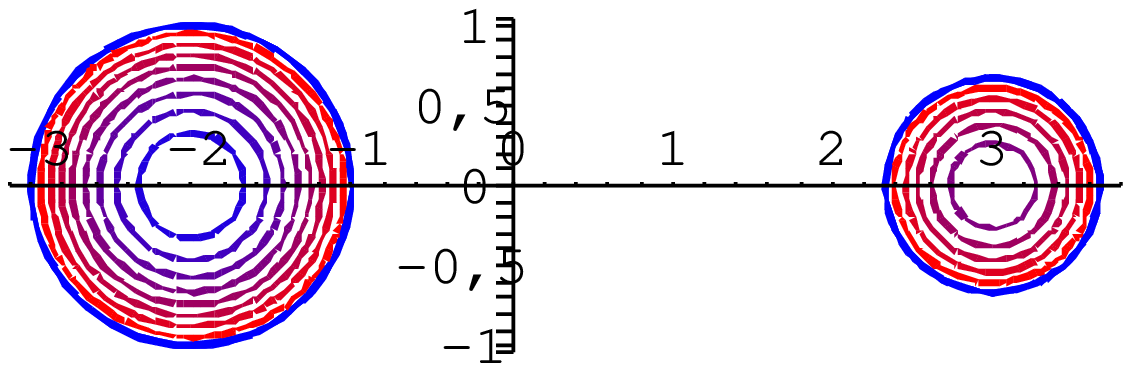}
\end{center}
\caption{Level curves of ${\Lambda_{100}}(x+iy)$, on $\{(x,y):-4\le x\le 4,-2\le y\le 2\}$, with $G$ as in Example~\ref{ex:2disks}.}
\label{fig:3.2disks}
\end{figure}

\begin{figure}[h]
\begin{center}
\includegraphics*[scale=0.6]{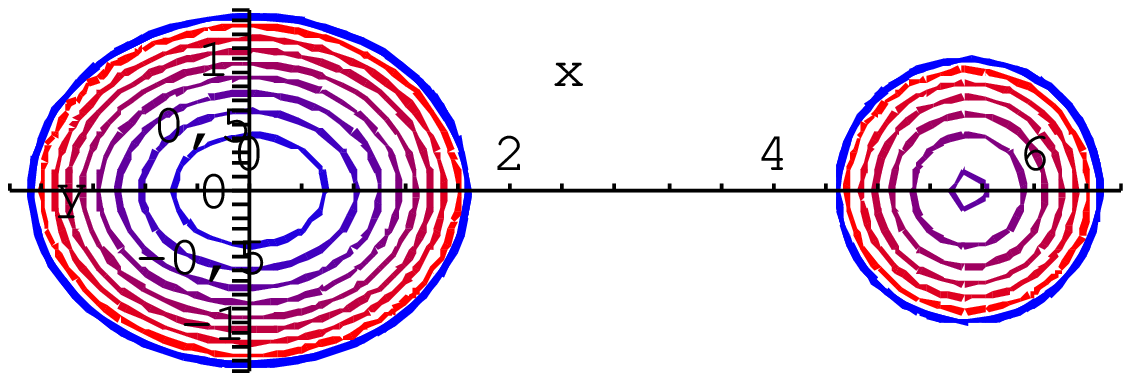}
\end{center}
\caption{Level curves of ${\Lambda_{100}}(x+iy)$, on $\{(x,y):-2\le x\le 8,-2\le y\le 2\}$, with $G$ as in Example~\ref{ex:diskell}, case (i).}
\label{fig:3.diskell}
\end{figure}

\begin{figure}[h]
\begin{center}
\includegraphics*[scale=0.6]{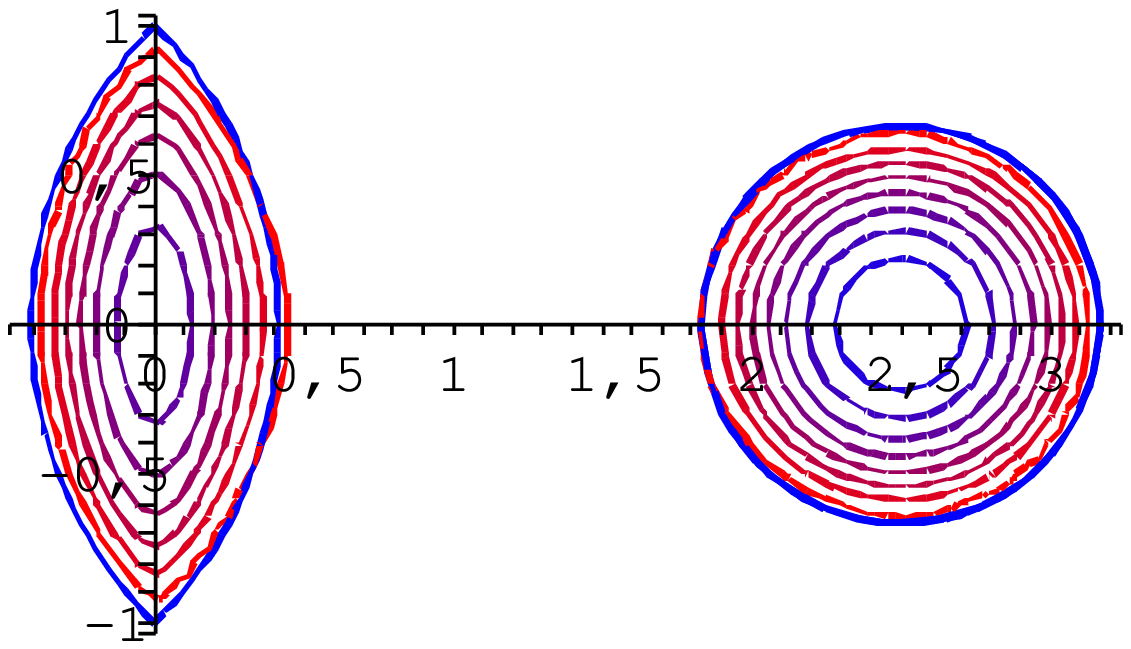}
\end{center}
\caption{Level lines of ${\Lambda_{100}}(x+iy)$, on $\{(x,y):-1\le x\le 4,-2\le y\le 2\}$, with $G$ as in Example~\ref{ex:disklens}.}
\label{fig:3.disklens}
\end{figure}

\begin{figure}[h]
\begin{center}
\includegraphics*[scale=0.5]{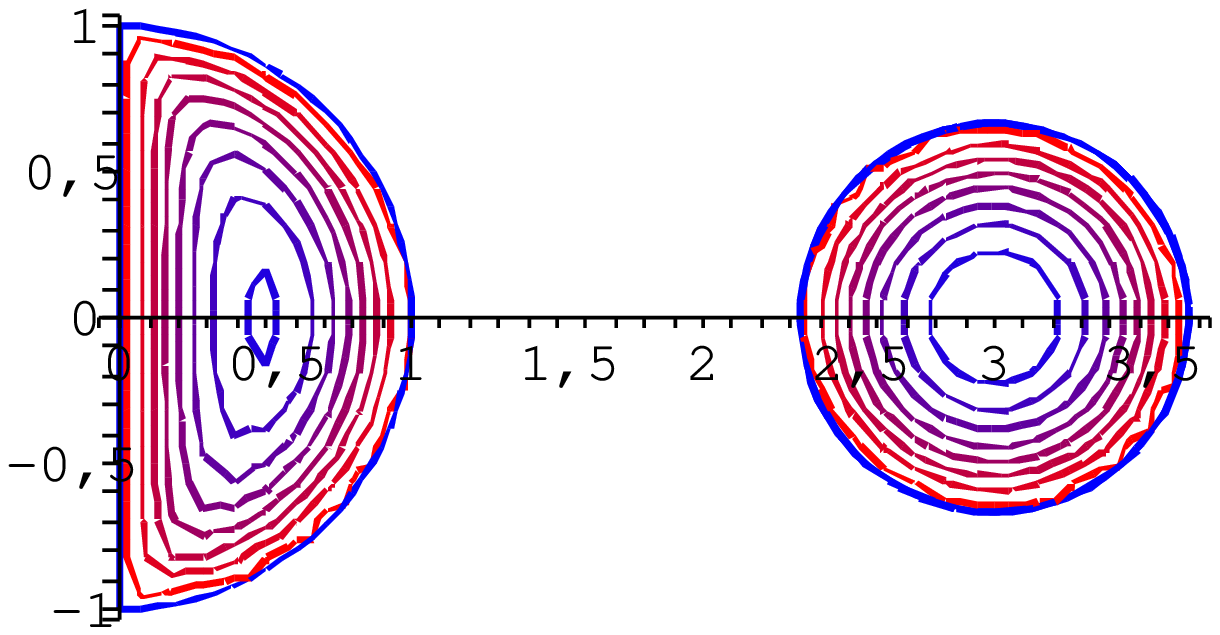}
\end{center}
\caption{Level lines of ${\Lambda_{100}}(x+iy)$, on $\{(x,y):-1\le x\le 6,-2\le y\le 2\}$, with $G$ as in Example~\ref{ex:HDdisk}.}
\label{fig:3.HDdisk}
\end{figure}

\begin{figure}[h]
\begin{center}
\fbox{\includegraphics*[scale=0.5]{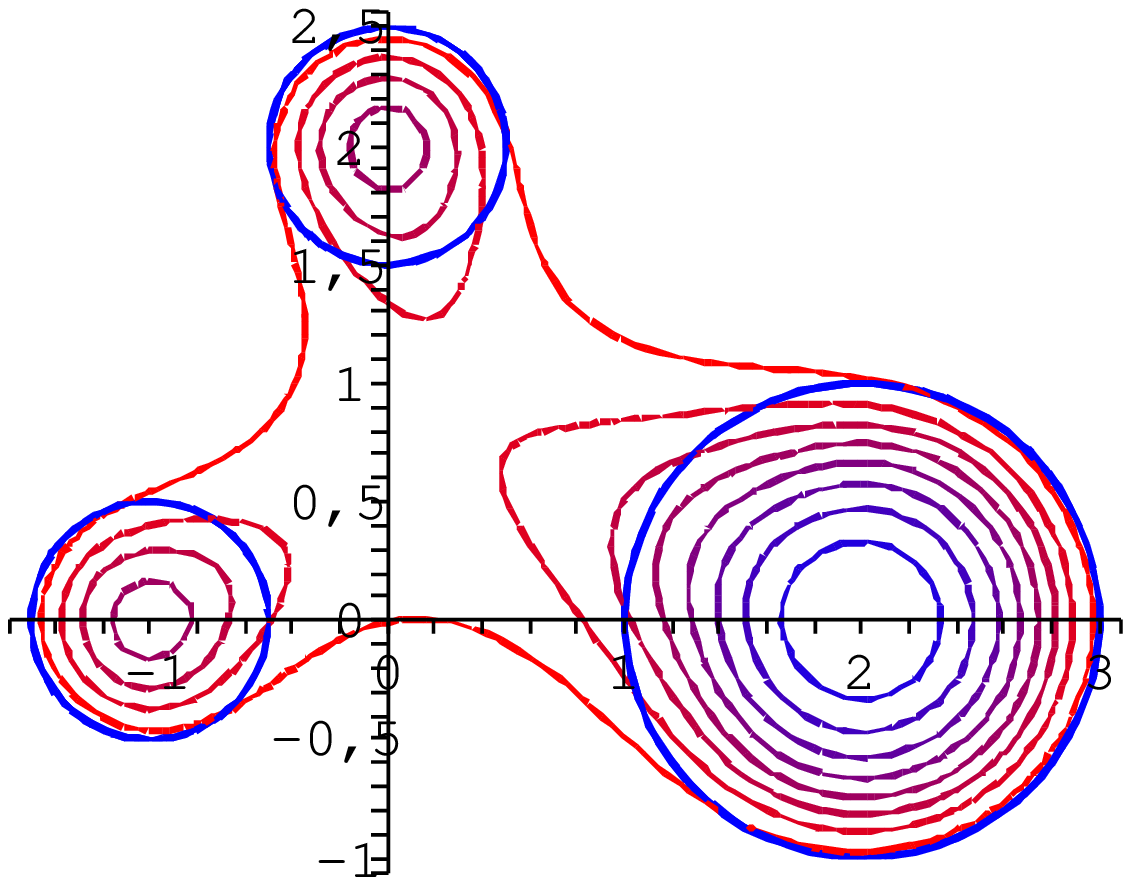}}\qquad
\fbox{\includegraphics*[scale=0.5]{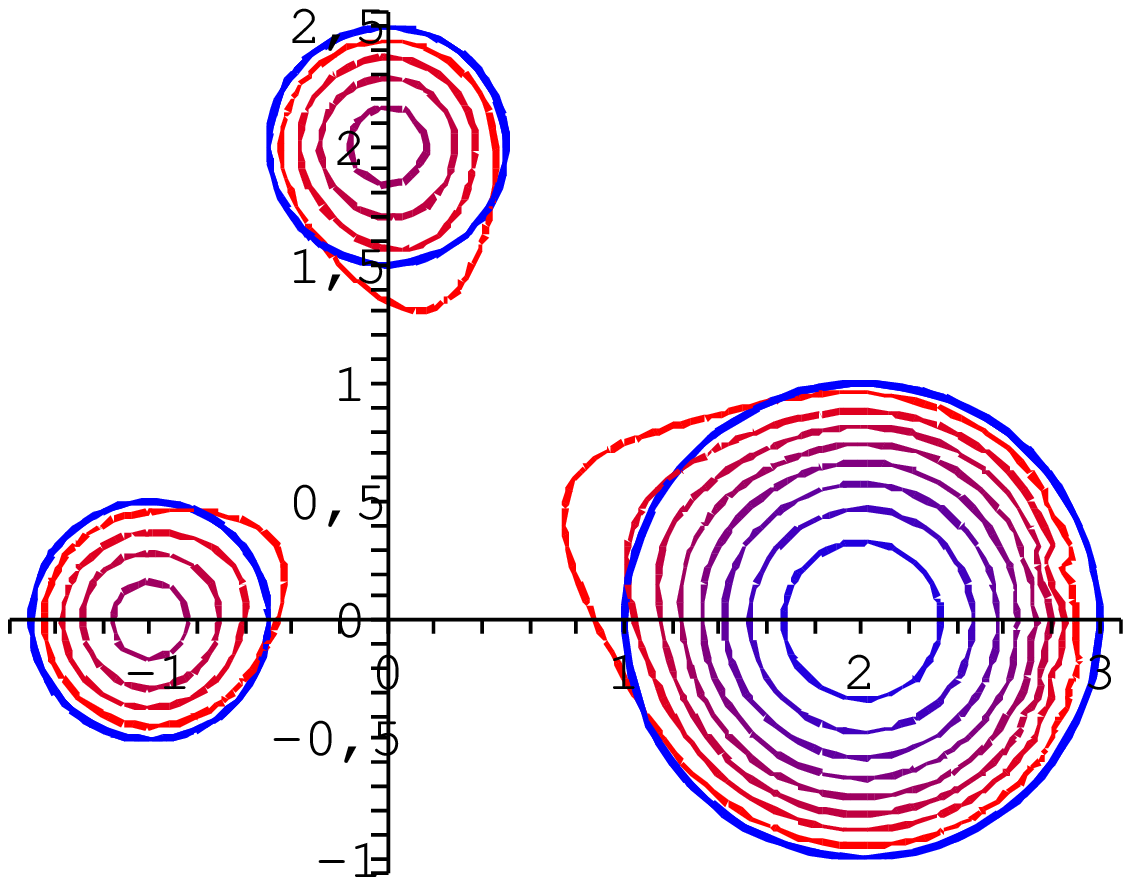}} \\ {\ } \\
\fbox{\includegraphics*[scale=0.5]{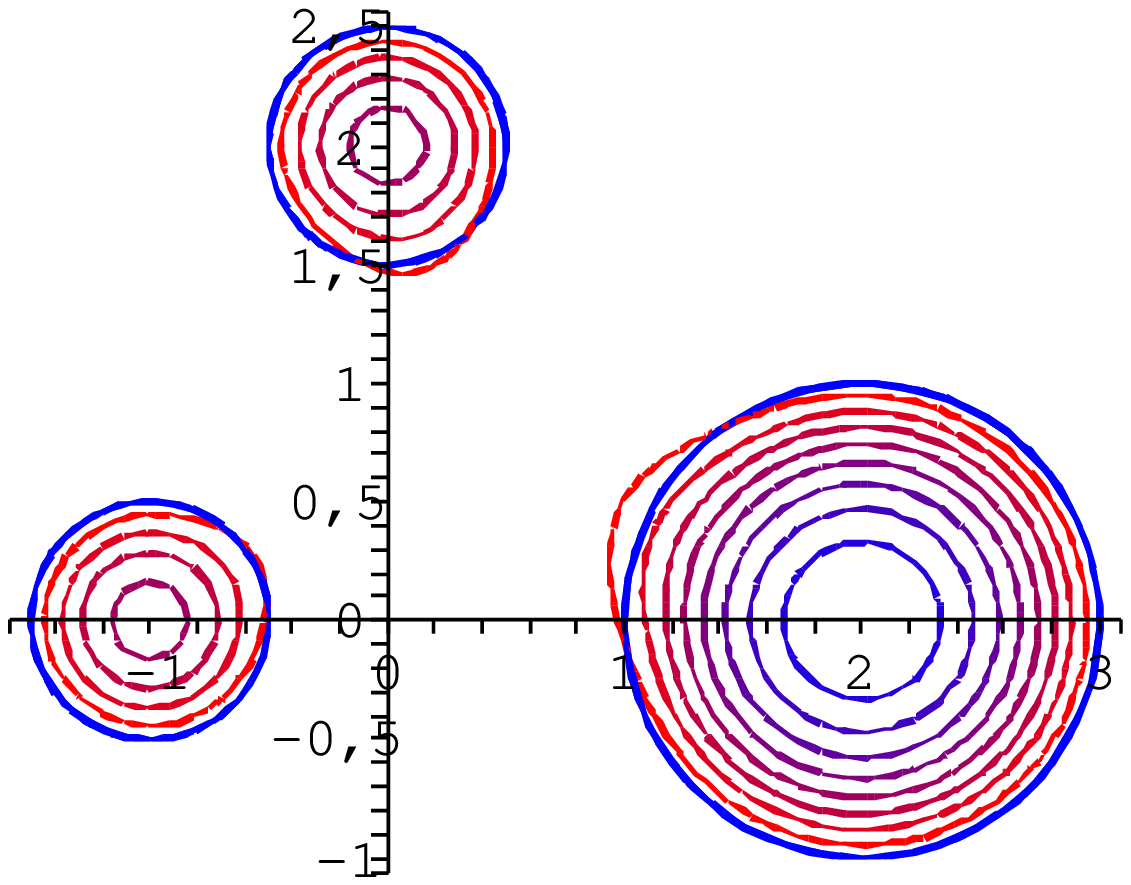}}\qquad
\fbox{\includegraphics*[scale=0.5]{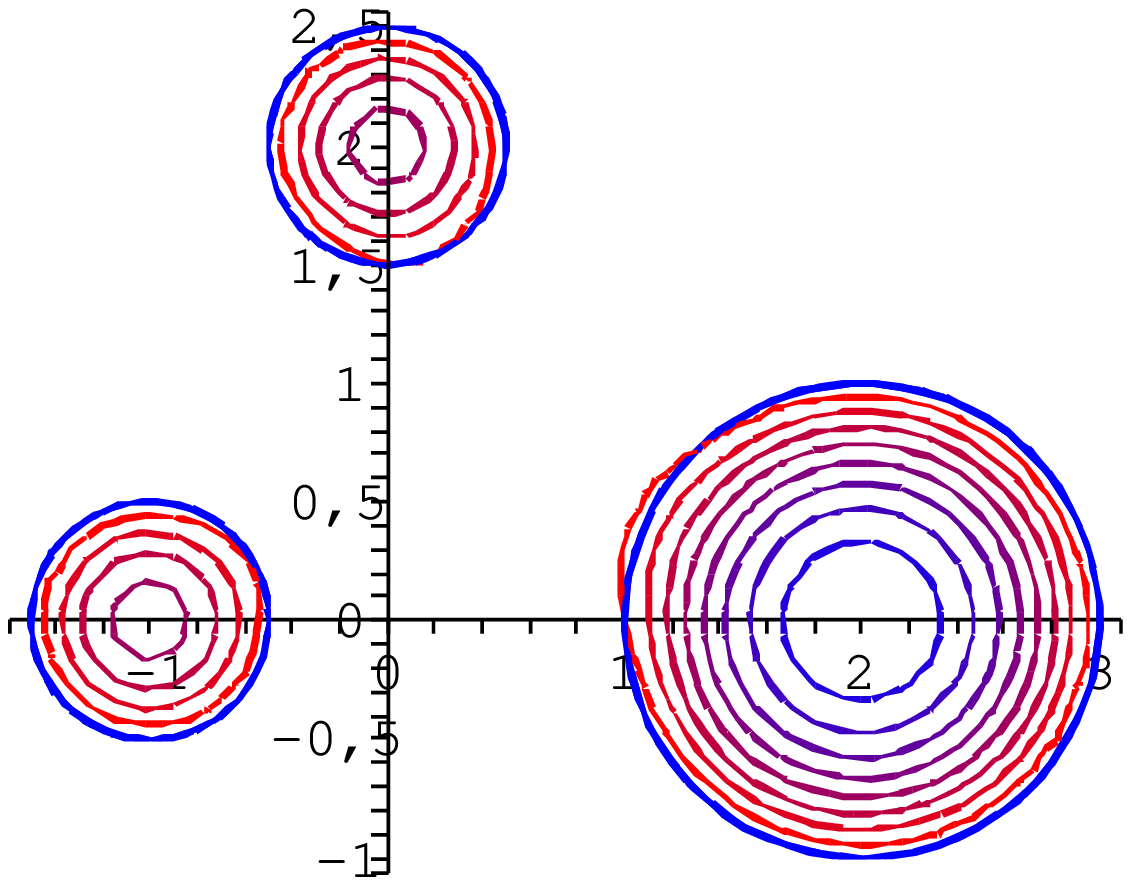}}
\end{center}
\caption{Level lines of ${\Lambda_{n}}(x+iy)$, for the values of $n$ (from left two right) $25,50,75$ and $100$, on $\{(x,y):-1\le x\le 4,-2\le y\le 3\}$, with $G$ formed by the three disjoint disks of Example~\ref{ex:3disks}.}
\label{fig:rec3disks}
\end{figure}
\clearpage

\section{Asymptotic behavior of zeros}\label{section:zeros}
\setcounter{equation}{0}
\subsection{General statements}
The first result of this section is our general theorem on the asymptotic behavior of the zeros of the Bergman polynomials $\{P_n\}_{n=1}^\infty$ on an archipelago of $N$ Jordan domains. It is established under the general assumptions made at the beginning of Section~\ref{subsec:bergman}.
In particular we note that, unlike the theory presented in Section~\ref{Estimates}, no extra smoothness is required for the boundary curves $\Gamma_j$ here. The result below, which is valid for any $N\geq 1$, requires some special attention for the single island case $N=1$.

\begin{theorem}\label{thm:ASgeneral}
Consider the following extension of the Green function $g_\Omega(\cdot, \infty)$ to all
$\overline{\mathbb{C}}$:
\begin{equation}\label{eq:h}
h(z) =
\begin{cases}
g_\Omega (z,\infty),\quad z\in\overline{\Omega},\\
-\log {\rho(K(\cdot,z))}, \quad  z\in G,
\end{cases}
\end{equation}
(recall (\ref{rho})) and define
\begin{equation}\label{eq:beta}
\beta=\beta_G:=\frac{1}{2\pi}\Delta h,
\end{equation}
where the Laplacian is taken in the sense of distributions. Let $\mathcal{C}$ denote the set of weak* cluster points of the counting measures $\{\nu_{P_n}\}_{n=1}^\infty$, i.e., the set of measures $\sigma$ for which there exists a subsequence $\mathcal{N}_\sigma\subset\mathbb{N}$ such that $\nu_{P_n}\sta\sigma$, as $n\to\infty$, $n\in\mathcal{N}_\sigma$. The following assertions hold.
\begin{itemize}
\item[(i)]
The function $h$ is harmonic in $\Omega$, subharmonic in all $\mathbb{C}$; hence $\beta$ is a positive unit measure with support contained in $\overline{G}$. In addition, if $N\geq 2$, then $h$ is continuous and bounded from below. If $N=1$, then $h$ can take the value $-\infty$ at most at two points, and outside these points $h$ is continuous.
\item[(ii)]
\begin{equation}\label{eq:potbeta}
U^\beta (z) = \log\frac{1}{{\rm cap\,}(\Gamma)}-h(z),  \,\, z\in\mathbb{C}.
\end{equation}
and balayage of $\beta$ onto $\Gamma$ gives the equilibrium measure $\mu_\Gamma$ of $\Gamma$:
\begin{equation}\label{eq:bal}
\begin{cases}
U^\beta\geq U^{\mu_\Gamma} \  { in\,\,} \mathbb{C},\\
U^\beta=U^{\mu_\Gamma} \ { in\,\,}\Omega.
\end{cases}
\end{equation}
\item[(iii)]
\begin{equation}\label{eq:limsupPn}
\limsup_{n\to\infty}\frac{1}{n}\log |P_n(z)|=h(z), \quad z\in\mathbb{C},
\end{equation}
\begin{equation}\label{eq:liminf}
\liminf_{n\to\infty}U^{\nu_{P_n}}(z)= U^\beta(z), \quad z\in\mathbb{C}.
\end{equation}
Moreover, in $\Cbar\setminus{\rm{Co}}(\Gbar)$ these equalities hold with $\limsup$ and $\liminf$ replaced by $\lim$.
\item[(iv)]
The set of cluster points $\mathcal{C}$ is nonempty, and for any  $\sigma\in\mathcal{C}$,
\begin{equation}\label{eq:usigma}
\begin{cases}
U^\sigma \geq U^\beta \quad { in\,\,} \mathbb{C},\\
U^\sigma= U^\beta\quad
{ in\,\,the\,\,unbounded\,\,component\,\,of\,\,}\overline{\mathbb{C}}\setminus{\rm supp\,}\beta.
\end{cases}
\end{equation}
\item[(v)]
The measure $\beta$ is the lower envelope of $\mathcal{C}$ in the sense that
$$
U^\beta={\rm lsc\,}{(\inf_{\sigma\in\mathcal{C}} U^\sigma)},
$$
where ``lsc" denotes lower semicontinuous regularization.
(This means that $U^\beta$ is the supremum of all lower semicontinuous
functions that are ${\leq\inf_{\sigma\in\mathcal{C}} U^\sigma}$.)
In addition, if $\mathcal{D}$ is any component of $\mathbb{C}\setminus{\rm supp\,}\beta$, then for any
$\sigma\in\mathcal{C}$ either $U^\sigma > U^\beta$ in $\mathcal{D}$ or $U^\sigma = U^\beta$ in $\mathcal{D}$; and there exists a $\sigma\in\mathcal{C}$ such that the latter holds.
\item[(vi)]
If $\mathcal{C}$ has only one element, then this is $\beta$ and
\begin{equation}
\nu_{P_n}\sta\beta, \quad n\to\infty,\ n\in\mathbb{N},
\label{eq:starfull}
\end{equation}
i.e., the full sequence converges to $\beta$.
\item[(vii)]
Assume that $\beta$ satisfies
  \begin{itemize}
  \item[(a)]
  ${\rm supp\,}\beta$ is a nullset with respect to area measure,
  \item[(b)]
  $\mathbb{C}\setminus{\rm supp\,}\beta$ is connected.
  \end{itemize}
Then $\beta$ is the unique element in $\mathcal{C}$; hence (\ref{eq:starfull}) holds. If \textup{(a)} holds and (in place of \textup{(b)})
  \begin{itemize}
  \item[(c)]
  $\mathbb{C}\setminus{\rm supp\,}\beta$ has at most two components,
  \end{itemize}
then $\beta\in\mathcal{C}$.
\end{itemize}
\end{theorem}

\begin{remark}
The measure $\beta=\beta_G$ is canonically associated to $G$ via the Bergman kernel. Constructive formulas for $\beta_G$ (or rather its potential) will be given in the proof (e.g. (\ref{eq:rhoPhi})--(\ref{eq:logrhoS})) and will be further elaborated in the examples of Section~\ref{subs:Casestudies}.
\end{remark}

\begin{remark}
Well-known properties for any $\sigma\in\mathcal{C}$ follow immediately from (ii) and (iv): That is, $U^\sigma = U^{\mu_\Gamma}$ in $\Omega$, ${\rm supp\,}\sigma\subset\overline{G}$ and balayage of $\sigma$ onto $\Gamma$ gives the equilibrium distribution $\mu_\Gamma$ (see e.g. \cite[Thm III.4.7]{ST}).
\end{remark}

\begin{remark}
We know of no example where $\beta$ isn't itself in $\mathcal{C}$. However it remains an open question whether it is always so.
\end{remark}

\begin{remark}
A measure $\beta$ satisfying (\ref{eq:bal}) together with (a) and
(b) in (vii) may be viewed as a potential theoretic skeleton for
$\mu_\Gamma$ (or ``\emph{Madonna body}'', in view of a common shape
of ${\rm supp\,}\beta$; cf. \cite{LSS, M-DSS}).
\end{remark}

\begin{remark}
When $N=1$, $h$ may take the value $-\infty$ at one or two points.
Note that, by (6.1), $h(a)=-\infty$ if and only if $K(z,a)$ is an entire function
of $z$.
With $G=\mathbb{D}$ we have $h(z)=\log |z|$, i.e., one pole for $h$.
An example with two poles is the following.

Choose a number $A>1$ and let $G$ be the image of the unit disk under the conformal map
$$
\psi (w) =\frac{1}{2}\log \frac{A+w}{A-w},
$$
the branch chosen so that $\psi(0)=0$. The inverse map is
$$
\varphi(z)=A \tanh z,
$$
which is meromorphic in the entire complex plane. Here $\psi$ maps the disk $|w|<A$ onto the strip $|{\rm Im\,}z|<\frac{\pi}{4}$. Hence $G$, which is the image of $|w|<1$, is a subdomain of that strip (a kind of an oval).

The function $\varphi$ does not attain the values $\pm A$ anywhere in the complex plane and the set
$\varphi|_\mathbb{C}^{-1}({1}/{\overline{\varphi(\zeta)}})$, which will play an important role in the proof of the theorem, may therefore be empty for up to two values of $\zeta\in G$. In fact, this occurs for $\zeta=\pm a\in G$, where $a=\frac{1}{2}\log\frac{A^2+1}{A^2-1}>0$. At these points, $h(\pm a)=-\infty$, $K(z,\pm a)=\frac{A^4-1}{\pi}e^{\pm 2z}$. One also finds that $\beta$ is a measure supported on the line segment $[-a,a]$ and hence is a Madonna body.
\end{remark}

\medskip
We call a boundary curve $\Gamma_j$ \emph{singular} if some conformal map $\varphi_j: G_j\to
\mathbb{D}$ does not extend analytically to a full neighborhood of $\overline{G}_j$, i.e., if $\rho(\varphi_j)=1$, or equivalently if $\rho(K(\cdot,z))=1$, $z\in G_j$; see (\ref{eq:KK1}) and (\ref{eq:Kjphij}). Clearly, this property is independent of the choice of the conformal map $\varphi_j$. Note that a boundary component that is not singular in the above sense still need not be fully smooth: it may be piecewise analytic but have certain kinds of corners so that  $\varphi_j$ extends analytically across $\Gamma_j$ but the extension is not univalent. This would be the case, for instance, if $G_j$ is a rectangle.

\begin{corollary}\label{cor:SS}
For each $j=1,\dots, N$, the following statements are equivalent:
\begin{itemize}
\item[(i)]
$\Gamma_j$ is singular.
\item[(ii)]
$\beta|_{\overline{G}_j}=\mu_\Gamma|_{\overline{G}_j}$.
\item[(iii)]
There is a subsequence $\mathcal{N}=\mathcal{N}_j\subset\mathbb{N}$ such that, with $V$ any neighborhood of $\overline{G}_j$ not meeting the other islands (e.g., $V=\mathcal{G}_{j,R_j}$),
\begin{equation}\label{eq:star1}
\nu_{P_n}|_V\sta\mu_\Gamma|_V, \quad n\to\infty,\ n\in\mathcal{N}.
\end{equation}
\end{itemize}
\end{corollary}

Clearly, under the conditions of the above corollary a certain proportion of the zeros of the Bergman polynomials converge to the part of the equilibrium measure located on $\Gamma_j$. By a reasoning as in deriving (\ref{eq:betaGR}) in the proof of Theorem~\ref{thm:ASgeneral} below, we conclude that this proportion is
$$
\int_{\Gamma_j}\, d\mu_\Gamma = b_j,
$$
where $b_j$ is the period in (\ref{eq:bj}). Thus, we easily deduce the following:

\begin{corollary}\label{co:zerosonGj}
If, for a particular $j=1,\dots,N$, $\Gamma_{j}$ is singular, then there is a exists a subsequence $\{P_n\}_{n\in \mathcal{N}}$, such that $P_n=Q_k R_k$,
$\mathrm{deg}(Q_k)=n_k$, where
\begin{equation}\label{eq:zerosonGj}
\frac{n_k}{n}\ \nu_{Q_k}\sta \mu_\Gamma|_{\Gamma_j}, \,\, \mathrm{as}\,\,\,n\rightarrow
\infty,\quad n\in \mathcal{N}
\end{equation}
and
$$
\frac{n_k}{n} \to b_j.
$$
\end{corollary}

As stated in (iv) of Theorem~\ref{thm:ASgeneral}, if $\sigma$ is a weak* cluster point of the  measures $\{\nu_{P_n}\}_{n=1}^\infty$ then: (a) ${\rm supp\,}\sigma\subset\overline{G}$ and (b) the balayage of $\sigma$ onto $\Gamma$ equals the equilibrium distribution $\mu_\Gamma$. The following corollary shows that the equilibrium distribution is also obtained if weak* convergence and balayage are applied in the opposite order.

\begin{corollary}\label{cor:bala}
Let ${\rm Bal\,}(\nu_{P_n})$ denote the measure obtained by balayage of $\nu_{P_n}|_G$ onto $\Gamma$ while keeping $\nu_{P_n}|_{\mathbb{C}\setminus G}$ unchanged. Then
$$
{\rm Bal\,}(\nu_{P_n})\sta \mu_\Gamma \quad {\rm as\,\,\,}
n\to\infty.
$$
\end{corollary}

\subsection{Case studies}\label{subs:Casestudies}

In this subsection we make more explicit Theorem~\ref{thm:ASgeneral} and its corollaries, and we illustrate them by means of a number of representative cases and examples.

\bigskip\noindent
\textit{Case I: Two singular boundaries.}

Here $N=2$ and $\rho (\varphi_j)=1$, $j=1,2$, for any two conformal maps
$\varphi_j:G_j\to\mathbb{D}$. By Corollary~\ref{cor:SS}, $\beta$ equals the equilibrium measure
$\mu_\Gamma$ of $G$ and there exists, for each island $G_j$, a subsequence of $\nu_{P_n}$ which converges to $\mu_\Gamma$ in a neighborhood of $\overline{G}_j$. However, we do not know whether
there necessarily exists a common subsequence for the two islands.

\bigskip\noindent
\textit{Case II: One singular boundary and one analytic boundary.}

Assume that $\Gamma_1$ is singular and $\Gamma_2$ is analytic. Then in terms of two specific conformal maps $\varphi_j:G_j\to \mathbb{D}$, $j=1,2$: (a) $\varphi_1$ has no analytic continuation beyond $\Gamma_1$, (b) $\varphi_2$ extends analytically as a univalent function to some domain containing  $\overline{G}_2$. 
Since $\Gamma_2$ is an analytic Jordan curve, it possesses a Schwarz function, which is given by
$$
S_2(z)=\overline{\varphi_2^{-1}({1}/{\overline{\varphi_2(z)}})}.
$$
In order to formulate a particular statement we assume further that $\varphi_2$ remains analytic and univalent throughout $\mathcal{G}_{2,R'}$. This implies that $g_\Omega(\cdot,\infty)$ extends by Schwarz reflection up
to the level line $L_{2,\frac{1}{R'}}$; see (\ref{uS}) and the terminology in Example~\ref{ex:Greentextened}. Moreover, the domain
$$
D_2:=\mathcal{G}_{2,R'}\setminus\overline{\mathcal{G}_{2,\frac{1}{R'}}}
$$
is connected and is a domain of involution of the Schwarz reflection $z\mapsto \overline{S_2(z)}$.

Set
$$
E={G}_1 \cup {\mathcal{G}_{2,\frac{1}{R'}}}.
$$
It follows that the multi-valued function
\begin{equation}\label{phidefinition}
\widehat{\Phi}(z):=\left\{\begin{array}{cl}
\Phi(z) & \mathrm{if}\ z\in\overline{\mathbb{C}}\setminus G, \\
          1\big/\conj{\Phi\left(\overline{S_2(z)}\right)}
        & \mathrm{if}\ z\in G_2\setminus{\mathcal{G}_{2,\frac{1}{R'}}}.\\
\end{array}\right.
\end{equation}
is (locally) analytic in $\mathbb{C}\setminus \overline{E}$ and (locally) continuous on ${\mathbb{C}}\setminus E$. It also follows from the expression (\ref{eq:rhoPhi2}) of $\rho(K(\cdot,z))$ appearing in the proof of Theorem~\ref{thm:ASgeneral}, by taking into account (\ref{rhom}) and (\ref{eq:rhoRj}), that
\begin{equation}\label{eq:rhoCaseII}
\rho(K(\cdot,z))=\left\{\begin{array}{cl}1 & {\rm if}\ z\in G_1, \\
\exp\{-g_\Omega(z,\infty)\} & {\rm if}\ z\in G_2\setminus\mathcal{G}_{2,\frac{1}{R'}},\\
R' & {\rm if} \ z\in\mathcal{G}_{2,\frac{1}{R'}}.\\
\end{array}\right.
\end{equation}

The relations in (\ref{phidefinition}) and (\ref{eq:rhoCaseII}) yield at once, in view of Proposition~\ref{pro:nthroot} and Corollary~\ref{cor2}, the $n$-th root asymptotic behavior of $\{P_n\}_{n=1}^\infty$ in $\mathbb{C}$:
\begin{equation}\label{nose25}
\limsup_{n\to\infty}|P_n(z)|^{1/n}= \left\{\begin{array}{cl}1 & {\rm if}\ z\in G_1, \\
|\widehat{\Phi}(z)| & {\rm if}\ z ,\in\overline{\mathbb{C}}\setminus E,\\
\frac{1}{R'} & {\rm if} \ z\in\mathcal{G}_{2,\frac{1}{R'}}.\\
\end{array}\right.
\end{equation}
In addition, these relations provide more detailed information for the potential $U^\beta$ of the canonical measure $\beta$, and thus for the counting measures $\{\nu_{P_n}\}_{n=1}^\infty$.

\begin{corollary}\label{cor:AS}
Under the assumption and notations of Case II, we have:
\begin{equation}\label{eq:betainAS}
U^\beta(z)=\left\{\begin{array}{ll}
\log\frac{1}{{\rm cap\,}(\Gamma)}  & { if}\ z\in G_1, \\
\log\frac{1}{{\rm cap\,}(\Gamma)}-g_\Omega(z,\infty)& { if}\ z\in{\mathbb{C}}\setminus E,\\
\log\frac{R'}{{\rm cap\,}(\Gamma)} & { if} \ z\in\mathcal{G}_{2,\frac{1}{R'}}.\\
\end{array}\right.
\end{equation}
In particular,
\begin{itemize}
\item[(i)]
${\rm supp\,}\beta=\partial{E}$.
\item[(ii)]
For any weak* cluster point $\sigma$ of $\{\nu_{P_n}\}$,
${\rm supp\,}\sigma\subset\overline{E}$,
and
\begin{equation*}
U^\sigma(z)=U^\beta(z),\quad z\in\overline{\mathbb{C}}\setminus E.
\end{equation*}
\item[(iii)]
There is a subsequence $\mathcal{N}\subset\mathbb{N}$ such that, with $V$ any neighborhood of $\overline{G}_1$ or $\overline{\mathcal{G}}_{2,\frac{1}{R'}}$ not meeting the other island,
\begin{equation}
\nu_{P_n}|_V\sta\beta|_V, \quad n\to\infty,\ n\in\mathcal{N}.
\end{equation}
Hence, every point of $\partial E=\Gamma_1\cup L_{2,\frac{1}{R'}}$ belongs to ${\rm supp\,}\sigma$, for some weak* cluster point $\sigma$ of $\{\nu_{P_n}\}_{n=1}^\infty$.
\end{itemize}
\end{corollary}

\begin{figure}[h]
\begin{center}
\includegraphics*[width=0.5\linewidth]{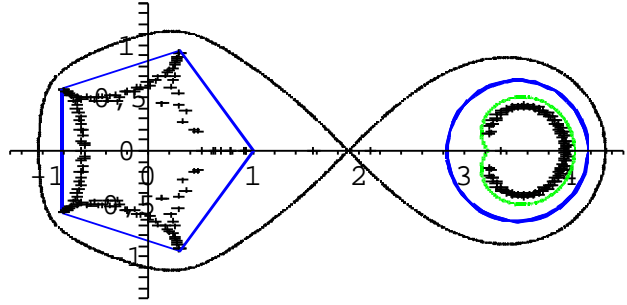}
\end{center}
\caption{Zeros of Bergman polynomials $P_n$ of
Example~\ref{ex:diskpen}, for $n=80$, $90$ and
$100$.}\label{fig:2.diskpen}
\end{figure}

The corollary is illustrated in the following example.
\begin{example}\label{ex:diskpen}
Bergman polynomials for $G=G_1\cup G_2$, with $G_1$ the canonical pentagon with vertices at the fifth roots of unity and $G_2=\{z:|z-7/2|<2/3\}$.
\end{example}
The zeros of the associated Bergman polynomials $P_n$, for $n=80$, $90$ and $100$ are shown in Figure~\ref{fig:2.diskpen}. In the same figure we also depict the critical line $L_{R'}$ and the curve $L_{2,\frac{1}{R'}}$. Note that $L_{2,\frac{1}{R'}}$ is simply the inverse image  of $L_{2,R^\prime}$ with respect to the circle $\{z: |z-7/2|=2/3\}$.

\bigskip\noindent
\textit{Case III: Two analytic boundary curves.}
This is the case $N=2$, where both $\Gamma_1$ and $\Gamma_2$ are analytic curves.

\begin{example}\label{ex:2disks}
Bergman polynomials for the union of the disks: $G_1=\{z:|z+2|<1\}$ and $G_2:=\{z:|z-3|<2/3\}$.
\end{example}
Let $S_1$ and $S_2$ denote the Schwarz functions defined by $\Gamma_1$ and $\Gamma_2$. (Note that the Schwarz function for the circle $\{z:|z-a|=r\}$ is simply $S(z)=r^2/(z-a)+\overline{a}$.)
Clearly, the Green function $g_\Omega$ extends by Schwarz reflection to the set
\begin{equation}\label{eq:gextD}
D=({\mathcal G}_{1,R'}\setminus\overline{\mathcal{G}_{1,\frac{1}{R'}}})\cup
({\mathcal G}_{2,R'}\setminus\overline{\mathcal{G}_{2,\frac{1}{R'}}}),
\end{equation}
and the multi-valued function
\begin{equation}\label{eq:Phimulti}
\widehat{\Phi}(z):=\left\{\begin{array}{cl}
  \Phi(z) & \mathrm{if}\ z\in\overline{\mathbb{C}}\setminus G, \\
1\big/\conj{\Phi\left(\overline{S_j(z)}\right)} & \mathrm{if}
  \ z \in G_j\setminus{\mathcal{G}_{j,\frac{1}{R'}}},\quad j=1,2,\\
\end{array}\right.
\end{equation}
is (locally) analytic in $\mathbb{C}\setminus \overline{E}$ and
(locally) continuous on ${\mathbb{C}}\setminus E$, where now
\begin{equation*}
E=\mathcal{G}_{1,\frac{1}{R'}}\cup \mathcal{G}_{2,\frac{1}{R'}}.
\end{equation*}
As in Case II, the extensions of $g_\Omega(z,\infty)$ and $\Phi(z)$ lead to the expressions
\begin{equation}\label{eq:rhoCaseIII}
\rho(K(\cdot,z))=\left\{\begin{array}{cl}
\exp\{-g_\Omega(z,\infty)\} & {\rm if}\ z\in G \setminus E,\\
R' & {\rm if} \ z\in E,\\
\end{array}\right.
\end{equation}
\begin{equation}\label{eq:PnEx2}
\limsup_{n\to\infty}|P_n(z)|^{1/n}= \left\{\begin{array}{cl}
|\widehat{\Phi}(z)| & {\rm if}\
 z\in\overline{\mathbb{C}}\setminus {E},\\
 \frac{1}{R'} & {\rm if} \ z \in {E},\\
\end{array}\right.
\end{equation}
and in parallel with Corollary~\ref{cor:AS}, to the conclusion
\begin{equation}\label{eq:betain2D}
U^\beta(z)=\left\{\begin{array}{ll}
\log\frac{1}{{\rm cap\,}(\Gamma)}-g_\Omega(z,\infty)& {\rm if}\ z\in{\mathbb{C}}\setminus E,\\
\log\frac{R'}{{\rm cap\,}(\Gamma)} & {\rm if} \ z\in E,\\
\end{array}\right.
\end{equation}
${\rm supp\,}\beta=\partial{E}$ and that every point of $\partial E=L_{1,\frac{1}{R'}}\cup L_{2,\frac{1}{R'}}$ attracts zeros of the sequence $\{P_n\}_{n=1}^\infty$. Furthermore, since the unbounded domains $\overline{\mathbb{C}}\setminus\overline{E}$ and $\Omega_{\frac{1}{R'}}$ coincide, it follows from (\ref{eq:gOmROm}), (\ref{eq:CapLRGam}) and (\ref{eq:betain2D}) that the same is true for the potentials $U^\beta$ and $U^{\mu_{\partial E}}$ in $\mathbb{C}$. Hence, the canonical measure $\beta$ is the equilibrium measure of $\partial E$. Therefore, by applying Corollary~\ref{cor:SS}~(ii) (with $E$ in the place of $G$), we conclude that for $j=1,2$, there is a subsequence $\mathcal{N}=\mathcal{N}_j\subset\mathbb{N}$ such that, with $V$ any neighborhood of $\overline{\mathcal{G}_{j,\frac{1}{R'}}}$ not meeting the other island,
\begin{equation}
\nu_{P_n}|_V\sta\mu_{\partial E}|_V, \quad n\to\infty,\ n\in\mathcal{N}.
\end{equation}

\begin{figure}[h]
\begin{center}
\includegraphics*[scale=0.7]{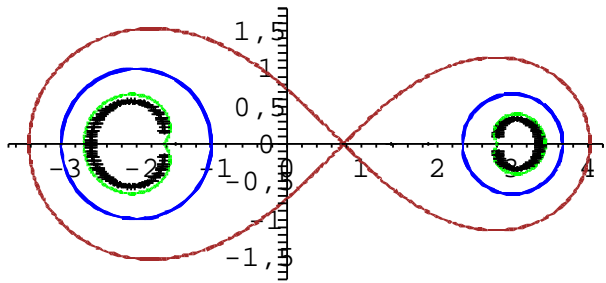}
\end{center}
\caption{Zeros of Bergman polynomials $P_n$ of
Example~\ref{ex:2disks}, for $n=140$, $150$ and
$160$.}\label{fig:2.2disks}
\end{figure}

The zeros of the associated Bergman polynomials $P_n$, for $n=140$, $150$ and $160$ are shown in
Figure~\ref{fig:2.2disks}. In the same figure we also depict the critical line $L_{R^\prime}$ and the curves $L_{1,\frac{1}{R'}}$ and $L_{2,\frac{1}{R'}}$. Note that $L_{j,\frac{1}{R'}}$ is the inverse image  of  $L_{j,R^\prime}$  with respect to the circle $\Gamma_j$, $j=1,2$.

\begin{figure}[h]
\begin{center}
\fbox{\includegraphics*[scale=0.6]{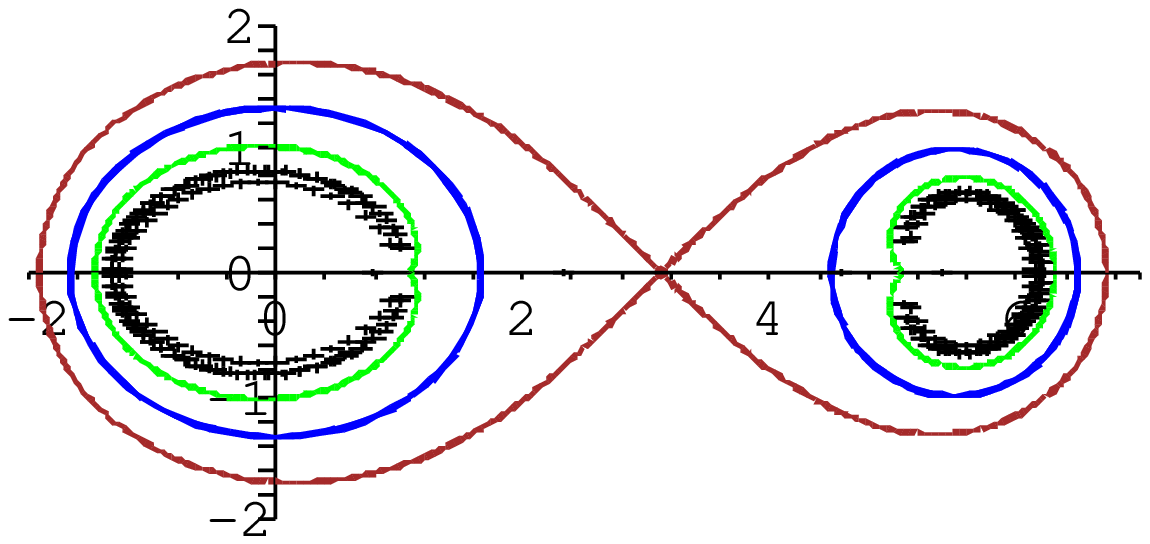}}\quad
\fbox{\includegraphics*[scale=0.58]{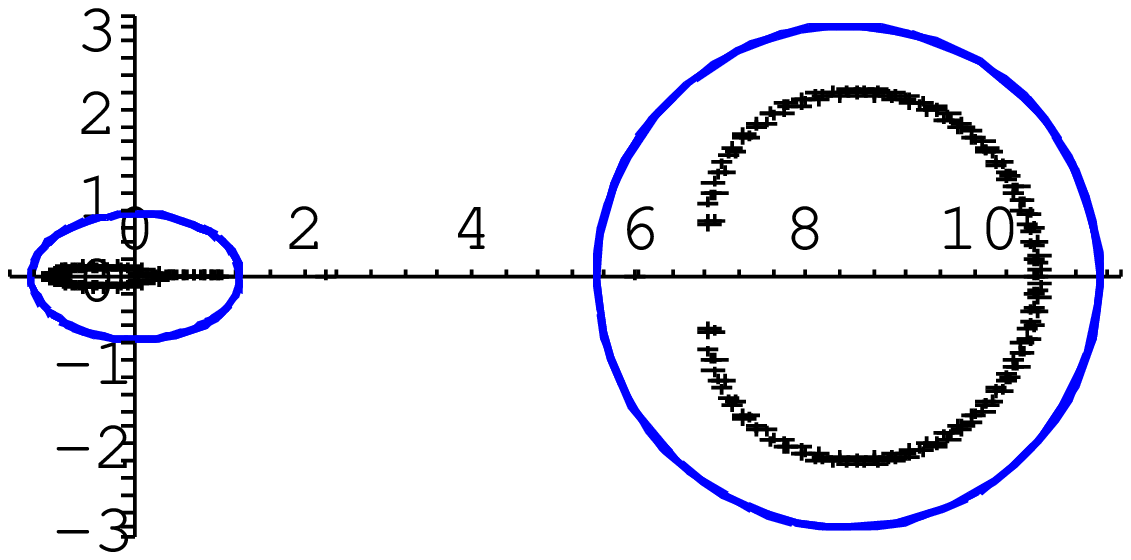}} \\
(i)$\qquad\qquad\qquad\qquad\qquad\qquad\qquad\qquad\qquad\qquad$ (ii)\\ {\ } \\
\fbox{\includegraphics*[scale=0.9]{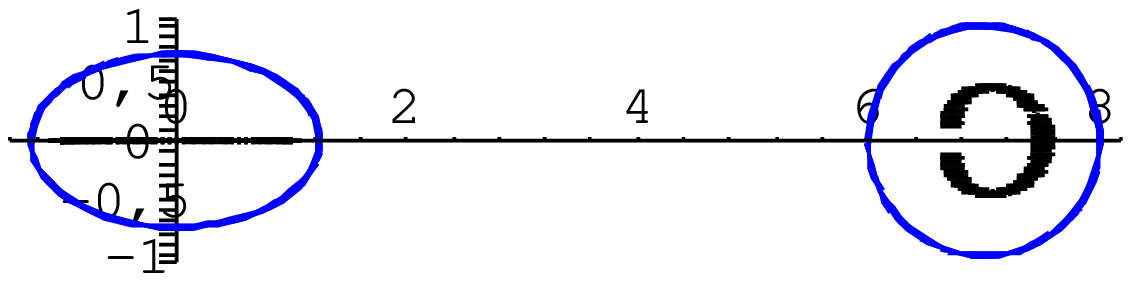}}\\
(iii)
\end{center}
\caption{Zeros of Bergman polynomials $P_n$ of
Example~\ref{ex:diskell}, for $n=80$, $90$ and
$100$.}\label{fig:2.diskell} \end{figure}

\begin{example}\label{ex:diskell}
Bergman polynomials for the union of an ellipse and a disk.
\end{example}

In Figure~\ref{fig:2.diskell} we plot the zeros of the Bergman polynomials $P_n$, for $n=80$, $90$ and $100$ of an ellipse (domain $G_1$) and a disk (domain $G_2$), in relative positions chosen to illustrate further the theory given above. To this end, let $S_1$ and $S_2$ denote the Schwarz function associated with the ellipse, respectively the circle. The three ellipses pictured in Figure~\ref{fig:2.diskell} have all focal segment [-1,1] and canonical equation
\begin{equation*}
\frac{x^2}{a^2} + \frac{y^2}{b^2} = 1,
\end{equation*}
with $a=5/3$, $b=4/3$, in (i) and $a=5/4$, $b=3/4$, in both (ii) and (iii).

For such ellipses the associated Schwarz function is given by
$$
S_1(z) = (2a^2-1)z - 2ab \sqrt{z^2-1},
$$
and the focal segment $[-1,1]$ is reflected to the confocal ellipse ${x^2}/{A^2}+{y^2}/{B^2} = 1,$ where $A=2a^2-1$ and $B=2ab$.
We denote by
$$
D_1 = \{ (x,y): \  \frac{x^2}{A^2} + \frac{y^2}{B^2} <1\} \setminus [-1,1],
$$
the maximal domain of involution for the Schwarz reflection and by $\gamma$ the outer boundary
of $D_1$, i.e.,
$$
\gamma =\{(x,y): \,\frac{x^2}{A^2}+\frac{y^2}{B^2}=1\}.
$$
Also, if $G_2$ is a disk centered at $z=z_0$, the reflection $z\mapsto\overline{S_2(z)}$ is an involution on the domain $D_2=\mathbb{C}\setminus\{z_0\}$.

The situations illustrated in Figure~\ref{fig:2.diskell} represent the three possible relative positions between the loop $L_{1,R'}$ of the singular level set $L_{R^\prime}$ and $\gamma$:
\begin{itemize}
\item
Figure~\ref{fig:2.diskell}~(i) corresponds to the case that $L_{1,R'}$ is interior to $\gamma$,
\item
Figure~\ref{fig:2.diskell}~(ii) corresponds to the case that $L_{1,R'}$ intersects $\gamma$,
\item
Figure~\ref{fig:2.diskell}~(iii) corresponds to the case that the ellipse $\gamma$ is interior to $L_{1,R'}$.
\end{itemize}

By specializing Theorem~\ref{thm:ASgeneral} to this example, we can conclude the following:

Case (i) is completely analogous to Example~\ref{ex:2disks}. That is,  ${\rm supp\,}\beta=\partial{E}=L_{1,\frac{1}{R'}}\cup L_{2,\frac{1}{R'}}$ and every point of $\partial E$ attracts zeros of the sequence $\{P_n\}_{n=1}^\infty$. More precisely, $\beta=\mu_{\partial E}$ and for any $j=1,2$, there exists a subsequence $\mathcal{N}=\mathcal{N}_j\subset\mathbb{N}$ such that, with $V$ any neighborhood of $\overline{\mathcal{G}_{j,\frac{1}{R'}}}$ not meeting the other island,
\begin{equation*}
\nu_{P_n}|_V\sta\mu_{\partial E}|_V, \quad n\to\infty,\ n\in\mathcal{N}.
\end{equation*}

In case (ii), the support of the canonical measure $\beta$ consists of three parts: the inverse image $L_{2,\frac{1}{R'}}$ of $L_{2,R'}$ with respect to the circle $\Gamma_2$, the reflection of $L_{1,R'}\cap D_1$ with respect to the ellipse $\Gamma_1$ and the part $[s,1]$ of the focal segment $[-1,1]$ of the ellipse that lies exterior to this reflection. In addition, every point of ${\rm supp\,}\beta$ attracts zeros of the sequence $\{P_n\}_{n=1}^\infty$.

Finally in case (iii), ${\rm supp\,}\beta=[-1,1]\cup L_{2,\frac{1}{R'}}$. Thus
$\mathbb{C}\setminus {\rm supp\,}\beta$ has exactly two components and it follows from (vii) of Theorem~\ref{thm:ASgeneral} that there exists is a subsequence $\mathcal{N}\subset\mathbb{N}$ such that
\begin{equation}
\nu_{P_n}\sta\beta, \quad n\to\infty,\ n\in\mathcal{N}.
\end{equation}

\bigskip\noindent
\textit{Case IV: One piecewise analytic non-singular boundary and one
analytic boundary curve.}

Assume that $\Gamma_2$ is analytic and $\Gamma_1$ is piecewise analytic and non-singular. By the latter we mean that any conformal map $\varphi_1:G_1\to \mathbb{D}$ has an analytic continuation to a neighborhood of $\overline{G}_1$,  but this continuation is not univalent in any neighborhood of $\overline{G}_1$. This occurs, for example, if $\Gamma_1$ consists of circular arcs and/or straight lines and all its interior corners are of the form $\pi/m$, $m\geq 2$ an integer.

\begin{example}\label{ex:HDdisk}
Bergman polynomials for the union of the half-disk
$
G_1=\{z:|z|<1,\,\Re(z)>0\}
$
and the disk
$
G_2=\{z:|z-3|<2/3\}.
$
\end{example}
\begin{figure}
\begin{center}
\includegraphics*[scale=0.7]{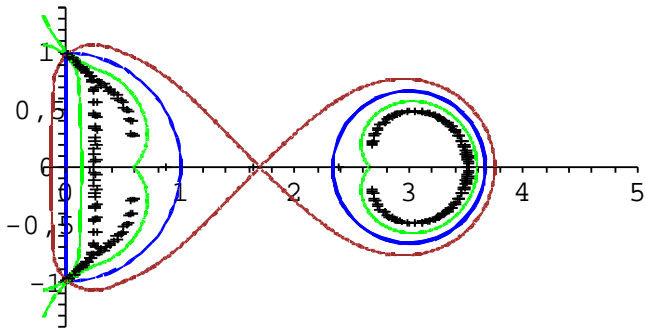}
\end{center}
\caption{Zeros of Bergman polynomials $P_n$ of
Example~\ref{ex:HDdisk}, for $n=80$, $90$ and
$100$.}\label{fig:2.HDdisk}
\end{figure}
In Figure~\ref{fig:2.HDdisk} we plot the zeros of the Bergman polynomials $P_n$ of $G$, for $n=80$, $90$ and $100$. In addition we depict:
\begin{itemize}
\item The critical level line $L_{R^\prime}$ of the Green function $g_\Omega(z,\infty)$.
\item
The part of the reflection (we denote it by $\Gamma_1^\prime$) of $L_{1,R^\prime}$ with respect to $\Gamma_1$ which lies in $G_1$.
\item
The inverse image $L_{2,\frac{1}{R'}}$ of $L_{2,R^\prime}$ with respect to the circle $\Gamma_2$.
\end{itemize}

By considering the symmetric and inverse images of the interior points of $G_1$ with respect to the two arcs forming $\Gamma_1$, in conjunction with the harmonic extension of the Green function inside $G_1$ defined by the Schwarz functions of these arcs, it is not difficult to see that the support of the canonical measure $\beta$ consists of three parts: the loop $\Gamma_1^\prime$ and two (symmetric) arcs that join together each one of the points $i$ and $-i$ with the nearest corner of $\Gamma_1^\prime$. In addition, every point of ${\rm supp\,}\beta$ attracts zeros of the sequence $\{P_n\}_{n=1}^\infty$.

\begin{example}\label{ex:disklens}
Bergman polynomials for the union of the symmetric lens domain $G_1$ formed by two circular arcs meeting at $-i$ and $i$ with interior angles $\pi/4$ and the disk $G_2=\{z:|z-5/2|<2/3\}$.
\end{example}
\begin{figure}
\begin{center}
\includegraphics*[scale=0.7]{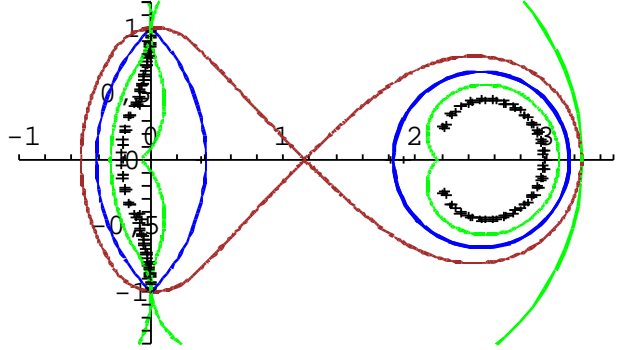}
\end{center}
\caption{Zeros of Bergman polynomials $P_n$ of
Example~\ref{ex:disklens}, for $n=80$, $90$
and $100$.}\label{fig:2.disklens}
\end{figure}
In Figure~\ref{fig:2.disklens} we plot the zeros of the Bergman polynomials $P_n$ of $G$, for $n=80$, $90$ and $100$. In addition we depict:
\begin{itemize}
\item
The critical level line $L_{R^\prime}$ of the Green function $g_\Omega(z,\infty)$.
\item
The part of the reflection (we denote it by $\Gamma_1^\prime$) of $L_{1,R^\prime}$ with respect to $\Gamma_1$ which lies in $G_1$.
\item
The inverse image $L_{2,\frac{1}{R'}}$ of $L_{2,R^\prime}$ with respect to the circle $\Gamma_2$.
\end{itemize}
As it is expected, identical conclusions to those of Example~\ref{ex:HDdisk} regarding the properties of the support of the canonical measure $\beta$ hold here.

\bigskip
\noindent
\textit{Case V: Three analytic boundaries.}

\begin{example}\label{ex:3disks}
Bergman polynomials for the union of the three disks
$G_1=\{z:|z+1|<1/2\}$, $G_2=\{z:|z-2|<1\}$ and $G_3=\{z:|z-2i|<1/2\}$.
\end{example}
In this example we have two critical Green level lines, $L_{R'}$ and $L_{R''}$, where $R'=R_2=R_3$ and $R''=R_1$. (See Figure~\ref{fig:greenlines} which depicts the present example.) On setting
$$
E'=\mathcal{G}_{2,\frac{1}{R'}}\cup \mathcal{G}_{3,\frac{1}{R'}}\quad\textup{and}\quad
E''=\mathcal{G}_{1,\frac{1}{R''}},
$$
we have
\begin{equation}\label{eq:Phiin3D}
\limsup_{n\to\infty}|P_n(z)|^{1/n}= \left\{\begin{array}{cl}
|\widehat{\Phi}(z)| & {\rm if}\
 z\in\overline{\mathbb{C}}\setminus {(E'\cup E'')},  \\
\frac{1}{R'} & {\rm if} \ z \in {E'},\\
\frac{1}{R''} & {\rm if} \ z \in {E''},\\
\end{array}\right.
\end{equation}
where $\widehat{\Phi}(z)$ is the multi-valued function defined as in (\ref{eq:Phimulti}), with $j=1,2,3$. From (\ref{eq:Phiin3D}) and (\ref{cor2eq}) conclusions can be drawn about the canonical measure $\beta$. In particular we note that ${\rm supp\,}\beta=\partial E'\cup\partial E''=L_{1,\frac{1}{R_1}}\cup L_{2,\frac{1}{R_2}}\cup L_{3,\frac{1}{R_3}}$ and that every point of $\partial E'\cup\partial E''$ attracts zeros of the sequence $\{P_n\}_{n=1}^\infty$.

\begin{figure}
\begin{center}
\includegraphics*[scale=0.5]{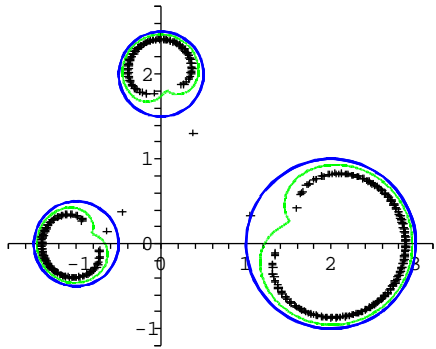}
\end{center}
\caption{Zeros of Bergman polynomials $P_n$ of
Example~\ref{ex:3disks}, for $n=80$, $90$
and $100$.}\label{fig:3disks}
\end{figure}
In Figure~\ref{fig:3disks} we plot the zeros of the Bergman polynomials $P_n$ of $G$, for $n=80$, $90$ and $100$. In order to illustrate the above observations regarding the zero distribution we also depict the inverse image $L_{j,\frac{1}{R_j}}$ of $L_{j,R_j}$ with respect to the circle $\Gamma_j$, $j=1,2,3$.

We end this section by noting that the critical level curves of the Green function depicted in the plots above, were computed by a simple modification of the MATLAB code {\tt manydisks.m} of Trefethen~\cite{Tre}. The original code {\tt manydisks.m} is designed for archipelagoes formed by circles; see also Remark~\ref{rem:manydisks}.

\section{An example: lemniscate islands}\label{section:lemniscate}
\setcounter{equation}{0}

Let $G:=\{z:\, |z^m-1|<r^m\}$, $m\ge 2$ an integer and $0<r<1$. Then $G$ consists of $m$ islands $G_1,G_2,\ldots,G_m$, where
\begin{equation}\label{eq:Gjcon}
G_j \mbox{ contains } e^{2\pi ji/m}, \quad j=1,2,\ldots,m.
\end{equation}
Let $P_n(z)=\lambda_n z^n+\cdots$ denote the (orthonormal) Bergman polynomial
of degree $n$ for the archipelago $G$, and write
\begin{equation*}
n=km+s,\quad 0\le s\le m-1.
\end{equation*}
By the rotational symmetry of $G$ and the uniqueness of the Bergman polynomials it is easy to see that
\begin{equation}  \label{eq:Pkm+s}
P_{km+s}(z)=z^sQ_{k,s}(z^m),\quad \deg Q_{k,s}=k.
\end{equation}
Then
\begin{equation}  \label{eq:pkm+s}
p_{km+s}(z):=\frac{P_{km+s}(z)}{\lambda_{km+s}}=z^sq_{k,s}(z^m)=z^{km+s}+%
\cdots,
\end{equation}
are the monic Bergman polynomials. Our first result concerns the asymptotic behavior of the
leading coefficient $\lambda_n$.

\begin{proposition}\label{pro:lemn1}
For each $s=0,1,\ldots,m-1$ there holds
\begin{equation}  \label{eq:proplemn1}
\lim_{k\to\infty}\lambda_{km+s}\,
r^{km+s+1}\sqrt{\frac{\pi}{km+s+1}}=\frac{1}{r^{m-s-1}}.
\end{equation}
\end{proposition}

\begin{remark}
Note that $r=\mathrm{cap}(\overline{G})={\mathrm{cap}}(\Gamma)$, where as above $\Gamma=\partial G$. Thus the sequence
\begin{equation*}
{\lambda_n\,{\mathrm{cap}}(\Gamma)^{n+1}}\sqrt{\frac{\pi}{n+1}},\quad
n\in\mathbb{N},
\end{equation*}
has exactly $m$ limit points, $\frac{1}{r^{m-1}}, \frac{1}{
r^{m-2}},\ldots,\frac{1}{r}, 1$.
\end{remark}

\begin{figure}[h]
\begin{center}
\mbox{\includegraphics*[scale=0.55]{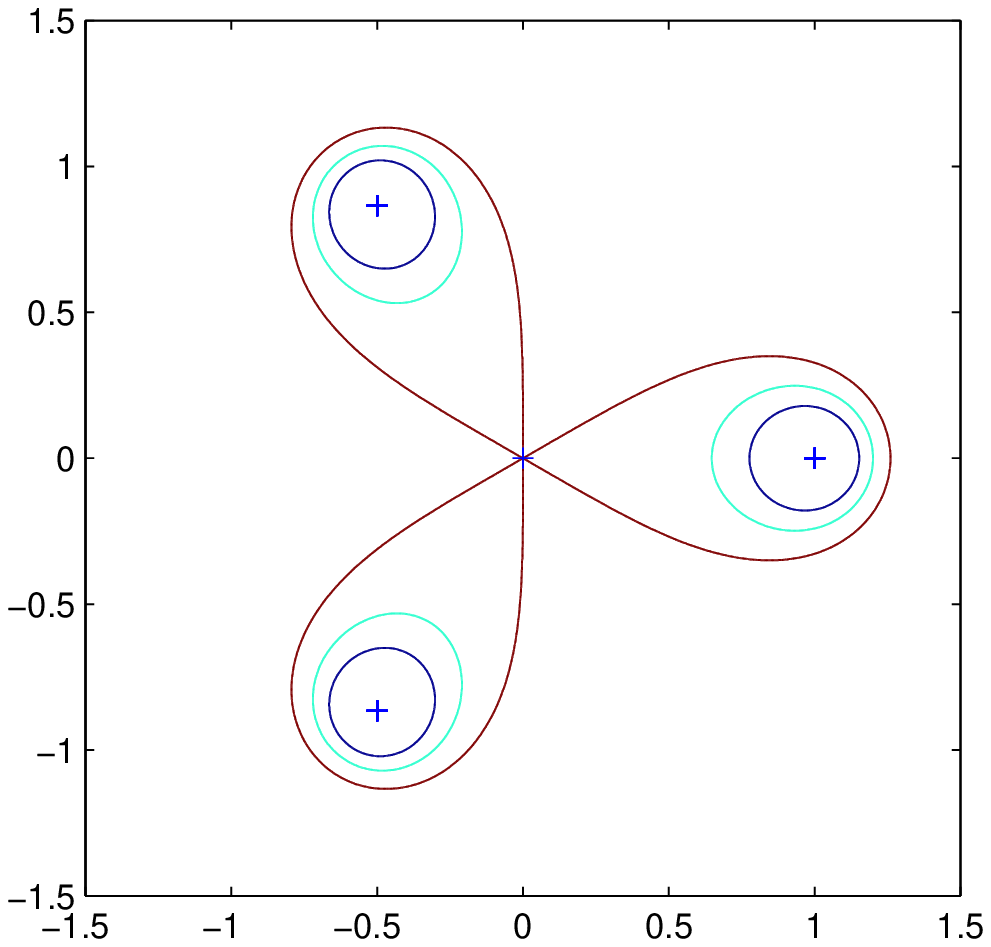}}
\mbox{\includegraphics*[scale=0.55]{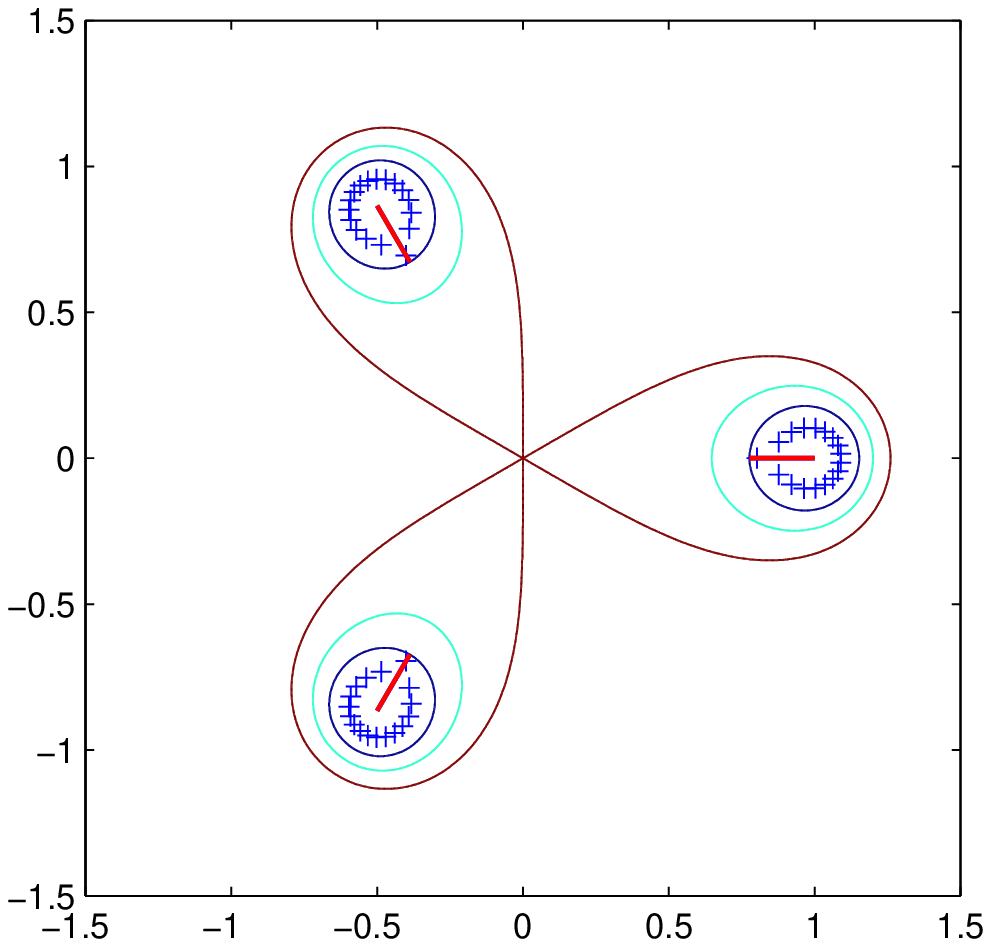}} \\ {\ } \\
\mbox{\includegraphics*[scale=0.55]{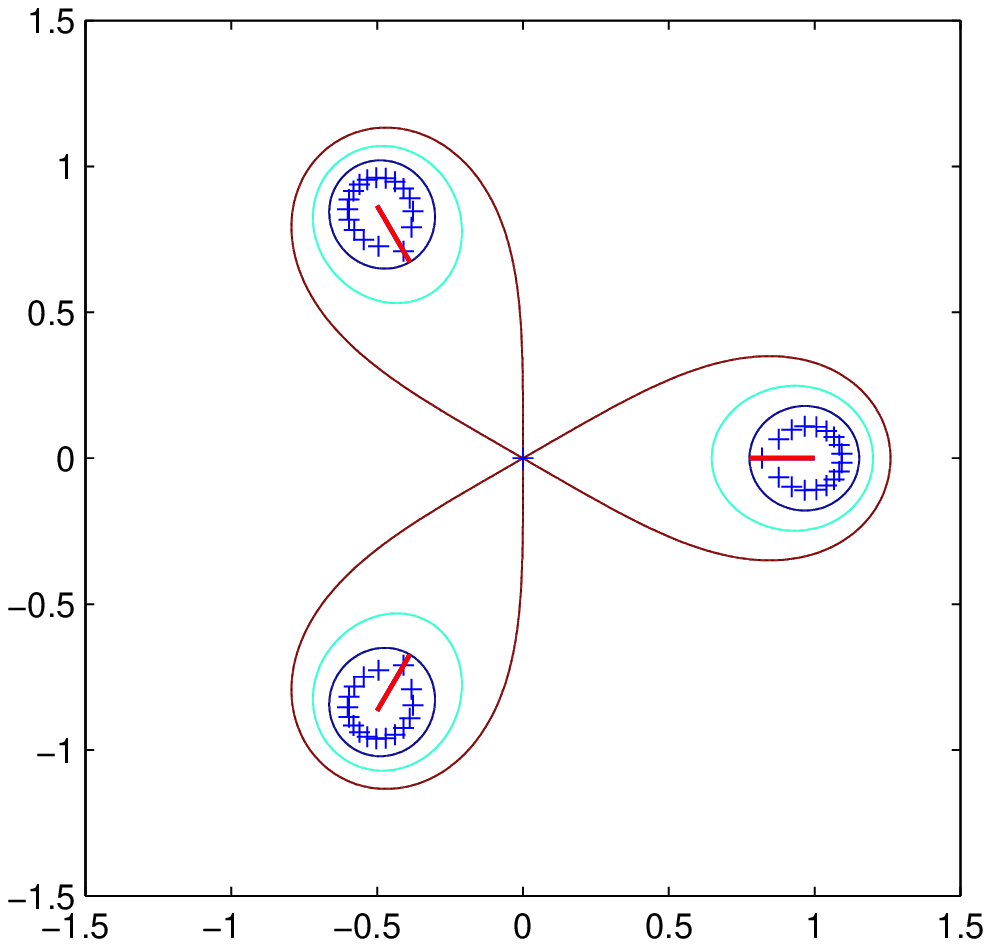}}
\end{center}
\caption{Zeros of the Bergman polynomials $P_n$ for the lemniscate case $m=3$ and $r=0.9$, for $n=50$, $51$ and $52$.}
\label{fig:berolemni}
\end{figure}

In Table~\ref{tab:coeff} we illustrate Proposition~\ref{pro:lemn1} for the lemniscate depicted in Figure~\ref{fig:berolemni}, where $m=3$ and $r=0.9$. More precisely, Table~\ref{tab:coeff} contains the computed values of the leading coefficients $\lambda_n$ correct to 6 decimal figures, for $n=38,\ldots,52$, together with the computed values of $\lambda_nr^{n+1}\sqrt{\frac{\pi}{n+1}}$. As predicted by the theory, the values of $\lambda_nr^{n+1}\sqrt{\frac{\pi}{n+1}}$ alternate, as $n$ increases, towards to the three limits
$$
1/r^2=1.234567 \ldots,\quad 1/r=1.111111\ldots,\quad 1.
$$
The coincidence for the values of $n=38,41,\ldots,50$ is explained in the proof of Proposition~\ref{pro:lemn1}.

\begin{table}[th]
\begin{center}
\begin{tabular}{|c|c|c|}
  \hline
${\ } $
  $n$ & $\lambda_n$ & $\lambda_nr^{n+1}\sqrt{\frac{\pi}{n+1}}$ \\ \hline
 38 &    214.535664 & 1.000000  \\
 39 &    305.078943 & 1.263740   \\
 40 &    305.314216 & 1.124276   \\
 41 &    305.396681 & 1.000000   \\
 42 &    433.231373 & 1.261795   \\
 43 &    433.526043 & 1.123400  \\
 44 &    433.629077 & 1.000000   \\
 45 &    613.834469 & 1.260094  \\
 46 &    614.205506 & 1.122633 \\
 47 &    614.334958 & 1.000000   \\
 48 &    868.011830 & 1.258593  \\
 49 &    868.481244 & 1.121956   \\
 50 &    868.644692 & 1.000000   \\
 51 &   1225.297855 & 1.257261  \\
 52 &   1225.894247 & 1.121355  \\ \hline
\end{tabular}
\end{center}
\medskip\caption{Illustrating Proposition~\ref{pro:lemn1} for the lemniscate case $m=3$ and $r=0.9$, for $n=38,\ldots,52$.}\label{tab:coeff}
\end{table}

\begin{proposition}\label{pro:lemn2}
The following representations hold for the monic
polynomials $p_{km+s}(z)$:
\begin{equation}  \label{eq:proplemn2.a}
p_{km+m-1}(z)=z^{m-1}\,(z^m-1)^k
\end{equation}
and for $s=0,1,\ldots,m-2$, we have for $k$ sufficiently large,
\begin{equation}  \label{eq:proplemn2.b}
\frac{p_{km+s}(z)\left(z^m-1+r^{2m}\right)}{z^s r^{m(k+1)}}=\pi_{k+1,s}(w)-
\frac{\pi_{k+1,s}(-r^m)}{\pi_{k,s}(-r^m)}\pi_{k,s}(w),
\end{equation}
where $w=(z^m-1)/r^m$ and $\pi_{n,s}(w)$ is the monic polynomial of degree $n
$ in $w$ that is orthogonal on the circle $|w|=1$ with respect to the weight
\begin{equation}  \label{eq:proplemn2.c}
\frac{|dw|}{|r^mw+1|^{2-\frac{2}{m}-\frac{2s}{m}}}.
\end{equation}
\end{proposition}

\begin{remark}
The representation formulas (\ref{eq:proplemn2.a})
and (\ref{eq:proplemn2.b}) have the same form as those found by Mi{\~n}a-D{%
\'{\i}}az \cite{MDthe}, who studied the simpler case when $r>1$, i.e.\ when $G$
consists of a single island.
\end{remark}

In our proof we utilize the following lemma that relates 'weighted' Bergman polynomials on the unit disk to Szeg\H{o} polynomials on the unit circle. This result is somewhat implicitly contained in \cite{MDthe}.
\begin{lemma}\label{lem:lemn}
Let $t_n(w)=w^n+\cdots$ be the monic polynomial orthogonal
with respect to the weight $|dw|/|\gamma w+1|^\tau$ on $|w|=1$, where $\tau$ is real, $\tau\neq 2,4,\ldots,2n$, and $|\gamma|<1$. Let $\beta_n(w)=w^n+\cdots$ be the monic polynomial orthogonal with respect to the weight $dA(w)/|\gamma w+1|^\tau$ over the unit disk $|w|<1$. If $t_n(-\overline{\gamma})\ne 0$, then
\begin{equation}  \label{eq:lemlemn}
(w+\overline{\gamma})\beta_n(w)=t_{n+1}(w)- \frac{t_{n+1}(-\overline{\gamma})%
}{t_n(-\overline{\gamma})}\,t_n(w).
\end{equation}
\end{lemma}
Our next result describes the fine asymptotics for the monic
Bergman polynomials.
\begin{proposition}\label{pro:lemn3}
Let
\begin{equation}  \label{eq:proplemn3.a}
\tau:=2-\frac{2}{m}-\frac{2s}{m}, \quad s=0,1,\ldots,m-1.
\end{equation}
Then for $|z^m-1|\geq r^{2m}, z^m-1 \neq -r^{2m}$, the monic Bergman polynomials satisfy for each $%
s=0,1,\ldots,m-1$
\begin{equation}  \label{eq:proplemn3.b}
\lim_{k\to\infty}\frac{p_{km+s}(z)}{z^s (z^m-1)^k}= \left(\frac{z^m-1+r^{2m}%
}{z^m-1}\right)^{\tau/2},
\end{equation}
where the branch of the power function on the right-hand side of (\ref%
{eq:proplemn3.b}) is taken to equal one at infinity, and the convergence is uniform on compact
subsets.\

Furthermore, for each $j=1,2,\ldots,m$ and $z\in G_j$ with $|z^m-1|<r^{2m},$ we have
\begin{equation}  \label{eq:proplemn4.a}
\lim_{k\to\infty}\frac{(-1)^{k+1}k^{2+\tau/2}}{r^{m(2k+4)}}\,p_{km+s}(z)=
\frac{e^{2\pi ij(s+1)/m}\,z^{m-1}\,\tau\Gamma(\tau/2)\,\sin(\tau\pi/2)} {%
2\pi(1-r^{2m})^{\tau/2}(z^m-1+r^{2m})^2}
\end{equation}
for each $s=0,1,\ldots,m-2$,  the convergence being uniform on closed subsets.
\end{proposition}
Observe that the lemniscate $|z^m-1|=r^{2m}$ is the reflection of the lemniscate $|z^m-1|=1$
in the bounding lemniscate of $G$.

\begin{remark}
From the first part of Proposition 7.3 we see that the Bergman polynomials for $G$ have no limit point of zeros in $|z^m-1|>r^{2m}$ other than at $z=0$. Furthermore, from the second part of the proposition, we deduce that, except for the subsequence (\ref{eq:proplemn2.a}), there are no limit points of the zeros of $P_n(z)$ in $|z^m-1|<r^{2m}$. Consequently, the only limit points of zeros of such $P_n(z)$ are at $z=0$ or on the lemniscate $|z^m-1|=r^{2m}$.
\end{remark}
In Figure~\ref{fig:berolemni}, we plot the zeros of the Bergman polynomials $P_n$, for $n=50$, $51$ and $52$, of $G:=\{z:\, |z^3-1|<0.9^3\}$. In each plot, we depict also the defining lemniscate $\Gamma=\{z:\, |z^3-1|=0.9^3\}$, the  reflection $\{z:\, |z^3-1|=0.9^6\}$ of $\{z:\, |z^3-1|=1\}$ in $\Gamma$ and, for the cases $n=51,52$, the branch cuts for the Schwarz function
$\displaystyle{S(z)=(\frac{z^3 -1+ 0.9^6}{z^3-1})^{1/3}}$ of $\Gamma$.

As a consequence of Proposition 7.3 we have the following:
\begin{corollary}
There are precisely two limit measures for the sequence $\{\nu_{P_n}\}_{n=1}^\infty$; namely
$$
\frac{1}{m}\sum_{j=1}^m\delta_{z_j},\quad z_j=\exp(2\pi i j/m),
$$
and the equilibrium measure for the lemniscate $|z^m-1| =r^{2m}$, which is given by the formula
$$
d\beta = \frac{|z|^{m-1}}{r^{2m}}|dz|.
$$
\end{corollary}

\section{Proofs}\label{section:proofs}
\setcounter{equation}{0}
The present section is devoted to the proofs of the results stated earlier in the article.

\medskip\noindent
{\bf Proof of Lemma~\ref{lem:schwarz}.} That $\Gamma$ is analytic is
clear, since $u$ is real analytic and $\nabla u\ne 0$.

All of $D$ is filled with integral curves of the gradient $\nabla
u$. These are disjoint and have no end points in $D$ since $\nabla
u\ne 0$. Hence they all end up on $\partial D$ (an integral curve
cannot be closed since $u$ is single-valued and increases along
it). These integral curves are at the same time level lines of any
locally defined harmonic conjugate of $u$.

Given $z\in D$ we want to define the reflected point
$\overline{S(z)}$ using only $u$. Assume for example that
$u(z)<0$. By the maximum principle, $|u|<c$ in $D$, so actually
$-c<u(z)<0$. There is a unique integral curve $\gamma$ of $\nabla
u$ passing through $z$, and $u$ increases along $\gamma$ with
limiting value $+c$ as $\gamma$ approaches $\partial D$. Thus
there is a unique point $w\in \gamma$ at which $u(w)= -u(z)$. In
terms of this we define
$$
\overline{S(z)}=w.
$$

The above procedure defines a function $S(z)$ in $D$. To see that
$S(z)$ is analytic, note that, in some neighborhood of $\gamma$,
$u$ has a single-valued harmonic conjugate $u^*$ and that $\gamma$
is a level line of $u^*$. The function $f=u+iu^*$ is analytic in a
neighborhood of $\gamma$, with $f'\ne 0$; hence $f$ can be used as
a new complex coordinate near $\gamma$, or $u$ and $u^*$ are new
real coordinates. In terms of these, the reflection map $z\mapsto
\overline{S(z)}$ just defined is given by
$$
u+iu^* \mapsto -u+iu^*,
$$
or $f(z)\mapsto -\overline{f(z)}$. This gives
$$
S(z)=\overline{f^{-1}(-\overline{f(z)})},
$$
which proves that $S(z)$ is analytic. It is also immediate that
$S(z)=\overline{z}$ on $\Gamma$, so that $S$ is indeed a Schwarz
function of $\Gamma$. \qed

\noindent\medskip
{\bf Proof of Lemma~\ref{lem:regular}.}
According to Theorem~3.2.3 of \cite{StTobo}, one criterion for $dA|_G$ to belong to the class \textbf{Reg} is that
\begin{equation}\label{eq:nthnorm}
\lim_{n\to\infty}\|P_n\|_{\Gbar}^{1/n}=1;
\end{equation}
note that $\Omega$ is regular with respect to the Dirichlet problem \cite[p.~92]{Ra}.
(Here and in the sequel $\|\cdot\|$ means the $\sup$ norm on the
subscripted set.)

The argument given in the proof of Lemma~4.3 of \cite{PSG},
when separately applied to each of the
Jordan regions $G_j$ yields
\begin{equation*}
\limsup_{n\to\infty}\|P_n\|_{\overline{G}_j}^{1/n}\le 1,\quad \jeqN.
\end{equation*}
Consequently, $\limsup_{n\to\infty}\|P_n\|_{\Gbar}^{1/n}\le 1$. But $\liminf_{n\to\infty}\|P_n\|_{\Gbar}^{1/n}\ge 1$, since
$\|P_n\|_{L^2(G)}=1$ for all $n$, and so (\ref{eq:nthnorm}) follows.
\qed

\subsection{The extremal problems}
We use $\mathcal{P}_n$ to denote the space of complex polynomials
of degree $n$. Recall that $K_n(z,\zeta)$ denotes the $n$-th finite section of $K(z,\zeta)$
$$
K_n(z,\zeta):=\sum_{k=0}^n \overline{P_k(\zeta)} P_k(z),
$$
and similarly set
$$
K^{G_j}_n(z,\zeta):=\sum_{k=0}^n \overline{P_{k,j}(\zeta)} P_{k,j}(z),
$$
where
$$
P_{n,j}(z) = \lambda_{n,j} z^n+ \cdots, \quad \lambda_{n,j}>0,\quad
n=0,1,2,\ldots,
$$
are the sequences of the Bergman polynomials associated with $G_j$,
$\jeqN$.

\begin{lemma}\label{lem:maxpzet}
For any $\zeta\in \mathbb C$,
$$
\max_{p\in\mathcal{P}_n}\frac{|p(\zeta)|}{\,\,\,\|p\|_{L^2(G)}}=\sqrt{K_n(\zeta,\zeta)},
\quad n=0,1,\ldots.
$$
\end{lemma}
\noindent
{\bf Proof.} Since for any $p\in\mathcal{P}_n$ and $\zeta\in\mathbb{C}$
$$
p(\zeta)=\langle p, K_n(\cdot,\zeta)\rangle,
$$
it follows
$$
|p(\zeta)|\le
\|p\|_{L^2(G)}\,\|K_n(\cdot,\zeta)\|_{L^2(G)}=\|p\|_{L^2(G)}\sqrt{K_n(\zeta,\zeta)}.
$$
Hence
$$
\frac{|p(\zeta)|}{\,\,\,\|p\|_{L^2(G)}}\le\sqrt{K_n(\zeta,\zeta)}
$$
with equality if $p(z)=c\,K_n(z,\zeta)$, for some constant $c\neq 0$.\qed

Obviously
$$
\|p\|_{L^2(G_j)}\le \|p\|_{L^2(G)},\quad \jeqN,
$$
therefore for $n=0,1,\ldots$,
$$
\max_{p\in\mathcal{P}_n}\frac{|p(\zeta)|}{\,\,\,\|p\|_{L^2(G_j)}}\ge
\max_{p\in\mathcal{P}_n}\frac{|p(\zeta)|}{\,\,\,\|p\|_{L^2(G)}},\quad
\jeqN,
$$
or
\begin{equation}\label{eq:KnDjge}
K^{G_j}_n(\zeta,\zeta)\ge K_n(\zeta,\zeta),\quad \jeqN,\quad \zeta\in\mathbb{C}.
\end{equation}
Furthermore, since for any $\zeta\in G_j$,
$$
K^{G_j}_n(\zeta,\zeta)\le K^{G_j}(\zeta,\zeta)=K(\zeta,\zeta),\quad
\jeqN,
$$
it follows from (\ref{eq:KnDjge}) that
\begin{equation}\label{eq:sqKnDjge}
\frac{1}{\sqrt{K_n(\zeta,\zeta)}}\ge
\frac{1}{\sqrt{K_n^{G_j}(\zeta,\zeta)}}\ge\frac{1}{\sqrt{K(\zeta,\zeta)}},
\quad \jeqN.
\end{equation}

\subsection{Proof of Theorem~\ref{th:GPSS_ln}}
The estimates from above require only a $C^{2+\alpha}$-smooth boundary and are based on comparison with corresponding estimates for the arc-length measure $|dz|$ and the Szeg\H{o} orthogonal polynomials. To this purpose, we compare the two extremal problems
\begin{equation}\label{eq:mGdA}
m_n^2(G,dA):=\min_{a_0,\ldots,a_{n-1}}\int_G
|z^n+a_{n-1}z^{n-1}\cdots+a_0|^2 dA(z),\quad n=0,1,2,\ldots,
\end{equation}
and
\begin{equation}\label{eq:mGarhodz}
m_n^2(\Gamma,\rho|dz|):=\min_{a_0,\ldots,a_{n-1}}
\int_\Gamma|z^n+a_{n-1}z^{n-1}\cdots+a_0|^2\rho(z)|dz|,\quad
n=0,1,2,\ldots,
\end{equation}
where $\rho$ is a positive smooth function on $\Gamma$. Recall from (\ref{eq:minimal1}) that
\begin{equation}\label{eq:mGda=}
m_n^2(G,dA)=\frac{1}{\lambda_n^2}=\int_G
|\frac{P_n(z)}{\lambda_n}|^2dA(z),
\end{equation}
where
$$
P_n(z) = \lambda_n z^n+ \cdots, \quad \lambda_n>0,\quad
n=0,1,2,\ldots,
$$
are the Bergman polynomials of $G$.

The asymptotic properties of $m_n(\Gamma,\rho |dz|)$ have been established by Widom in \cite[Thm~9.1]{Wi69}. In particular, the next estimate for $\rho=1$ and some constant $C>0$ follows from Theorems~9.1 and 9.2 of \cite{Wi69}:
\begin{equation}\label{eq:widom}
m_n^2(\Gamma,|dz|)\geq C \,{\rm cap}(\Gamma)^{2n}.
\end{equation}
On the other hand, Suetin's lemma (Lemma~\ref{lem:suetin} above) applied to each island separately gives
$$
m_n(G,dA)^2 = \int_G | \frac{P_n(z)}{\lambda_n} |^2 dA \geq
\frac{C}{{n+1}} \int_\Gamma | \frac{P_n(z)}{\lambda_n} |^2 |dz|
\geq \frac{C}{{n+1}}\, m_n(\Gamma, |dz|)^2,
$$
where $C>0$ is a another positive constant.

Combining the above two estimates we conclude
$$
m_n(G,dA)\geq C\,\frac{{\rm cap\,}(\Gamma)^n}{\sqrt{n}},
$$
which yields the upper inequality in Theorem~\ref{th:GPSS_ln}.

For estimates from below we require analyticity of the boundary. The main technical aid is provided by a family of polynomials $\omega_n$ constructed by Walsh in \cite{Wa63}, which we thereby refer to as \textit{Walsh polynomials}.

\begin{lemma}\label{lem:WaSa}
Assume that each $\Gamma_j$, $\jeqN$, is analytic. Then, there
exists a sequence of monic polynomials $\omega_n(z)=z^{n}+\cdots$,
$n=1,2,\dots$, with all zeros on a fixed compact subset $E\subset
G$, and a constant $C$ such that
\begin{equation}\label{eq:ome_nle}
\|\omega_{n}\|_{L^2(G)} \leq \frac{C}{\sqrt{n}}\ \capGm^n.
\end{equation}
\end{lemma}
From this we deduce the lower inequality in Theorem~\ref{th:GPSS_ln}:
\begin{corollary}
If each $\Gamma_j$, $\jeqN$, is analytic then
\begin{equation}\label{eq:lamnge}
C\frac{\sqrt{n}}{\capGm^{n}}\le \lambda_n .
\end{equation}
\end{corollary}

\noindent{\bf Proof of Lemma~\ref{lem:WaSa}.}
Since each $\Gamma_j$, $\jeqN$, is analytic, the Green function $g_\Omega(z,\infty)$ extends harmonically across $\partial G$ by Schwarz reflection. Choose first a number $0<\tau<1$ such that $\frac{1}{\tau}<R'$ (see Subsection~\ref{subsec:green} for
the definition of $R'$) and such that $g_\Omega(z,\infty)$ extends into each component of $G$, at least to the negative level $\log\tau$. Since $g_\Omega (z,\infty)$ has no critical points in
$\mathcal{G}_{R'}\setminus G$ it follows that the extended Green function has no critical points in
$D=\mathcal{G}_{\frac{1}{\tau}}\setminus\overline{\mathcal{G}}_{\tau}=
\Omega_\tau\setminus\overline{\Omega}_{1/\tau}$. The latter open set has $N$ components, each of which is a domain of involution for the Schwarz reflection (see Lemma~\ref{lem:schwarz}).

Now choose a number $\rho$ in the interval
$$
\tau<\rho< 1.
$$
For any $R\geq\rho$,
\begin{equation}\label{eq:gomeR}
g_{\Omega_R}(z,\infty):=g_{\Omega}(z,\infty)-\log R,
\end{equation}
is the Green function of $\Omega_R$ with pole at infinity. Hence,
\begin{equation}\label{eq:capGR}
\mathrm{cap} (L_R)=R\, \mathrm{cap} (\Gamma).
\end{equation}
Choose the compact set $E\subset G$ in the statement of the lemma to
be $E=L_\tau$. By a  theorem of Walsh \cite{Wa63} (see also \cite[p.~515]{Sa69}), there exists a
sequence of monic polynomials $\omega_n(z)=z^{n}+\cdots$,
$n=1,2,\dots$, with zeros approximately equidistributed with respect
to the conjugate function of $g_\Omega(z,\infty)$ and such that
\begin{equation}\label{eq:walsh}
|g_{\Omega_{\tau}}(z,\infty)+\log{\mathrm{cap\,}}
(L_\tau)-\frac{1}{n}\log|\omega_n(z) || \leq \frac{C}{n} \quad {\rm
in}\quad \Omega_\rho.
\end{equation}

Note that $g_{\Omega_{R}}(z,\infty)+\log{\mathrm{cap\,}}(L_R)$ is
independent of $R$, hence in (\ref{eq:walsh}) $\tau$ can be replaced
by  any number $R>\tau$. For $z\in L_R$ and $R\geq \rho$ this gives
$$
|\log{\mathrm{cap}}(L_R)-\frac{1}{n}\log|\omega_n(z) || \leq
\frac{C}{n},
$$
or after exponentiating and using (\ref{eq:capGR})
\begin{equation}\label{omega-bounds}
e^{-C}\leq \frac{|\omega_n(z)|}{{R^n\mathrm{cap}}(\Gamma)^n}\leq e^{C}
\quad(z\in L_R,\, \rho\leq R< \infty).
\end{equation}
In particular, from the maximum principle,
\begin{equation}\label{eq:omegan}
|\omega_n (z)|\leq C\, R^n \capGm^n \quad (z\in \mathcal{G}_R,\,
\rho\leq R< \infty),
\end{equation}
for another constant $C$.

Next we estimate the $L^2(G)$-norm of $\omega_n$. On decomposing
$$
\int_G|\omega_n|^2 \,dA =\int_{G_\rho}|\omega_n|^2 \,dA
+\int_{G\setminus  G_\rho}|\omega_n|^2 \,dA,
$$
the first term can be directly estimated by means of (\ref{eq:omegan}):
$$
\int_{G_\rho}|\omega_n|^2 \,dA\leq C \max_{z\in L_\rho}|\omega_n(z)|^2
\leq C \rho^{2n}\capGm^{2n}.
$$

For the second term we foliate $G\setminus G_\rho$ by the level
lines $L_R$ of $g_\Omega (z,\infty)$, or $|\Phi(z)|=\exp[g_\Omega
(z,\infty)]$, and use the coarea formula. Since $\nabla g_\Omega
(z,\infty)$, and hence $\nabla |\Phi(z)|$, is bounded away from zero
on $G\setminus G_\rho$ we obtain by using once more (\ref{eq:omegan})
\begin{eqnarray*}
\int_{G\setminus  G_\rho}|\omega_n|^2 \,dA &=& \int_\rho^1\int_{L_R}
              \frac{|\omega_n(z)|^2}{|\nabla |\Phi{(z)}||} \,|dz|dR\\
&\leq& C\,\int_\rho^1 \max_{z\in L_R}|\omega_n(z)|^2\, dR \leq C\,\capGm^{2n}
\int_\rho^1 R^{2n}\,dR \\
&\leq& C\, \capGm^{2n}\,\frac{1-\rho^{2n+1}}{2n+1}\leq C\frac{\capGm^{2n}}{n},
\end{eqnarray*}
for various positive constants $C$. Thus altogether we have
$$
\int_G|\omega_n|^2 \,dA\leq C(\rho^{2n}+\frac{1}{n})\capGm^{2n},
$$
and since $\rho<1$, this gives (\ref{eq:ome_nle}). \qed

The corollary is an immediate consequence of the lemma and the
definition of $\lambda_n$:
$$
\frac{1}{\lambda_n}=m_n (G,dA) \leq \| \omega_n \|_{L^2(G)} \le C
\frac{\sqrt{n}}{\capGm^n}.
$$

\subsection{Proof of Theorem~\ref{pro:ChrinC}}
We turn now our attention to the problem of determining the rate of convergence of
$\Lambda^{G_j}_n$ as compared to $\Lambda_n$. The solution will obviously depend on a set of numerical constants which reflect the global configuration of $G$.

In the case of a single island $N=1$ we have $\Lambda^{G_1}_n\equiv \Lambda_n$, hence both  (\ref{eq:ChrinC0}) and (\ref{eq:ChrinC}) hold trivially with $m=1$. For the case $N\ge 2$, we assume that $\Gamma_j$ is analytic, for some fixed $j\in\{1,2,\ldots,N\}$. Let $\mathcal{X}$ denote the characteristic function of $\overline{G}_j$ in $\overline{G}$ and set
\begin{equation}\label{eq:gamma_n}
\gamma_n:=\inf_{p\in\mathcal{P}_n}\frac{\|\mathcal{X}p\|_{L^2(G)}}{\|p\|_{L^2(G)}}.
\end{equation}
(Note that $\|\mathcal{X}p\|_{L^2(G)}=\|p\|_{L^2(G_j)}$, hence $0<\gamma_n<1$.)

By considering the Bergman polynomial $P_{n,j}$ of $G_j$, as a competing polynomial in (\ref{eq:gamma_n}) and using Carleman asymptotics (Theorem~\ref{th:Carleman}) for $P_{n,j}$ in $G\setminus G_j$ in conjunction with the fact $|\Phi_j(z)|>|\Phi(z)|$, $z\in\Omega$ (subordinate principle for the Green function; see e.g. \cite[p.~108]{Ra}), we conclude that there exist constants $C>0$ and $R>R_j\ (>1)$ such that, for any $n\in\mathbb{N}$,
$$
\frac{1}{\gamma_n}\ge 1+C\ \sqrt{n}\ R^n.
$$
Hence for large values of $n$,
$$
\gamma_n < \alpha^n,
$$
where $0<\alpha<1$. Since $\mathcal{X}$ has an analytic continuation up to $L_{R'}$ in $\Omega$, it follows
from Walsh's theorem of maximal convergence \cite[Thm IV.5]{Wa} that for any $n\in\mathbb{N}$, there exist a constant $m \geq 1$ and a polynomial $q_{m(n)}\in\mathcal{P}_{m(n)}$, where $m(n)=mn$, with the property,
\begin{equation}\label{eq:qmn-X}
\|q_{m(n)}-\mathcal{X}\|_{\overline{G}}<\gamma_n.
\end{equation}
Then we have:
\begin{lemma}
Assume that $\Gamma_j$, $j\in\{1,2,\dots, N\}$, is analytic. Then for any $n\in\mathbb{N}$
$$
\sqrt{K_n^{G_j}(\zeta,\zeta)}\le\frac{2}{1-\gamma_n}\sqrt{K_{n+m(n)}(\zeta,\zeta)},\quad \zeta\in \overline{G}_j.
$$
\end{lemma}
\noindent
{\bf Proof.} Take $\zeta\in \overline{G}_j$ and let $h\in\mathcal{P}_n$ be an extremal polynomial for
$$
\max_{p\in\mathcal{P}_n}\frac{|p(\zeta)|}{\,\,\,\|p\|_{L^2(G_j)}}.
$$
Then from Lemma~\ref{lem:maxpzet}
$$
\sqrt{K_n^{G_j}(\zeta,\zeta)}=\frac{|(\mathcal{X}h)(\zeta)|}{\|\mathcal{X}h\|_{L^2(G)}}.
$$
It holds,
$$
|(\mathcal{X}h)(\zeta)|\le
\frac{1}{1-\gamma_n}\,|(q_{m(n)}h)(\zeta)|,
$$
because from (\ref{eq:qmn-X}),
$$
(1-\gamma_n)\mathcal{X}(\zeta)\le |q_{m(n)}(\zeta)|.
$$
Also
\begin{eqnarray*}
\|q_{m(n)}h\|_{L^2(G)}&\le&\|\mathcal{X}h\|_{L^2(G)}+\|(\mathcal{X}-q_{m(n)})h\|_{L^2(G)}\\
  &\le& \|\mathcal{X}h\|_{L^2(G)}+\gamma_n\|h\|_{L^2(G)}\
  \le\ 2 \|\mathcal{X}h\|_{L^2(G)},
\end{eqnarray*}
where in the last inequality we made use of the defining property of $\gamma_n$.
Finally,
\begin{eqnarray*}
 \frac{|(\mathcal{X}h)(\zeta)|}{\|\mathcal{X}h\|_{L^2(G)}}&\le&
\frac{2}{1-\gamma_n}\, \frac{|(q_{m(n)}h)(\zeta)|}{\|q_{m(n)}h\|_{L^2(G)}} \\
  &\le&\frac{2}{1-\gamma_n}\,
  \max_{f\in\mathcal{P}_{n+m(n)}}\frac{|f(\zeta)|}{\,\,\,\|f\|_{L^2(G)}},
\end{eqnarray*}
and the result follows from Lemma~\ref{lem:maxpzet}. \qed

This yields Inequality~(\ref{eq:ChrinC}) in Theorem~\ref{pro:ChrinC}. The other inequality (\ref{eq:ChrinC0}) follows immediately from (\ref{eq:KnDjge}).

\subsection{Proof of Theorem~\ref{th:ChrinGj}}
Keeping in mind Lemma~\ref{lem:dense}, it is clear from its definition that the functions $\Lambda_n$ converge uniformly on compact subsets of $G$ to $\Lambda$. By imposing analyticity of the boundary, we will be able to estimate jointly the rate of convergence of $\Lambda_n(z)$ on $\Gamma$ and in a neighborhood of $\Gamma$ in the interior. In view of the reduction to a single island established in the previous subsection, we will assume in the first part of the proof that $N=1$. In order to simplify further the notation, we will simply write $G=G_1$, $\Phi=\Phi_1$ and so forth.

Thus, we deal now with a Jordan domain $G$ with analytic boundary $\Gamma$. The normalized external conformal mapping $\Phi$ analytically extends to the level set $\mathcal G_\rho$, with $\rho<1$. According to Theorem~\ref{th:Carleman}, the Bergman orthogonal polynomials satisfy:
$$
P_n(z) = \sqrt{\frac{n+1}{\pi}} \Phi(z)^n \Phi'(z)\{1+A_n(z)\}, \ \
z \in G \setminus \overline{\mathcal G_\rho},
$$
where $A_n(z) = O((\frac{\rho}{r})^n),$ whenever $z \in \Gamma_r$, and $\rho<r<1.$
Fix a $z \in G \setminus \overline{\mathcal G_\rho}$ and denote $t = |\Phi(z)|^2$. Then
\begin{eqnarray}\label{eq:KnRn}
K_n(z,z)&=&\sum_{k=0}^n |P_k(z)|^2=\frac{|\Phi'(z)|^2}{\pi} \sum_{k=0}^n(k+1)t^k+R_n(z)  \\
   &=&\frac{|\Phi'(z)|^2}{\pi}\frac{1-(n+2)t^{n+1}+(n+1)t^{n+2}}{(1-t)^2}+R_n(z).\nonumber
\end{eqnarray}
Similarly,
$$ K(z,z) = \frac{|\Phi'(z)|^2}{\pi} \frac{1}{(1-t)^2} + R(z).$$
The convergence of $R_n(z)$ to $R(z)$, for $\rho^2<r^2\leq t<1$,
is uniformly dominated by a convergent geometric series.

In view of (\ref{eq:LamdazGamma}) we set $\Lambda(z)=0$ for all $z\in\Gamma$. Since
$$
0<\Lambda_n(z -\Lambda(z) =\frac{1}{\sqrt{K_n(z,z)}} - \frac{1}{\sqrt{K(z,z)}},
$$
we are led to the estimate
$$
\Lambda_n(z)-\Lambda(z)\leq C(1-t)[ \frac{1}{\sqrt{1-(n+2)t^{n+1}+(n+1)t^{n+2}}}-1].
$$
In its turn, elementary algebra yields:
\begin{eqnarray*}
(1-t) [\frac{1}{\sqrt{1-(n+2)t^{n+1}+(n+1)t^{n+2}}} -1] &=&
\frac{1}{\sqrt{\sum_{k=0}^n(k+1)t^n}}
\frac{(n+2)t^{n+1}-(n+1)t^{n+2}}{1+\sqrt{1-(n+2)t^{n+1}+(n+1)t^{n+2}}}\\
  &\le&\frac{n+1}{t^{n/2} \sqrt{1+2+...+(n+1)}} t^{n+1} [1-\frac{1}{n+1} -t]  \\
  &\le& C t^{n/2}(1-t + \frac{1}{n}),
\end{eqnarray*}
which implies Inequality (\ref{eq:ChrinGj}) in Theorem \ref{th:ChrinGj}, since for $z$ near $\Gamma$:
$$
1-|\Phi(z)|^2\asymp 1-|\Phi(z)|\asymp \textup{dist}(z,\Gamma).
$$

Using (\ref{eq:KnRn}), which holds for $z\in\Gamma$ with $R_n(z)=O(n^2\sqrt{n}\rho^n)$, we derive easily (\ref{eq:limnChr}), which is the limit of the exact form of (\ref{eq:ChrinGj}).

We resume now our general assumption $G=\cup_{j=1}^N G_j$ and we turn our attention to deriving (\ref{eq:CrhinGm}). The lower bound emerges at once by combining (\ref{eq:limnChr}) with (\ref{eq:ChrinC0}). To obtain the upper bound we apply (\ref{eq:ChrinC}) to $\Lambda_k(z)$, for large $k$, with $k=[k/m]m+r$, where $0\le r<m-1$, and $[k/m]$ is the integral part of the fraction, and then we use again (\ref{eq:limnChr}).

In order to estimate $\Lambda_n$ in the exterior of $\overline{G}$ we employ the Walsh polynomials: From Lemma~\ref{lem:maxpzet},
$$
\Lambda_n(z)=\min_{p\in\mathcal{P}_n}\frac{\,\|p\,\|_{L^2(G)}}{|p(z)|}
$$
and therefore,
$$
\Lambda_n(z)\le\frac{\,\|\omega_n\|_{L^2(G)}}{|\omega_n(z)|}
\leq C \frac{1}{\sqrt{n}|\Phi(z)|^n};
$$
where we made use of Lemma~\ref{lem:WaSa} and (\ref{omega-bounds}).

Finally, the lower estimate for $\Lambda_n(z)$ for $z$ exterior to $\overline{G}$ is directly derived from the upper estimates for the orthogonal polynomials appearing in Theorem~\ref{th:PnOme}. \qed

\subsection{Proof of Theorem~\ref{th:PnOme}}
Our aim is to derive estimates for $P_n(z)$, for $z$ in the exterior of the archipelago. To do so, we assume that every curve constituting  $\Gamma$ is analytic and we rely, once more, to the Walsh polynomials $\omega_n$.

We fix a positive integer $n$ and consider the rational function $\frac{P_n(z)}{\omega_{n+1}(z)}$, whose poles lie in a compact subset of $G$ and which vanishes at infinity. With $z\notin \overline{G}$, Cauchy's formula yields:
$$
\frac{P_n(z)}{\omega_{n+1}(z)}=\frac{-1}{2\pi i}\int_\Gamma
\frac{P_n(\zeta)d\zeta}{\omega_{n+1}(\zeta)(\zeta-z)},
$$
whence, from (\ref{omega-bounds}),
$$
|P_n(z)|\leq\frac{C}{{{\rm dist}(z,\Gamma)}}\,\frac{|\omega_{n+1}(z)|}{\capGm^{n+1}}\,\,\|P_n\|_{L^1(\Gamma)},
$$
where $\|\cdot\|_{L^1(\Gamma)}$ denotes the $L^1$-norm on $\Gamma$ with respect to $|dz|$.

Since the $L^1$-norm is dominated by a constant times the $L^2$-norm, Lemma~\ref{lem:suetin} gives $\|P_n\|_{L^1(\Gamma)}\le C\,\sqrt{n}$ and one more application of (\ref{omega-bounds}) yields
$$
|P_n(z)| \leq \frac{C}{{{\rm dist}(z,\Gamma)}} \sqrt{n} |\Phi(z)|^n.
$$
(In the above we use $C$ to denote positive constants, not necessarily the same in all instances.)

In order to obtain the estimates from below, we have to restrict the point $z$ to the complement of the convex hull ${\rm Co} (\overline{G})$. On that set, including the point at infinity, the sequence of rational functions $R_n(z)=\frac{P_n(z) {\rm cap}(\Gamma)^n}{\sqrt{n}\omega_{n}(z)}$ has no zeros, and by the above estimate, it is equicontinuous on compact subsets of $U= \overline{\mathbb C}\setminus{\rm Co}(\overline{G})$. Thus $\{R_n\}_{n=0}^\infty$ forms a normal family on $U$ and the possible limit functions are either identically zero, or zero free. The normalization at infinity was chosen so that, in view of (\ref{eq:GPSS_lnge}) and (\ref{omega-bounds}), $\inf_{n\in\mathbb{N}} R_n(\infty) >0$. Thus, every limit point of the sequence $R_n$ is bounded away from zero, on compact subsets of $U$.

\subsection{Distribution of Zeros}
{ $\,$}

\medskip
\noindent
{\bf Proof of Theorem~\ref{thm:ASgeneral}.} To prove (i) we need to figure out the general
structure of $\rho(K(\cdot,z))$. We have already remarked, cf.\ (\ref{cor2eq}), that for $\zeta\in G_j$,
$$
\rho (K(\cdot,\zeta))=\min\{R_j, \rho(K^{G_j}(\cdot,\zeta))\}.
$$
Recall (\ref{eq:Kjphij}), that is, in terms of any conformal mapping $\varphi_j:G_j\to \mathbb{D}$,
$$
K^{G_j}(z,\zeta)=\frac{\varphi'_j(z)\conj{\varphi'_j(\zeta)}}
{\pi\,\left[1-\varphi_j(z)\conj{\varphi_j(\zeta)}\right]^2},\quad z,\zeta\in G_j.
$$
Conversely, if (given $\zeta\in G_j$) $\varphi_j$ is chosen so that
$\varphi_j (\zeta)=0$, then
$$
\varphi_j(z)=\frac{\pi}{\overline{\varphi'_j(\zeta)}}\int_\zeta^z
K^{G_j} (t,\zeta) dt.
$$
Hence, for a general $\varphi_j$,
$$
\frac{\varphi_j (z)-\varphi_j (\zeta)}{1-\varphi_j
(z)\overline{\varphi_j (\zeta)}}=
\frac{\pi(1-|\varphi_j(\zeta)|^2)}{\overline{\varphi_j'(\zeta)}}
\int_\zeta^z K^{G_j} (t,\zeta) dt.
$$
It follows therefore that, given a $\zeta\in G_j$ and a simply connected region $D$ with $G_j\subset D\subset \mathcal{G}_{j,R_j}$, $K^{G_j}(z,\zeta)$ has an analytic extension to $D$ as a function of $z$ if and only if $\varphi_j(z)$ has a meromorphic extension to $D$ and does not attain the value
$1/\overline{\varphi_j(\zeta)}$ there.

We introduce a meromorphic version of the function $\rho$ defined in (\ref{rho}) by setting, for $f$ meromorphic in $G$,
\begin{equation}\label{rhom}
\rho_m (f):=\sup\left\{R\geq 1:f\ \mathrm{\ has \ a\ meromorphic\
continuation\ to}\ \mathcal{G}_R\right\}.
\end{equation}
Next we extend each $\varphi_j$ to all $G$ by setting $\varphi_j=0$ in $G\setminus G_j$. Clearly the so extended $\varphi_j$ cannot be meromorphic in $\mathcal{G}_{j,R}$ for any $R>R_j$, hence
\begin{equation}\label{eq:rhoRj}
1\le\rho_m(\varphi_j)\leq R_j.
\end{equation}
(This is vacuous statement if $N=1$, thus we simply set $R_1=+\infty$ in such a case.) The largest $R$ for which $\varphi_j$ does not take the value ${1}/{\overline{\varphi_j(\zeta)}}$ in
$\mathcal{G}_{j,R}$ is
$
\displaystyle{\inf \{|\Phi (\varphi_j|_{\mathcal{G}_{j,\rho_m(\varphi_j)}}^{-1}
({1}/{\overline{\varphi_j (\zeta)}}))| \}\ (\ge 1)},
$
where the infinmum is taken over all points in the preimage
$\varphi_j|_{\mathcal{G}_{j,\rho_m(\varphi_j)}}^{-1}({1}/{\overline{\varphi_j
(\zeta)}})$, which is a subset of $\mathcal{G}_{j,\rho_m(\varphi_j)}\setminus G_j$. (We assign the value $+\infty$ for the infimum of the empty set.)

Putting things together we get, in view of (\ref{eq:rhoRj}),
\begin{equation}\label{eq:rhoPhi}
\rho(K(\cdot,\zeta))=\min\left\{\rho_m (\varphi_j),\ \inf\{ |\Phi
(\varphi_j|_{\mathcal{G}_{j,\rho_m(\varphi_j)}}^{-1}({1}/{\overline{\varphi_j
(\zeta)}}))|\} \right\},\ \zeta\in G_j,
\end{equation}
or, by taking the logarithm,
\begin{equation}\label{eq:logrho}
\log \rho(K(\cdot,\zeta))=\min\left\{\log \rho_m(\varphi_j),\ \inf\{
g_\Omega (\varphi_j|_{\mathcal{G}_{j,\rho_m(\varphi_j)}}^{-1}
({1}/{\overline{\varphi_j (\zeta)}}),\infty)\} \right\},\ \zeta\in G_j.
\end{equation}
This may look messy, but in principle it means that we have expressed $\log \rho(K(\cdot,\zeta))$ as the infimum of some harmonic functions. This is the basic argument telling that $\log\rho(K(\cdot,\zeta))$ is superharmonic as a function of $\zeta$ in $G_j$.

Now, if $\varphi_j$ has a singularity on $\Gamma_j$, then $\rho_m(\varphi_j)=1$ and $\rho(K(\cdot,\zeta))=1$, $\zeta\in G_j$. In the complementary case, i.e., if $\varphi_j$ has an analytic continuation across $\Gamma_j$, then for any $\zeta\in G_j$,   $\varphi_j|_{\mathcal{G}_{j,\rho_m(\varphi_j)}}^{-1}({1}/{\overline{\varphi_j (\zeta)}})$ is either void or it defines a (possibly) multi-valued reflection map in $\Gamma_j$, i.e., the conjugate of a (possibly) multi-valued Schwarz function of $\Gamma_j$. By our assumption that the infimum of the empty set is $+\infty$, we only need to concentrate on the latter case. Denoting $\varphi_j|_{\mathcal{G}_{j,\rho_m(\varphi_j)}}^{-1}({1}/{\overline{\varphi_j (\zeta)}})$ by $S_{j,{\rm multi}}(\zeta)$ we can write (\ref{eq:logrho}) somewhat more handily as
\begin{equation}\label{eq:logrhoS}
\log \rho(K(\cdot,\zeta))=\min\left\{\log \rho_m(\varphi_j),\ \inf\{
g_\Omega (\overline{S_{j,{\rm multi}}(\zeta)},\infty) \}\right\},\ \zeta\in G_j,
\end{equation}
where the infinmum is taken over all branches of $S_{j,{\rm multi}}(\zeta)$. One step further, this reflection map gives a multi-valued analytic extension of the Walsh function $\Phi$ into $G_j$:
$$
\Hat{\Phi}_{\rm multi}(\zeta)=1\big/\conj{\Phi\left(\overline{S_{j,{\rm multi}}(\zeta)}\right)},
\quad \zeta\in G_j
$$
(where we have used hat to emphasize the analytic extension). Inserting the latter into (\ref{eq:rhoPhi}) gives the following, more direct, description of $\rho(K(\cdot,\zeta))$:
\begin{equation}\label{eq:rhoPhi2}
\rho(K(\cdot,\zeta))=\min\left\{\rho_m (\varphi_j),\ \inf\{1/|\Hat{\Phi}_{\rm multi}(\zeta)|\}\right\},\ \zeta\in G_j,
\end{equation}
the infimum is taken, again, over all (local) branches.

In order to make the above considerations more rigorous we take (\ref{eq:logrho}) as our starting point. We first treat the case $N\geq 2$, which is somewhat simpler because in this case (\ref{eq:rhoRj}) gives an upper bound for $\log\rho(K(\cdot,\zeta))$ in (\ref{eq:logrho}). Let $\zeta\in G_j$. Then ${1}/{\overline{\varphi_j (\zeta)}}$ is outside the closed unit disk, and the preimage
$\varphi_j|_{\mathcal{G}_{j,\rho_m(\varphi_j)}}^{-1}({1}/{\overline{\varphi_j (\zeta)}})$ is either empty or is a finite or infinite subset of $\mathcal{G}_{j,\rho_m(\varphi_j)}\setminus
G_j$. If it is an infinite set, then all cluster points will be on the boundary of $\mathcal{G}_{j,\rho_m(\varphi_j)}$, where $g_\Omega(\cdot, \infty)$ is larger, than near $\Gamma_j$. This means that only finitely many of the points in the preimage will be serious candidates in the competition for the infimum in (\ref{eq:logrho}). We may also vary $\zeta$ within a small disk, compactly contained in $G_j$, and there will still be only finitely many branches of the multivaled analytic function $\varphi_j|^{-1}$ involved, when forming the infimum. Within such a disk there will also be only finitely many branch points (where two or more preimages coincide).

Thus, locally away from the mentioned branch points, $\log\rho(K(\cdot,\zeta))$ is the infimum of finitely many harmonic functions, hence is continuous and superharmonic. At the branch points
$\log \rho(K(\cdot,\zeta))$ is still continuous, and since the set of branch points is discrete (in $\mathcal{G}_{j,\rho_m(\varphi_j)}\setminus G_j$) they make up a removable set for continuous superharmonic functions; see e.g.\ \cite[Thm~3.6.1]{Ra}. It follows, therefore, that $\log\rho(K(\cdot,\zeta))$ is superharmonic (and continuous) in all $G_j$.

We apply now the above inferences to $h(z)=-\log\rho(K(\cdot,z))$, for $z\in G$. If $\rho_m (\varphi_j)=1$, for some $j$, then $h(z)=0$, for $z\in G_j$, hence the transition across $\Gamma_j$ to $g_\Omega(z,\infty)$ is continuous and subharmonic. If $\rho_m(\varphi_j)>1$ and $\varphi_j$ remains univalent in a neighborhood of $\overline{G}_j$, then it is easy to see that $h(z)$ defines the harmonic continuation of $g_\Omega (z,\infty)$ across $\Gamma_j$ (in fact, $\Gamma_j$ turns out to be analytic and thus $S_{j,{\rm multi}}$ is the associated ordinary single-valued Schwarz function). Finally, if $\rho_m (\varphi_j)>1$ but $\varphi_j$ is not univalent in any neighborhood of $\overline{G}_j$ then locally, away from finitely many branch points on $\Gamma_j$, $h$ is still the ordinary harmonic continuation of $g_\Omega (z,\infty)$. At the branch points $h$ is still continuous and the set of branch points is too small to affect the overall subharmonicity. Hence, in all possible situations $h(z)=-\log\rho(K(\cdot,z))$ is continuous and subharmonic in $G$.

Therefore, we have established so far that in the case $N\geq 2$, $h$ is subharmonic (and continuous) in $\mathbb{C}$ and since it coincides with the Green function in $\Omega$, $\beta$ is a positive measure, with support contained in $\overline{G}$. Moreover, from Gauss' theorem (see e.g. \cite[p.~83]{ST}), and the singularity of the Green function at infinity, we have for any $R>1$:
\begin{equation}\label{eq:betaGR}
\beta(G_R)=\frac{1}{2\pi }\int_{L_R}\frac{\partial h}{\partial n}\,ds
=\frac{1}{2\pi }\int_{L_R}\frac{\partial g_\Omega(z,\infty)}{\partial n}\,ds=1.
\end{equation}
Hence $\beta$ is a unit measure and this completes the proof of (i), for $N\geq 2$.

In order to derive (ii), we observe that the Riesz decomposition theorem for subharmonic functions applied to $h$ in $\mathbb{C}$ (see e.g. \cite[p.~76]{Ra}) gives,
$$
h(z)=-U^\beta(z)+v(z),  \,\, z\in\mathbb{C},
$$
where $v$ is harmonic in $\mathbb{C}$. Then, by considering the expansions near infinity of $U^\beta(z)$ and  $h(z)=g_\Omega(z,\infty)$, we see that $v(z)=-\log{\rm cap\,}(\Gamma)$, which yields (\ref{eq:potbeta}). Relation (\ref{eq:bal}) is an immediate consequence of (\ref{eq:potbeta}) the fact that $h$ coincides with the Green function in $\Omega$, in conjunction with the relations (\ref{UmuinK})--(\ref{UmuonE}).

When $N\geq 2$,  $U^\beta$ is bounded from above because of
(\ref{eq:rhoRj}):
$$
U^\beta  \leq \log\frac{\max_j\{R_j\}}{{\rm cap\,}(\Gamma)}<\infty.
$$

Statement (iii) of the theorem is just a juxtaposition of  Proposition~\ref{pro:nthroot} and Corollary~\ref{cor2} along with (\ref{eq:potbeta}).

As for (iv), $\mathcal{C}$ is nonempty by general compactness principles for measures and the known fact that all counting measures $\nu_{P_n}$ have support within a fixed compact set; see Remark~\ref{rem:Fejer}. Let $\sigma\in\mathcal{C}$. Then there is a subsequence $\mathcal{N}=\mathcal{N}_\sigma\subset \mathbb{N}$ such that
\begin{equation}
\nu_{P_n}\sta\sigma, \quad n\to\infty,\ n\in\mathcal{N}.
\end{equation}
Using the lower envelope theorem \cite[Thms~I.6.9 ]{ST} and (\ref{eq:liminf}) we get
\begin{equation}\label{eq:starsigma2}
U^\sigma(z)=\liminf_{\underset{n\in\mathcal{N}}{n\to\infty}} U^{\nu_{P_n}}(z)
\ge\liminf_{n\to\infty} U^{\nu_{P_n}}(z)=U^\beta(z),
\end{equation}
where the first equality holds only quasi everywhere in $\mathbb{C}$. However the relation between $U^\sigma$ and $U^\beta$ persists everywhere in $\mathbb{C}$, since both members are potentials.

Let $\mathcal{D}$ be any component of $\mathbb{C}\setminus{\rm supp\,}\beta$. Applying the minimum principle to $u=U^\sigma-U^\beta\geq 0$, which is superharmonic in $\mathcal{D}$, gives that either $u>0$ in all $\mathcal{D}$ or $u=0$ in all $\mathcal{D}$. Since $u$ vanishes at $\infty$ (recall that $\sigma$ and $\beta$ are unit measures) it follows that it vanishes in the entire unbounded component of $\mathbb{C}\setminus{\rm supp\,}\beta$. From this and the observations above follow all parts of (iv).

Turning to (v), let
$$
U={\rm lsc\,}{(\inf_{\sigma\in\mathcal{C}} U^\sigma)}.
$$
By (iv), $U^\beta\leq U$ in $\mathbb{C}$. To prove the opposite inequality, choose an arbitrary point $z\in\mathbb{C}$. Then there is subsequence $\mathcal{N}_z\subset\mathbb{N}$, such that the $\liminf$ in (\ref{eq:liminf}) is realized at $z$, i.e.
\begin{equation}\label{eq:liminflim}
\lim_{\underset{n\in\mathcal{N}_z}{n\to\infty}}U^{\nu_{P_n}}(z)= U^\beta (z).
\end{equation}
By weak* compactness there exists a further subsequence $\mathcal{N}'_z\subset\mathcal{N}_z$ and a measure $\sigma=\sigma_z\in\mathcal{C}$ such that
\begin{equation}
\nu_{P_n}\sta\sigma, \quad n\to\infty,\ n\in\mathcal{N}_z'.
\label{eq:starsigma}
\end{equation}
Then, by the principle of descent (see \cite[Thm~I.6.8]{ST}) and (\ref{eq:liminflim}),
\begin{equation}\label{Ubeta}
U^\beta (z)=\liminf_{\underset{n\in\mathcal{N}_z'}{n\to\infty}} U^{\nu_{P_n}}(z)\geq
U^\sigma(z).
\end{equation}
Since $z\in\mathbb{C}$ was
arbitrary,
$$
U^\beta \geq \inf_{\sigma\in\mathcal{C}} U^\sigma \quad {\rm in\,\,}\mathbb{C},
$$
by which $U^\beta\geq U$ follows in all $\mathbb{C}$.

To finish the proof of (v), we let again $\mathcal{D}$ be a component of $\mathbb{C}\setminus{\rm supp\,}\beta$. By choosing above $z\in\mathcal{D}$ we get a measure $\sigma=\sigma_z\in\mathcal{C}$ with $U^\sigma (z)=U^\beta (z)$ (since equality necessarily holds in (\ref{Ubeta})). Thus $U^\sigma=U^\beta$ in all $\mathcal{D}$ because, as we have already proved, the other alternative would be $U^\sigma>U^\beta$ in all $\mathcal{D}$.

Regarding (vi), if $\mathcal{C}$ consists of only one point, say $\sigma$, then $U^\beta=U^\sigma$ by
(v), and from the unicity theorem for logarithmic potentials (see \cite[Thm~II.2.1]{ST}) we must have $\beta=\sigma$. Clearly, the full sequence must converge to $\beta$, because otherwise one could extract a subsequence converging to something else, which would be a different element in $\mathcal{C}$.

The assertions in (vii) are easy consequences of (iv) and (v): Since, for any $\sigma\in\mathcal{C}$, $U^\sigma=U^\beta$ in the unbounded component of $\mathbb{C}\setminus{\rm supp\,}\beta$ we get in the case of (a) plus (b) that (for any $\sigma\in\mathcal{C}$) $U^\sigma=U^\beta$, almost everywhere with respect to the area measure in $\mathbb{C}$. This and the unicity theorem yield $\beta=\sigma\in\mathcal{C}$. In the case of (a) plus (c), there exists (by (v)) at least one $\sigma\in\mathcal{C}$ satisfying $U^\sigma=U^\beta$ in the bounded component of $\mathbb{C}\setminus{\rm supp\,}\beta$, and for this $\sigma$ we have the same conclusion:
$U^\sigma=U^\beta$ almost everywhere in $\mathbb{C}$ and, as above, $\beta=\sigma\in\mathcal{C}$.

So far we have assumed that $N\geq 2$. Let us indicate the modifications needed for $N=1$. Equation
(\ref{eq:logrho}) may be written
\begin{equation}\label{eq:logrholim}
\log \rho(K(\cdot,\zeta))=\lim_{M\to+\infty}\min\left\{M,\ \log\rho_m(\varphi_j),\ \inf\{g_\Omega
(\varphi_j|_{\mathcal{G}_{j,\rho_m(\varphi_j)}}^{-1}({1}/{\overline{\varphi_j (\zeta)}}),\infty)\}\right\},
\end{equation}
that is,  by introducing an auxiliary upper bound $M$, which finally tends to infinity. Before passing to the limit we can work with the corresponding quantities
$$
h_M =\sup\{h, - M\}, \quad \beta_M = \frac{1}{2\pi}\Delta h_M
$$
(etc.) as before. Since a decreasing sequence of  subharmonic functions is subharmonic, $\displaystyle{h=\lim_{M\to \infty} h_M}$ will be again subharmonic. It is however not clear that it will be continuous, only upper semicontinuity is automatic. If $\rho_m(\varphi_j)<\infty$, then the bound $M$ is not needed and everything will be as in the case $N\geq 2$. So assume
$\rho_m(\varphi_j)=\infty$. This means that $\varphi_j$ is meromorphic in the entire complex plane and hence (\ref{eq:logrho}) reads
\begin{equation}\label{eq:logrhoC}
\log \rho(K(\cdot,\zeta))= \inf\{ g_\Omega
(\varphi_j|_{\mathbb{C}}^{-1} ({1}/{\overline{\varphi_j (\zeta)}}),\infty)\},\ \zeta\in G_j.
\end{equation}

Problems concerning the lower boundedness and continuity of $h$ could conceivably occur at points $\zeta\in G$ at which the inverse image above is either empty or is an infinite set. The first case
can, by Picard's theorem, occur for at most two values of $\zeta\in G$. At such points the infimum in (\ref{eq:logrhoC}) is $+\infty$, and hence $h(\zeta)=-\infty$. In particular, $h$ will not be bounded
from below, but it will still be subharmonic and upper semicontinuous. Moreover, it will be continuous at all other points, which is enough for the reasoning in the proof (above) of (iv), where we used the continuity of $h$ (or $U^\beta$).

The second conceivable problem, that $\varphi_1|_{\mathbb{C}}^{-1}({1}/{\overline{\varphi_1 (\zeta)}})$ is an infinite set, presents no actual difficulty because the only cluster points can be at infinity, hence all but finitely many branches of $\varphi_1|_{\mathbb{C}}^{-1} ({1}/{\overline{\varphi_1 (\zeta)}})$
will be ruled out when taking the infimum in (\ref{eq:logrhoC}).
\qed

\medskip
\noindent
{\bf Proof of Corollary \ref{cor:SS}.} As already remarked, the boundary curve $\Gamma_j$ is singular if and only if $\rho_m(\varphi_j)=1$, which by the proof of the theorem (e.g., Equation (\ref{eq:logrho}))
occurs if and only if $h=0$ in $G_j$. This, in view of (\ref{eq:potbeta}), is equivalent to
\begin{equation*}
U^\beta(z)=\log\frac{1}{{\rm cap\,}(\Gamma)},\ z\in G_j.
\end{equation*}
Also from (\ref{eq:potbeta}),
\begin{equation*}
U^\beta(z)=\log\frac{1}{{\rm cap\,}(\Gamma)}-g_{\Omega}(z,\infty),\ z\in\mathcal{G}_{j,R_j}\setminus G_j.
\end{equation*}
It follows that $U^\beta$ is harmonic in $\mathcal{G}_{j,R_j}\setminus\Gamma_j$, thus   ${\rm supp\,}\beta\subset\Gamma_j$. It also follows that the logarithmic potentials of $\beta$ and $\mu_\Gamma$ coincide in the domain $\mathcal{G}_{j,R_j}$, hence the equation
$\beta|_{\overline{G}_j}=\mu_\Gamma|_{\overline{G}_j}$ holds as a result of the unicity theorem (see e.g. \cite[p.~97]{ST}). This proves the equivalence of (i) and (ii).

By assertion (v) of the theorem, there exists a $\sigma\in\mathcal{C}$ such that $U^\sigma=U^\beta$ in $G_j\ (=\mathcal{D})$. The equation persists on $\Gamma_j$, because of the continuity of logarithmic potentials in the fine topology and in view of (\ref{eq:usigma}), it also holds in any neighborhood of $\overline{G}_j$ not meeting the other islands. Thus, from the unicity theorem $\sigma=\beta$, in such a neighborhood. As $\sigma$ is a cluster point of $\{\nu_{P_n}\}$, we conclude that (iii) follows from (ii).

If (iii) holds then by selecting a further subsequence we conclude  $\sigma|_V=\mu_\Gamma|_V$, for some $\sigma\in\mathcal{C}$. Then $U^\sigma=U^{\mu_\Gamma}$ in $V$, which in view of (\ref{eq:bal}) and (\ref{eq:usigma}) yields the relation $U^\beta=U^{\mu_\Gamma}$ in $V$. Therefore $\beta|_{\overline{G}_j}=\mu_\Gamma|_{\overline{G}_j}$. \qed

\medskip
\noindent
\textbf{Proof of Corollary \ref{cor:bala}.} Set $\mu_n={\rm Bal\,}(\nu_{P_n})$. Then
\begin{equation}\label{eq:supplambda}
{\rm supp\,}\mu_n \subset\mathbb{C}\setminus G,
\end{equation}
\begin{equation}\label{eq:Ulambda}
U^{\nu_{P_n}}= U^{\mu_n}\quad {\rm in}\quad \Omega.
\end{equation}
Let $\mu$ be any weak* cluster point of $\{\mu_n\}$ and let $\mathcal{N}\subset\mathbb{N}$ be a subsequence with $\mu_n\sta\mu$, $n\in\mathcal{N}$. By refining $\mathcal{N}$ we may assume also that $\nu_{P_n}\sta\sigma$, $n\in\mathcal{N}$, for some measure $\sigma$. Then in view of (\ref{eq:Ulambda}) we have $U^\sigma=U^{\mu}$ in $\Omega$.

On the other hand, $U^\sigma=U^{\mu_\Gamma}$ in $\Omega$ by Theorem~\ref{thm:ASgeneral}, thus $U^\mu=U^{\mu_\Gamma}$ in $\Omega$. But $U^{\mu_\Gamma}$ is harmonic in $\Omega\setminus\{\infty\}$ and ${\rm supp\,}\mu \subset\mathbb{C}\setminus G$ by (\ref{eq:supplambda}), hence ${\rm supp\,}\mu\subset\Gamma$. Now Carleson's unicity theorem \cite[p.~123]{ST}, shows that $\mu=\mu_\Gamma$. Since $\mu$ was an arbitrary cluster point of $\mu_n$ it follows that
$\mu_n\sta\mu_\Gamma$ for the full sequence. \qed

\medskip
\noindent
\textbf{Proof of Corollary \ref{cor:AS}.}
The expression for $U^\beta$ follows immediately after uploading (\ref{eq:rhoCaseII}) into Theorem~\ref{thm:ASgeneral}~(ii). From this expression and the unicity theorem for logarithmic potentials we gather that ${\rm supp\,}\beta$ must be contained in $\partial E$. To show that eventually  ${\rm supp\,}\beta=\partial E$ we can argue as in \cite[pp.~215--216]{M-DSS}. That is, by assuming that a point $z_0\in\partial E$ does not belong to ${\rm supp\,}\beta$, hence the potential $U^\beta$ is harmonic in a small disk centered at $z_0$, we arrive to a contradiction by comparing the resulting harmonic extension of $U^\beta$ with the one  given in (\ref{eq:betainAS}).

In view of the connectedness of the complement of $E$ and the fact that the support of $\beta$ is contained in $\overline{E}$ the equality $U^\sigma(z)=U^\beta(z)$, for $z\in\overline{\mathbb{C}}\setminus\overline{E}$, is immediate from Theorem~\ref{thm:ASgeneral}~(iv). Hence ${\rm supp\,}\sigma\subset\overline{E}$. Furthermore, since the boundary of the domain $\overline{\mathbb{C}}\setminus\overline{E}$ in the fine topology coincides with its boundary in the Euclidean topology (see e.g. \cite[Cor.~I.5.6]{ST}), we conclude that the equality between the potentials persists in $\overline{\mathbb{C}}\setminus{E}$.

The last assertion in the corollary can be deduced from Theorem~\ref{thm:ASgeneral}~(iv)--(v), because this guarantees the existence of a cluster point $\sigma$ of the sequence $\nu_{P_n}$ such that $U^\sigma = U^\beta$ on both sides of $\Gamma_1$. More precisely, $U^\sigma = U^\beta$ in $V\setminus\Gamma_1$, where $V$ is a neighborhood of $\overline{G}_1$ not meeting the other islands, and therefore $\sigma=\beta$ in such a neighborhood.  Similarly we argue for $L_{2,\frac{1}{R'}}$.
\qed

\subsection{The lemniscate example} $\, $

\medskip
\noindent
\textbf{Proof of Lemma~\ref{lem:lemn}.} Let $(\gamma w+1)^{\tau/2}$ denote
the analytic branch in $\mathbb{D}=\{w: |w|<1\}$ that equals $1$
at $w=0$. Then applying Green's formula we have, for $j=0,1,
\ldots,n-1$,
\begin{eqnarray*}
0&=&\int\limits_{\mathbb{D}}\beta_n(w)\overline{(\gamma
w+1)}^j\frac{\mathit{d A(w)}}{|\gamma
w+1|^\tau}=\int\limits_{\mathbb{D}}\frac{\beta_n(w)}{(\gamma
w+1)^{\tau/2}}\overline{(\gamma w+1)}^{j-\tau/2}\mathit{d A(w)}\\
&=&\int\limits_{|w|=1}\frac{\beta_n(w)}{|\gamma
w+1|^\tau}\overline{(\gamma w+1)}^{j+1}w|\mathit{d
w}|=\int\limits_{|w|=1}\beta_n(w)(\bar{\gamma
}+w)\frac{\overline{(\gamma w+1)}^j}{|\gamma w+1|^\tau}|\mathit {d
w}|,
\end{eqnarray*}
where we have ignored nonzero constants, and in the last equality, we used that $\overline{(\gamma w+1)}=(\bar{\gamma}/w+1)$ for $|w|=1$. Consequently, $\beta_n(w)(\bar{\gamma}+w)$ is a monic
polynomial of degree $n+1$ that vanishes at $w=-\bar{\gamma}$ and is orthogonal to all polynomials of degree less than $n$ with respect to $|\mathit{d w}|/|\gamma w+1|^\tau$. The same is true of the right-hand side of (\ref{eq:lemlemn}) and hence the difference of these two polynomials (which is of degree $\leq n$) must be a multiple of $t_n(w)$ that vanishes at $-\bar{\gamma}$. Since
$t_n(-\bar{\gamma})\neq 0$, the difference of the left and right-hand sides of (\ref{eq:lemlemn}) must be identically zero.\qed

\begin{remark}
It is essential that the cases $\tau=2,4,\ldots,2n$ be excluded in Lemma~\ref{lem:lemn}. Indeed for $\tau=2j$, where $j$ is a positive integer, it is well-known (cf. \cite{Sz}, \S 11.2) that $t_n(w)=w^{n-j}(w+\bar{\gamma})^j$ for $n\geq j$, so that $t_n(-\bar{\gamma})=0$ in this case. There appears, however, to be no simple formula\footnote[2]{For the weight $\mathit{d A}/|\gamma w+1|^2$, we have $$\beta_1(w)=w+\frac{1}{\gamma}+\frac{\bar{\gamma}}{\ln(1-|\gamma|^2)}.$$}
for the polynomials $\beta_n(w)$ for such values of $\tau$. We shall show in Lemma \ref{lemmpi} that if $\tau$ is not an even integer, then $t_n(-\bar{\gamma})\neq 0$ for all $n$ sufficiently large.
\end{remark}

\medskip
\noindent
\textbf{Proof of Proposition~\ref{pro:lemn2} .} Here we use the minimality property of the monic Bergman polynomials $p_{km+s}(z)=z^sq_{k,s}(z^m)$. More precisely, $q_{k,s}$ solves the extremal
problem
\begin{eqnarray}\label{e1}
I_{k,s}:=\min\{\int\limits_G|z^sq(z^m)|^2\mathit{d A}:
q(t)=t^k+\cdots \in \mathcal{P}_k\}.
\end{eqnarray}
Clearly,$$\int\limits_G|z^sq(z^m)|^2\mathit{dA}=m\int\limits_{G_m}|z^sq(z^m)|^2\mathit{d A},$$ and the change of variables $w=(z^m-1)/r^m$, which maps $G_m$ conformally onto the unit disk $\mathbb{D}$ in the $w$-plane, yields
$$
\int\limits_{G_m}|z^sq(z^m)|^2\mathit{d A(z)}=\frac{r^{2m}}{m^2}
\int\limits_\mathbb{D}\frac{|q(r^mw+1)|^2}{|r^mw+1|^\tau}\mathit{d A(w)},
$$
where
\begin{eqnarray}\label{e2}
\tau:=2-\frac{2}{m}-\frac{2s}{m}.
\end{eqnarray}
Consequently,
\begin{eqnarray}\label{e3}
I_{k,s}=\frac{r^{2m}}{m}\min\{\int\limits_\mathbb{D}\frac{|q(r^mw+1)|^2}{|r^mw+1|^\tau}\mathit{d
A(w)}: q(t)=t^k+\cdots\in \mathcal{P}_k \},
\end{eqnarray}
and, moreover, $r^{-mk}q_{k,s}(r^mw+1)$ is the monic (in $w$) orthogonal polynomial with respect to the weight $\mathit{dA(w)}/|r^mw+1|^\tau$ on $\mathbb{D}$. Applying Lemma~\ref{lem:lemn} then
yields formulas (\ref{eq:proplemn2.a}) and (\ref{eq:proplemn2.b}), provided that $\pi_{k,s}(-r^m)$ is not zero. In the next lemma we show that this condition is indeed satisfied for $k$ sufficiently large.\qed

\begin{lemma}\label{lemmpi}
Let $\pi_{k,s}(w)$ be as in Proposition~\ref{pro:lemn2} and $\tau$ be given by (\ref{e2}). Then, for each $s=0,1,\ldots,m-2$, we have
\begin{eqnarray}\label{e4}
(-1)^k\frac{k^{\tau/2}}{r^{mk}}\pi_{k,s}(-r^m)=\sin(\tau\pi/2)
\Big[\frac{1}{\pi}\Gamma\Big(\frac{\tau}{2}\Big)+\frac{b_s}{k}+
\mathit{O}\Big(\frac{1}{k^2}\Big)\Big]
\end{eqnarray}
as $k\rightarrow\infty$, where $b_s$ is a constant independent of $k$.
\end{lemma}
\textbf{Proof.} As in \cite{MDthe}, we utilize the results of \cite{MMS} for Szeg\H{o} polynomials with respect to an analytic weight on $|w|=1$. For the weight $|w+r^m|^{-\tau}=1/|r^mw+1|^\tau$, we have, imitating the notation of \cite{MMS}, the following formulas for the exterior and interior Szeg\H{o} functions $D_{e,\tau}(w)$ and $D_{i,\tau}(w)$, respectively,
\begin{equation}\label{e5}
D_{e,\tau}(w)=\Big(\frac{w+r^m}{w}\Big)^{\tau/2}, \quad
D_{i,\tau}(w)=(1+r^mw)^{-\tau/2},
\end{equation}
where the branches of the square roots are chosen so that $D_{e,\tau}(\infty)=D_{i,\tau}(0)=1$. The scattering function $S_\tau(w)$ is given by
\begin{eqnarray}\label{e6}
S_\tau(w)=D_{e,\tau}(w)D_{i,\tau}(w)=\Big(\frac{w+r^m}{w}\Big)^{\tau/2}(1+r^mw)^{-\tau/2}
\quad \text{for} \quad r^m<|w|<r^{-m}.
\end{eqnarray}

As shown in \cite{MMS} (see Equations (16), (25), and (39)), we have for
$|w|<\eta$, where $r^m<\eta<1$,
\begin{eqnarray}\label{e7}
D_{i,\tau}(w)\pi_{k,s}(w)=\frac{1}{2\pi
i}\oint\limits_{|t|=1}\frac{t^kS_\tau(t)}{t-w}\mathit{d
t}+\mathit{O}(\eta^{3k}), \quad \text{as} \quad k\rightarrow
\infty.
\end{eqnarray}
For $w=-r^m$, we can deform the unit circle in the integral in
(\ref{e7}) so that the integration takes place along each side of
the branch cut of $D_{e,\tau}(w)$ joining $-r^m$ to $0$ to obtain
\begin{eqnarray}\label{e8}
I_k:=\oint\limits_{|t|=1}\frac{t^kS_\tau(t)}{t+r^m}\mathit{dt}=
\left(\int\limits_{[-r^m,0]}+\int\limits_{[0,-r^m]}\right)\frac{x^kS_\tau(x)}{x+r^m}\mathit{dx},
\end{eqnarray}
where we utilize the limiting values from below for $S_\tau$ in
integrating from $-r^m$ to $0$ and the limiting values of $S_\tau$
from above in integrating from $0$ to $-r^m$. Thus we get (cf.~(\ref{e6}))
$$
I_k=2i\sin(\tau\pi/2)\int\limits_{-r^m}^0
\frac{x^k(1+r^mx)^{-\tau/2}}{|x|^{\tau/2}(x+r^m)^{1-\tau/2}}\mathit{d x},
$$
and on making the change of variable $x=-r^m(1+\cos\theta)/2$ we
find that
\begin{eqnarray}\label{e9}
I_k=\frac{ir^{mk}}{2^{k-1}}\sin(\tau\pi/2)(-1)^k\int\limits_{0}^\pi
e^{-kp(\theta)}q(\theta)\mathit{d \theta},
\end{eqnarray}
where $p(\theta):=-\log(1+\cos\theta)$ and
\begin{eqnarray}\label{e10}
q(\theta):=\Big[1-\frac{r^{2m}}{2}(1+\cos\theta)\Big]^{-\tau/2}
(1+\cos\theta)^{1-\tau}\theta^{\tau-1}
\Big(\frac{\sin\theta}{\theta}\Big)^{\tau-1}.
\end{eqnarray}

We now apply Laplace's method  to deduce the
asymptotic behavior of the integral in (\ref{e9}). Since
$$p(\theta)=-\log 2 +\sum_{j=0}^\infty p_j\theta^{j+2}=-\log 2+\frac{1}{4}\theta^2+\cdots$$
and $$q(\theta)=\sum_{j=0}^\infty
q_j\theta^{j+\tau-1}=(1-r^{2m})^{-\tau/2}2^{1-\tau}\theta^{\tau-1}+q_2\theta^{\tau+1}+\cdots,$$
(note that $q_1=0$) we obtain from \cite[Ch. 3, Thm 8.1]{Ol97}, that, as
$k\rightarrow \infty$,
\begin{eqnarray}\label{e11}
\int\limits_0^\pi e^{-kp(\theta)}q(\theta)\mathit{d
\theta}=2^k\Big[\Gamma\Big(\frac{\tau}{2}\Big)\frac{(1-r^{2m})^{-\tau/2}}{k^{\tau/2}}+
\frac{a_{2,\tau}}{k^{\tau/2+1}}+\mathit{O}\Big(\frac{1}{k^{\tau/2+2}}\Big)\Big],
\end{eqnarray}
where $a_{2,\tau}$ is a constant independent of $k$.
From (\ref{e7})--\ref{e11})
(taking $\eta$ such that $\eta^3<r^m<\eta$) we deduce (\ref{e4}).\qed

As an immediate consequence of the preceding lemma we obtain that
\begin{eqnarray}\label{e12}
\frac{\pi_{k+1,s}(-r^m)}{\pi_{k,s}(-r^m)}=-
r^m\Big[1-\frac{\tau}{2k}+\mathit{O}\Big(\frac{1}{k^2}\Big)\Big]\quad
\text{as} \quad k\rightarrow \infty.
\end{eqnarray}

\medskip
\noindent
\textbf{Proof of Proposition~\ref{pro:lemn3} } For $s=m-1$ the assertion is
obvious from (\ref{eq:proplemn2.a}). For $|z^m-1|>r^{2m}$ and $s=0,1,\ldots,m-2,$ we appeal to the
well-known fact regarding exterior asymptotics of Szeg\H{o}
polynomials (see e.g. \cite{MMS}, Proposition 1) that for $|w|>r^m$ we
have
\begin{eqnarray}\label{e13}
\lim_{k\rightarrow
\infty}\frac{\pi_{k,s}(w)}{w^k}=D_{e,\tau}(w)=\Big(\frac{w+r^m}{w}\Big)^{\tau/2},
\end{eqnarray}
where the convergence is locally uniform and takes place with a
geometric rate. Thus from (\ref{e12}) and the representation (\ref{eq:proplemn2.b})
we deduce (\ref{eq:proplemn3.b}) for $|z^m-1|>r^{2m}$.

For $|z^m-1|\leq r^{2m}$, we begin with the asymptotic analysis of
$\pi_{k,s}(w)$, for $s=0,1,\ldots,m-2$ and $|w|\leq r^m$. Assume
at first that $w\notin [-r^m,0]$, and consider the integral in the
representation (\ref{e7}). For each $\epsilon >0$ sufficiently
small, we can write
\begin{eqnarray}\label{e14}
J_k(w):=\frac{1}{2\pi
i}\oint\limits_{|t|=1}\frac{t^kS_\tau(t)}{t-w}dt=
\frac{1}{2\pi i}\left(\oint\limits_{|t-w|=\epsilon}+\int\limits_{[-r^m,0]}+\int\limits_{[0,-r^m]}\right)
\frac{t^kS_\tau(t)}{t-w}\mathit{d
t},
\end{eqnarray}
where integration along both sides of the branch cut from $-r^m$
to $0$ is as in the proof of Lemma \ref{lemmpi}. From Cauchy's formula and
the representation of $S_\tau(t)$ along each side of the branch
cut, we deduce that
$$J_k(w)=w^kS_\tau(w)+\frac{1}{\pi}\sin(\tau\pi/2)\int\limits_{-r^m}^0
\frac{x^k(1+r^mx)^{-\tau/2}(x+r^m)^{\tau/2}}{|x|^{\tau/2}(x-w)}\mathit{d x},$$
which, upon performing the change of variable
$x=-r^m(1+\cos\theta)/2$, yields
\begin{eqnarray}\label{e15}
J_k(w)=w^kS_\tau(w)+\frac{1}{\pi}\sin(\tau\pi/2)(-1)^{k+1}\frac{r^{m(k+1)}}{2^{k+1}}\int\limits_0^\pi
e^{-kp(\theta)}\hat{q}(\theta)\mathit{d \theta},
\end{eqnarray}
where $p(\theta)=-\log(1+\cos\theta)$ and
$$
\hat{q}(\theta):=\frac{[1-\frac{r^{2m}}{2}(1+\cos\theta)]^{-\tau/2}
(\frac{\sin\theta}{\theta})^{\tau+1}\theta^{\tau+1}}{[\frac{r^m}{2}(1+\cos\theta)+w]
(1+\cos\theta)^\tau}.
$$
Since
$$\hat{q}(\theta)=\sum_{j=0}^\infty\hat{q}_j\theta^{j+(\tau+2)-1}=
\frac{(1-r^{2m})^{-\tau/2}}{(r^m+w)2^\tau}\theta^{\tau+1}+\hat{q}_3\theta^{\tau+3}+\cdots$$
(note that $\hat{q}_1=0$), Laplace's method yields
$$
\int\limits_0^\pi e^{-kp(\theta)}\hat{q}(\theta)\mathit{d \theta}=2^k\Big[\frac{\tau
\Gamma(\frac{\tau}{2})(1-r^{2m})^{-\tau/2}}{r^m+w}\frac{1}{k^{1+\tau/2}}+
\frac{\hat{b}_s(w)}{k^{2+\tau/2}}+\mathit{O}\Big(\frac{1}{k^{3+\tau/2}}\Big)\Big],
$$
as $k\rightarrow \infty$, where $\hat{b}_s(w)$ is a constant
independent of $k$. Thus, from (\ref{e15}) and (\ref{e7}), we obtain
\begin{eqnarray}\label{e16}
D_{i,\tau}(w)\pi_{k,s}(w)\frac{k^{1+\tau/2}(-1)^{k+1}}{r^{m(k+1)}}=
\frac{\sin(\tau\pi/2)\tau\Gamma(\tau/2)}{2\pi(1-r^{2m})^{\tau/2}(w+r^m)}
\Big[1+\frac{\hat{b}_s(w)}{k}+\mathit{O}\Big(\frac{1}{k^2}\Big)\Big],
\end{eqnarray}
as $k\rightarrow \infty$, provided $|w|<r^m$ and
$\eta^3<r^m<\eta$, while for $|w|=r^m$, $w\neq -r^m$, we obtain
\begin{eqnarray}\label{e17}
D_{i,\tau}(w)\frac{\pi_{k,s}(w)}{w^k}=S_\tau(w)+\mathit{O}\Big(\frac{1}{k^{1+\tau/2}}\Big),
\end{eqnarray}
as $k\to\infty$, where we take $r^m<\eta<1$.

Combining (\ref{e12}) with (\ref{e16}) and (\ref{e17}), we deduce
from the representation (\ref{eq:proplemn2.b}) that (\ref{eq:proplemn4.a}) holds for
$|z^m-1|<r^{2m}$, $z^m\notin[1-r^{2m},1]$, and that (\ref{eq:proplemn3.b}) holds
for $|z^m-1|=r^{2m}$, except for the $m$ roots $(1-r^{2m})^{1/m}$. In
deriving (\ref{eq:proplemn4.a}) we used the fact that
$(z^m)^{\tau/2}z^s=z^{m-1}e^{2\pi ij(s+1)/m}$ for $z\in G_j$
(recall (\ref{eq:Gjcon})). Finally, by a slight modification of the above
analysis, it is easy to see that (\ref{e16}) is valid also for
$w\in (-r^m,0]$ and so (\ref{eq:proplemn4.a}) holds for all $z$ satisfying
$|z^m-1|<r^{2m}$. \qed

\medskip
\noindent
\textbf{Proof of Proposition~\ref{pro:lemn1}.} We use the obvious fact that
\begin{eqnarray}\label{e18}
\lambda^{-2}_{km+s}=\int\limits_G|p_{km+s}(z)|^2\mathit{d A}(z).
\end{eqnarray}
For $s=m-1$, we have from (\ref{eq:proplemn2.a}),
\begin{eqnarray*}
\lambda^{-2}_{km+m-1}&=&\int\limits_G|z^{m-1}(z^m-1)^k|^2d A(z)=
m\int\limits_{G_m}|z^{m-1}(z^m-1)^k|^2\mathit{d A(z)}\\
&=&\frac{r^{2m(k+1)}}{m}\int\limits_\mathbb{D}|w|^{2k}d A(w)=\frac{\pi r^{2m(k+1)}}{m(k+1)},
\end{eqnarray*}
where, as in the proof of Lemma \ref{lemmpi}, we have made the change of
variables $w=(z^m-1)/r^m$. Thus
\begin{eqnarray}\label{e19}
\lambda_{km+m-1}=\sqrt{\frac{m(k+1)}{\pi r^{2m(k+1)}}}.
\end{eqnarray}

Now suppose that $0\leq s<m-1$. Then, on utilizing the formula (\ref{eq:proplemn2.b}) we
deduce that, for $k$ sufficiently large,
\begin{eqnarray}\label{e20}
\lambda^{-2}_{km+s}&=&m\int\limits_{G_m}|z^sq_{k,s}(z^m)|^2\mathit{d A(z)}
=\frac{r^{2m}}{m}\int\limits_\mathbb{D}
\frac{|q_{k,s}(r^mw+1)|^2}{|r^mw+1|^\tau}\,\mathit{d A(w)}\nonumber\\
&=&\frac{r^{2m(k+1)}}{m}\int\limits_\mathbb{D}
\frac{|\pi_{k+1}(w)-\frac{\pi_{k+1}(-r^m)}{\pi_{k}(-r^m)}\pi_k(w)|^2}
{|w+r^m|^2\,|r^mw+1|^\tau}\,\mathit{d A(w)},
\end{eqnarray}
where for simplicity of notation we have written
$\pi_k=\pi_{k,s}$. On using the orthogonality property of the
$\pi_k$'s we can simplify the last integral in (\ref{e20}) to
obtain
\begin{eqnarray}\label{e21}
\lambda^{-2}_{km+s}=\frac{-\pi_{k+1}(-r^m)r^{2mk+m}}{\pi_k(-r^m)2m(k-\frac{\tau}{2}+1)}
\int\limits_{|w|=1}\frac{|\pi_k(w)|^2}{|r^mw+1|^\tau}|\mathit{d
w}|.
\end{eqnarray}

Finally we note that the integral on the right-hand side of (\ref{e21}) equals
$\mu^{-2}_{k,s}$, where $\mu_{k,s}$ is the leading coefficient of
the orthonormal polynomial with respect to the weight $|\mathit{d
w}|/|r^mw+1|^\tau$ on the unit circle. As is well-known (see e.g. \cite{MMS}, Corollary 2)
$$|\mu^2_{k,s}-\frac{1}{2\pi}|=\mathit{O}(\eta^{2k}) \quad \text{as}\,\, k\rightarrow\infty,$$
where $r^m<\eta<1$. Combining this fact with (\ref{e21}) and (\ref{e12}) yields
(\ref{eq:proplemn1}).\qed

\bibliographystyle{amsplain}

\providecommand{\bysame}{\leavevmode\hbox to3em{\hrulefill}\thinspace}
\providecommand{\MR}{\relax\ifhmode\unskip\space\fi MR }
\providecommand{\MRhref}[2]{%
  \href{http://www.ams.org/mathscinet-getitem?mr=#1}{#2}
}
\providecommand{\href}[2]{#2}
\begin{thebibliography}{}

\end{thebibliography}


\begin{thebibliography}{10}

\bibitem{Am95}
A.~Ambroladze, \emph{On exceptional sets of asymptotic relations for general
  orthogonal polynomials}, J. Approx. Theory \textbf{82} (1995), no.~2,
  257--273.

\bibitem{AB}
V.~V. Andrievskii and H.-P. Blatt, \emph{Erd{\H o}s-{T}ur\'an type theorems on
  quasiconformal curves and arcs}, J. Approx. Theory \textbf{97} (1999), no.~2,
  334--365.

\bibitem{Ca23}
T.~Carleman, \emph{\"{U}ber die {A}pproximation analytisher {F}unktionen durch
  lineare {A}ggregate von vorgegebenen {P}otenzen}, Ark. Mat., Astr. Fys.
  \textbf{17} (1923), no.~9, 215--244.

\bibitem{Da74}
P.~J. Davis, \emph{The {S}chwarz function and its applications}, The
  Mathematical Association of America, Buffalo, N. Y., 1974, The Carus
  Mathematical Monographs, No. 17.

\bibitem{FCB}
M.~D. Finn, S.~M. Cox, and H.~M. Byrne, \emph{Topological chaos in inviscid and
  viscous mixers}, J. Fluid Mech. \textbf{493} (2003), 345--361.

\bibitem{Ga}
D.~Gaier, \emph{Lectures on complex approximation}, Birkh\"auser Boston Inc.,
  Boston, MA, 1987, Translated from the German by Renate McLaughlin.

\bibitem{GGMPV}
G.~Golub, B.~Gustafsson, C.~He, P.~Milanfar, M.~Putinar, and J.~Varah,
  \emph{Shape reconstruction from moments: theory, algorithms, and
  applications}, SPIE Proccedins (F.~T. Luk, ed.), Advanced Signal Processing,
  Algorithms, Architecture, and Implementations X, vol. 4116, 2000,
  pp.~406--416.

\bibitem{GHMP}
B.~Gustafsson, C.~He, P.~Milanfar, and M.~Putinar, \emph{Reconstructing planar
  domains from their moments}, Inverse Problems \textbf{16} (2000), no.~4,
  1053--1070.

\bibitem{Ha89}
W.~K. Hayman, \emph{Subharmonic functions. {V}ol. 2}, London Mathematical
  Society Monographs, vol.~20, Academic Press, London, 1989.

\bibitem{HKZ}
H.~Hedenmalm, B.~Korenblum, and K.~Zhu, \emph{Theory of {B}ergman spaces},
  Graduate Texts in Mathematics, vol. 199, Springer-Verlag, New York, 2000.

\bibitem{LSS}
A.~L. Levin, E.~B. Saff, and N.~S. Stylianopoulos, \emph{Zero distribution of
  {B}ergman orthogonal polynomials for certain planar domains}, Constr. Approx.
  \textbf{19} (2003), no.~3, 411--435.

\bibitem{MMS}
A.~Mart{\'{\i}}nez-Finkelshtein, K.~T.-R. McLaughlin, and E.~B. Saff,
  \emph{Szeg{\H o} orthogonal polynomials with respect to an analytic weight:
  canonical representation and strong asymptotics}, Constr. Approx. \textbf{24}
  (2006), no.~3, 319--363.

\bibitem{MS}
V.~Maymeskul and E.~B. Saff, \emph{Zeros of polynomials orthogonal over regular
  {$N$}-gons}, J. Approx. Theory \textbf{122} (2003), no.~1, 129--140.

\bibitem{MDthe}
E.~Mi{\~n}a-D{\'{\i}}az, \emph{Asymptotics for {F}aber polynomials and
  polynomials orthogonal over regions in the complex plane}, Ph.D. thesis,
  Vanderbilt {U}niversity, August 2006.

\bibitem{MD08}
\bysame, \emph{An asymptotic integral representation for {C}arleman orthogonal
  polynomials}, Int Math Res Notices \textbf{2008} (2008), article ID rnn066,
  35 pages.

\bibitem{M-DSS}
E.~Mi{\~n}a-D{\'{\i}}az, E.~B. Saff, and N.~S. Stylianopoulos, \emph{Zero
  distributions for polynomials orthogonal with weights over certain planar
  regions}, Comput. Methods Funct. Theory \textbf{5} (2005), no.~1, 185--221.

\bibitem{Ol97}
F.~W.~J. Olver, \emph{Asymptotics and special functions}, AKP Classics, A K
  Peters Ltd., Wellesley, MA, 1997, Reprint of the 1974 original (Academic
  Press, New York).

\bibitem{PSG}
N.~Papamichael, E.~B. Saff, and J.~Gong, \emph{Asymptotic behaviour of zeros of
  {B}ieberbach polynomials}, J. Comput. Appl. Math. \textbf{34} (1991), no.~3,
  325--342.

\bibitem{PaWa}
N.~Papamichael and M.~K. Warby, \emph{Stability and convergence properties of
  {B}ergman kernel methods for numerical conformal mapping}, Numer. Math.
  \textbf{48} (1986), no.~6, 639--669.

\bibitem{Ra}
T.~Ransford, \emph{Potential theory in the complex plane}, London Mathematical
  Society Student Texts, vol.~28, Cambridge University Press, Cambridge, 1995.

\bibitem{Sa69}
E.~B. Saff, \emph{Polynomials of interpolation and approximation to meromorphic
  functions}, Trans. Amer. Math. Soc. \textbf{143} (1969), 509--522.

\bibitem{Sa90}
\bysame, \emph{Orthogonal polynomials from a complex perspective}, Orthogonal
  polynomials (Columbus, OH, 1989), Kluwer Acad. Publ., Dordrecht, 1990,
  pp.~363--393.

\bibitem{SaSt08}
E.~B. Saff and N.~S. Stylianopoulos, \emph{Asymptotics for polynomial zeros:
  Beware of predictions from plots}, Comput. Methods Funct. Theory \textbf{8}
  (2008), no.~2, 185--221.

\bibitem{ST}
E.~B. Saff and V.~Totik, \emph{Logarithmic potentials with external fields},
  Springer-Verlag, Berlin, 1997.

\bibitem{Sh92}
H.~S. Shapiro, \emph{The {S}chwarz function and its generalization to higher
  dimensions}, University of Arkansas Lecture Notes in the Mathematical
  Sciences, 9, John Wiley \& Sons, New York, 1992.

\bibitem{StTo}
H.~Stahl and V.~Totik, \emph{$n$th root asymptotic behavior of orthonormal
  polynomials}, Orthogonal polynomials (Columbus, OH, 1989), Kluwer Acad.
  Publ., Dordrecht, 1990, pp.~395--417.

\bibitem{StTobo}
\bysame, \emph{General orthogonal polynomials}, Cambridge University Press,
  Cambridge, 1992.

\bibitem{NSpre}
N.~S. Stylianopoulos, \emph{The use of orthogonal {B}ergman polynomials for
  recovering planar domains from their moments}, preprint.

\bibitem{Su66}
P.~K. Suetin, \emph{Order comparison of various norms of polynomials in a
  complex region}, Ural. Gos. Univ. Mat. Zap. \textbf{5} (1966), no.~tetrad 4,
  91--100, (in Russian).

\bibitem{Su74}
\bysame, \emph{Polynomials orthogonal over a region and {B}ieberbach
  polynomials}, American Mathematical Society, Providence, R.I., 1974,
  Translated from the Russian by R. P. Boas.

\bibitem{Sz}
G.~Szeg{\H{o}}, \emph{Orthogonal polynomials}, fourth ed., Colloquium
  Publications, Vol. XXIII, American Mathematical Society, Providence, R.I.,
  1975.

\bibitem{Topre}
V.~Totik, \emph{Christoffel functions on curves and domains}, preprint.

\bibitem{To08}
\bysame, \emph{Orthogonal polynomials}, Surv. Approx. Theory \textbf{1} (2005),
  70--125.

\bibitem{Tre}
L.~N. Trefethen, \emph{{Ten-digits algorithms}}, Report 05/13, Oxford
  University Computing Laboratory, 2005.

\bibitem{Wa63}
J.~L. Walsh, \emph{A sequence of rational functions with application to
  approximation by bounded analytic functions}, Duke Math. J. \textbf{30}
  (1963), 177--189.

\bibitem{Wa}
\bysame, \emph{Interpolation and approximation by rational functions in the
  complex domain}, fourth ed., Colloquium Publications, Vol. XX, American
  Mathematical Society, Providence, R.I., 1965.

\bibitem{Wi67}
H.~Widom, \emph{Polynomials associated with measures in the complex plane}, J.
  Math. Mech. \textbf{16} (1967), 997--1013.

\bibitem{Wi69}
\bysame, \emph{Extremal polynomials associated with a system of curves in the
  complex plane}, Advances in Math. \textbf{3} (1969), 127--232 (1969).

\end{thebibliography}
\def\cprime{$'$}
\providecommand{\bysame}{\leavevmode\hbox to3em{\hrulefill}\thinspace}
\providecommand{\MR}{\relax\ifhmode\unskip\space\fi MR }
\providecommand{\MRhref}[2]{%
  \href{http://www.ams.org/mathscinet-getitem?mr=#1}{#2}
}
\providecommand{\href}[2]{#2}

\end{document}